\documentclass[a4paper,
10pt
]{article}

\usepackage[english]{babel}

\usepackage[a4paper,left=3cm,right=3cm,top=2.5cm,bottom=2.5cm]{geometry}

\usepackage{amsfonts} 
\usepackage{hyperref}
\usepackage{bm}
\usepackage{authblk}
\usepackage{siunitx}
\usepackage{graphicx}
\usepackage[english]{babel}
\usepackage{amsmath}
\usepackage{xcolor}
\usepackage{subcaption}
\usepackage{float}
\usepackage{multirow}
{\color{blue}}

\usepackage{comment}
\usepackage{stackengine,graphicx}
\usepackage{tikz}
\usepackage[linesnumbered,ruled,vlined]{algorithm2e}
\usepackage{esint}
\usepackage{amsthm} % if not already included

\theoremstyle{remark}
\newtheorem{remark}{Remark}[section]

\title{Coupling Physics Informed Neural Networks with External Solvers}
\author[1]{Rahul Halder\footnote{rhalder@sissa.it}}
\author[2]{Giovanni Stabile\footnote{giovanni.stabile@santannapisa.it}}
\author[1]{Gianluigi Rozza\footnote{grozza@sissa.it}}

\affil[1]{Mathematics Area, mathLab, SISSA, via Bonomea 265, I-34136 Trieste, Italy}
\affil[2]{The Biorobotics Institute, Sant'Anna School of Advanced Studies, V.le R. Piaggio 34, 56025, Pontedera, Pisa - Italy}
\date{} % Leave empty to omit a date

\begin{document}
\maketitle
\newcommand{\GS}[1]{\noindent \textcolor{red}{GS: #1}}
%\maketitle
\begin{center}
\bf ABSTRACT
\end{center}

The current work aims to incorporate physics-based loss in Physics Informed Neural Network (PINN) directly using the numerical residual obtained from the governing equation in any dicretized forward solver. PINN's major difficulties for coupling with external forward solvers arise from the inability to access the discretized form (Finite difference, finite volume, finite element, etc.) of the governing equation directly through the network and to include them in its computational graph. This poses a significant challenge to conventional automatic-differentiation-based derivative computation of physics-based loss terms concerning the neural network hyperparameters if gradient-based optimization techniques are adopted. Therefore, we propose modifying the physics-based loss term to account for the residual arising from the external solver and to compute the derivative required for the optimization machinery. The proposed methodologies are demonstrated on benchmark full-order and reduced-order systems. 

\maketitle


\begin{center}
\textbf{Code available at:} \\
\url{https://github.com/rahulhalderAERO/PINNFOAM}
\end{center}

\section{\bf{Introduction}}\label{sec:Intro}

In recent years, the integration of Computational Fluid Dynamics (CFD) with Machine Learning (ML) techniques has attracted significant research attention. This integration spans various domains, including Reduced Order Model (ROM) methods, Physics-Informed Neural Networks (PINNs), and the enhancement of existing numerical techniques using ML-based strategies. The primary objective of combining ML techniques with CFD is to enable rapid exploration of the parametric space, which plays a crucial role in a wide range of applications, such as design optimization, control, and real-time prediction. We will examine a brief survey of different approaches and contextualise our current research direction within the general framework of scientific Machine Learning \cite{Quarteroni_2025}. 
 
The development of ROMs typically involves two distinct computational stages. The first is the offline stage, where high-fidelity simulation datasets are generated using high-performance computing (HPC) resources. The second is the online stage, where the reduced-order model is deployed for rapid prediction on comparatively less powerful computational platforms. Depending on the learning strategy adopted for the parametric space or unsteady dynamics, ROMs are broadly categorized into two classes: intrusive and non-intrusive approaches. In intrusive ROMs, the governing equations are projected onto a reduced space, and the resulting low-dimensional system is solved numerically as a time-marching problem. This approach often requires significant modification of the existing numerical solver. On the other hand, non-intrusive ROMs are entirely data-driven, where surrogate models are constructed using datasets from high-fidelity numerical simulations or experiments, thereby requiring minimal or no changes to the computational solver. Although intrusive ROMs, such as the projection-based Proper Orthogonal Decomposition (POD)-Galerkin method \cite{STABILE2018273,KaratzasStabileNouveauScovazziRozza2019b,StabileHijaziMolaLorenziRozza2017,stabile2020efficient}, have shown excellent performance in approximating system dynamics on a low-dimensional manifold, they often encounter numerical stability issues. This is particularly evident when the solution manifold exhibits strong nonlinearity, necessitating a large number of modes for accurate representation and increasing the computational cost of the reduced-order scheme. A detailed review of projection-based ROM techniques for parameterized dynamical systems is provided by Benner et al. \cite{BennerROM}. Furthermore, handling nonlinear terms in intrusive ROM frameworks demands additional computational effort, often addressed through hyper-reduction techniques such as the Empirical Interpolation Method (EIM) \cite{BARRAULT2004667}, Discrete Empirical Interpolation Method (DEIM) \cite{BDEIM2010}, GNAT \cite{carlberg2013gnat}, and other structure-preserving methods \cite{farhat2015structure}. Since intrusive ROMs enforce the governing equations in the reduced space, they generally provide accurate predictions even with a sparse training dataset. In recent years, significant progress has also been made in non-intrusive ROM approaches \cite{benner2021system,demo1,demo2019complete,GadallaCianferraTezzeleStabileMolaRozza2020}, particularly in handling ill-posed problems such as inverse modeling, regression, and classification in computational mechanics, as discussed by Schmidhuber et al. \cite{schmidhuber2015deep}. Several deep learning (DL) approaches have been proposed for such problems. For instance, Maulik and San \cite{maulik2019subgrid} and Duraisamy et al. \cite{parish2016paradigm} developed turbulence closure models using machine learning. Other notable contributions include the combination of POD with Long Short-Term Memory (LSTM) networks for modeling incompressible flows by Wang et al. \cite{wang2018model}, and the use of dimensionality reduction with deep learning for learning feature dynamics from noisy datasets by Lui and Wolf \cite{lui2019construction}. However, non-intrusive ROMs generally require a large amount of training data to achieve high predictive accuracy. LSTM networks have also been successfully applied to various problems in computational physics, such as aeroelastic and hydrodynamic applications by Halder et al. \cite{AHalder2020,halder2023deep}, and for modeling atmospheric turbulence by Rahman et al. \cite{rahman2019nonintrusive}.

Recent years have witnessed significant surge in the application of deep learning techniques for the improvement of the classical numerical techniques. In the areas of deep learning augmented Finite Element Method (FEM) techniques, notable works are \cite{samaniego2020energy,bhattacharya2020model}. Works on the application of deep learning methods for acceleration of linear solver includes Greenfeld et al., Xu et al. \cite{greenfeld2019learning, xu2023multigrid} to accelerate the multigrid solver, Wang et al. \cite{wang2022learning} for learning effective preconditioners for linear systems of equations. Recently, neural network-based correction strategies have been proposed to enhance the accuracy of classical numerical schemes. For instance, neural networks have been employed to learn Godunov-type corrections for approximate Riemann solvers within a bi-fidelity learning framework, leading to improved shock capturing and reduced numerical dissipation in hyperbolic conservation laws \cite{bar2019neural}.

Physics Informed Neural Network (PINN) and Neural Operators have gained the most notable research direction in the the recent times. In the past century, early work on solving partial differential equations
(PDE) using neural networks starts with Dissanayake and Phan-Thien
\cite{dissanayake1994neural}. With recent advancements in research of neural networks, this
topic got renewed attention with the work of Raissi et al.\cite{CRAISSI2019686}.
Physics-informed neural networks can circumvent the problem of large training
dataset requirements by introducing additional physics-based loss terms. This
physics-based loss term in neural networks can be computed from the governing
equations. Different types of neural networks can be considered in the context of PINN, such as artificial neural networks (ANN) by Raissi et al. \cite{CRAISSI2019686}, Waheed et al.\cite{bin2021pinneik}, Tartakovsky et al. \cite{tartakovsky2020physics} and
Schiassi et al.\cite{schiassi2021extreme}. Similarly, Gao et al. \cite{gao2021phygeonet} and Fang \cite{fang2019physics} coupled convolution neural networks (CNN) with physics constraints. A review of several developments in PINN is summarized by Cuomo et
al.\cite{cuomo2022scientific}. Zhang et
al.\cite{cheng2021deep} used a tensor differentiator to determine the derivatives of state space variables to couple the LSTM network with physics information arising from the governing equations and applied to nonlinear-structural problems. PINN usually considers an automatic differentiation (AD) based approach \cite{baydin2018automatic} for the computation of the gradients.

Now, we aim to focus on the integration of an existing forward CFD solver with the physics-informed neural network without much modification of the external CFD solver. Two notable works here are by \cite{shukla2025neurosem}, where spectral element-based solver Nektar++ is coupled with Physics Informed Neural Network, and by \cite{aulakh2022generalized}, where finite volume-based solver OpenFOAM is connected with the PINN. In both works, a PINN-based solution for a particular governing equation or a particular region of the computational domain is exported to the external CFD solver. In the current work, we propose a different route, importing the numerical residual from the governing equations discretized in an external CFD solver and considering it as a physics-based loss term in the PINN. We additionally propose a correction of the physics-based loss term in the case where the discretized form of the governing equations can not be included in the computational graph of the neural network. Additionally, we saw that large dimensionality often acts as a bottleneck for effective PINN training. Therefore, Chen et al. \cite{CHEN2021110666} proposed coupling of reduced order model approach with the PINN to circumvent this issue. However, they did not address the issue of whether the reduced order system is an external solver and can not be included in the computational graph of the neural network. Halder et al. \cite{halder2025physics} coined the term DisPINN for the first time based on a discretized physics-informed neural network. The main contribution of the current work is as follows:

\begin{itemize}
\item{ The current work proposes a novel approach where the existing forward solvers can be used in the PINN framework. %Therefore, the solvers can be used for ill-posed and inverse problems.
The motivation of the current work is to bridge the PINN
environment with the existing open-source CFD solvers for handling complex geometry and practical applications.}

\item{We propose an augmentation of the PINN algorithm, which incorporates an additional correction of the physics-based loss term in conventional PINN to seamlessly integrate it with any external solver (full order and reduced order form).}
\end{itemize}

The current work is organized in the following way: first, in \autoref{sec:GE}, the general formulation of the discretized governing equation is presented, and a simple neural network architecture is briefly recalled in \autoref{NN-coupling}. Furthermore, we propose an additional step of importing the physics-based residual and the Jacobian of the residual vector with respect to the output variable in the PINN solver if the external CFD solver is completely detached from the PINN environment to facilitate the inclusion of the discretized form of equation in the computational graph of the neural network. In \autoref{sec:testcase}, the coupling procedure of PINN with a finite volume formulation of governing equations in the forward solver and demonstrated on the nonlinear transport equation. Furthermore, the coupling approaches of the PINN with a reduced-order system are discussed next on viscous incompressible flow past a cylinder. Finally, the numerical experiments of the proposed methodologies on the test cases are shown in \autoref{sec:PINN}.

\section{\bf{Discretized Governing Equations}}\label{sec:GE}

    The current work proposes a novel approach to integrate an external CFD solver with a physics-informed neural network as $(\text{PINN + X})$, mentioned in \cite{shukla2025neurosem}, where X is the CFD solver. We introduce a library called PINNFOAM for coupling PINN with the external CFD solver OpenFOAM \cite{openfoam}. To establish a foundation, we discuss the discretization approach to the governing equations in both full-order and reduced-order forms. Various discretization methods, such as finite volume, finite difference, and finite element approaches, are commonly used in different CFD solvers. Here, the finite volume approximation is briefly recalled; for more details, the reader may see \cite{barth2018finite,ferziger2002computational}. We particularly demonstrate the cell-centred finite volume ($V_p$ centred at $p$)  formulation for a general transport equation of a scalar property $\phi$ between time instant $t$ to $t+\Delta t$ as follows:

\begin{equation} \label{eq:gen_eqn}    
\begin{aligned}
\int_t^{t+\Delta t}\quad\left[\int_{V_p} \frac{\partial \rho \phi}{\partial t} d V+\int_{V_p} \nabla \cdot(\rho \mathbf{v} \phi) d V-\int_{V_p} \nabla \cdot\left(\rho \Gamma_\phi \nabla \phi\right) d V\right] d t \quad= \\
\int_t^{t+\Delta t}\left(\int_{V_p} S_\phi(\phi) d V\right) d t ,
\end{aligned}
\end{equation}

where $\rho$ is the density, $\textbf{v}$ is the velocity vector (considered fixed in this context, without dependence on $\phi$), $\Gamma_\phi$ is the diffusivity, $S_\phi$ is the source term. Various spatial and temporal discretization schemes can be employed for $\phi$. However, for demonstration, we will consider the implicit Euler time-discretization approach here. Now, the Gauss divergence theorem is applied to convert volume integrals into surface integrals. Consequently, the discretized form of \autoref{eq:gen_eqn} can be written as follows:

\begin{equation} \label{eq_dis}
\begin{aligned}
\rho {\phi}^{(n+1)} = \rho {\phi}^{(n)} + dt(\sum_f A_f \cdot(\rho \mathbf{v}\phi^{n+1})_f-\sum_f\left(\rho \Gamma_\phi\right)_f A_f \cdot(\nabla \phi^{n+1})_f-S_u V_p+S_p V_p \phi^{n+1}_p),
\end{aligned}
\end{equation}

where $A_f$ is an outward pointing area vector for the cell face, where $f$ is the value of the quantity at the centre of the cell face. The source term can be linearized as \( S_\phi(\phi) = S_u + S_p \phi \), where \( S_u \) and \( S_p \) are the constant and linear coefficients, respectively. This formulation allows the variable \( \phi \) to be treated using its value from the subsequent time step. Alternatively, the source term can be treated explicitly using its value at the current time level, i.e., \( S_\phi(\phi^n) \). The superscript $n$ is the current time step for the unsteady problem at hand, where we aim to find the solution $\phi$ at the subsequent time step $\left(n+1\right)$ by solving \autoref{eq_dis}, which can be further simplified to \autoref{eq_lin}.

\begin{equation} \label{eq_lin}
\begin{aligned}
a_P \phi_{p}^{n+1}= \sum_{k} a_i \phi_{i}^{n+1} + b_{p},
\end{aligned}
\end{equation}

here, \( p \) and \( i \) represent the values of \( \phi \) corresponding to the cell centered at \( p \) and its \( k \) neighboring cells. The coefficients \( a_p \) and \( a_i \) are constant, while \( b_p \) on the right-hand side includes the values of \( \phi \) at the \( n^{th} \) time step, along with the constant part of the source terms. \autoref{eq_lin} for all cells can be collectively expressed as a system of linear equations, \eqref{linear_Mat_GEN}, where matrix \( \mathbf{A} \) consists of the coefficients \( a_p \) and \( a_i \). Vector \( \mathbf{\Phi}^{(n+1)} \) contains the unknown values \( \phi^{n+1} \), and vector \( \mathbf{b} \) holds the terms \( b_p \) for all cells in the computational domain. For various choices of temporal and spatial discretization schemes applied to $\phi$ in \autoref{eq_dis}, one can still derive a linear system of equations of the form given in \autoref{linear_Mat_GEN}, albeit with different forms of the coefficient matrix $\mathbf{A}$ and the right-hand side vector $\mathbf{b}$.

\begin{equation}\label{linear_Mat_GEN}
\mathbf{A} {\mathbf{\Phi}^{n+1}}=\mathbf{b}.
\end{equation}

The vector, $\mathbf{\Phi}$, can be computed using different types of linear solvers as mentioned in \cite{quarteroni2010numerical} or by minimizing the residual $\mathbf{R}$ in \autoref{linear_Mat_res}, using an optimization approach \cite{bottou2018optimization} in a neural network environment.

\begin{equation}\label{linear_Mat_res}
\mathbf{R}_{(n+1)} = \mathbf{A} {\mathbf{\Phi}^{n+1}}-\mathbf{b}.
\end{equation}

In the following sections, we will describe how a neural network can be trained considering $\mathbf{R}$ as a physics-based loss term in addition to a data-driven loss term as mentioned in \cite{CRAISSI2019686}. 

\section{Neural Network and Integration with External CFD solver} \label{NN-coupling}

Missing data in spatial or temporal location can cause a trained neural network (NN) to deviate from the ground truth. In the context of computational fluid dynamics (CFD), missing data within the solution domain may result from a coarse mesh region or larger time steps in the temporal direction. Additionally, training data loss can occur in experimental settings due to sparse probe locations or collecting prob data at large time intervals. The accuracy of the neural network depends on the fidelity of the data. If the resolution of the available data is not enough to represent all the features of the ground truth solution, the NN prediction will be erroneous. A sparse spatio-temporal dataset also causes erroneous prediction at the future time steps or outside the training parameter zone if the available training data is not sufficient to capture all the features of the spatio-temporal dynamics. This problem can be addressed by the introduction of the a apriori knowledge from governing equations to augment the learning from the sparse available data. Therefore, a hybrid methodology is presented combining physics-based learning with data-driven learning. We demonstrate an effective neural network training with physics-based loss terms from an external CFD solver while considering sparse numerical datasets in the spatio-temporal locations for the data-driven loss term. First, a neural network architecture such as ANN will be briefly introduced, followed by the integration with the discretized governing equation in the framework of the physics-informed neural network \cite{CRAISSI2019686} termed as ANN-DisPINN. ANN formulates a nonlinear functional relationship between the input and output.
In \autoref{eq:ann}, $\textbf{z}$ can be considered as input of an unsteady problem consisting of $N_t$ time instants (i.e., $\textbf{t} = \left [ t_1, t_2,...., t_{N_{t}}  \right ]$) and $m$ parameter values of $\mathbf{\Gamma}_{\phi}$ (i.e., $\mathbf{\Gamma}_{\phi} = \left[ {\Gamma}_{\phi,1}, {\Gamma}_{\phi,2},....,  {\Gamma}_{\phi,m}\right]$), whereas the variable
$\mathbf{\Phi}$ can be the neural network output. Therefore, the size of input and the output will be $\mathbf{z} \in \mathbb{R}^{mN_t \times 2}$ and $\mathbf{\Phi} \in \mathbb{R}^{mN_t \times \text{N}_\phi}$ where $\text{N}_\phi$ is the mesh size of the computational domain. The input-output formulations of an Artificial Neural Network (ANN) can be expressed mathematically as follows:

\begin{equation}\label{eq:ann}
\begin{aligned}
\mathbf{\Phi}=f_l\left(\boldsymbol{W}_l f_{l-1}\left(\cdots f_1\left(\boldsymbol{W}_1 \boldsymbol{z}+\boldsymbol{c}_1\right)+\cdots+\boldsymbol{c}_{l-1}\right)+\boldsymbol{c}_l\right),
\end{aligned}
\end{equation}

where $l$ is the number of layers. $f_{1}$, $f_{2}$ and $f_{l}$ are the activation function. 
the $\boldsymbol{W}_1$, $\boldsymbol{W}_2$ and $\boldsymbol{W}_l$ are the weight matrices and
$\boldsymbol{c}_1$, $\boldsymbol{c}_2$ and $\textbf{c}_l$ are the bias vectors. For large values of $\text{N}_\phi$, the neural network often has a problem of scalability. Therefore, in the following sections, we will demonstrate how to couple the Reduced Order Model approach \cite{STABILE2018273} with the neural network. The classical artificial neural network is trained by minimizing only the data-driven loss terms, which can be defined as follows: 

\begin{equation}\label{eq:ldata}
\begin{aligned}
\mathrm{L}_{\text {Data }}=\frac{1}{N_{\text {Data }}} \sum_{i=1}^{N_{\text {Data }}}\left(\mathbf{\Phi}_{i,\text{pred}}-\mathbf{\Phi}_{i,\text{actual}}\right)^2,
\end{aligned}
\end{equation}

where, $\mathbf{\Phi}_{i,\text{pred}}$ is the prediction from the neural network as shown in \autoref{eq:ann}, and $\mathbf{\Phi}_{i,\text{actual}}$ are the benchmark numerical data available at total number of time instants, $N_{\text{Data}}$. In a general network, the training algorithm aims to optimize the set of hyperparameters, $\theta=\left[b_i, W_i\right]_{i=1}^l$, by minimizing the loss term $\mathrm{L}_{\text {Data }}$. The optimized $\theta$, i.e. $\theta^*$ is defined as: $\theta^*=\arg \min_{\theta} [\text{L}_{\text{Data}}(\theta)]$. In the present work, we have considered a first-order gradient-based optimization approach, ADAM \cite{adam2014method}. As previously discussed, when the number of available data points $N_\text{Data}$ is limited, conventional neural networks struggle to capture all the spatio-temporal features. To address this limitation, a hybrid approach is introduced, which incorporates a physics-driven loss term alongside the data-driven loss term. The physics-driven loss term $\mathrm{L}_{\text{eqn}}$ can be derived as follows: 

\begin{equation}\label{eq:leqn}
\begin{aligned}
\mathrm{L}_{\text{eqn}}=\frac{1}{N_{\text{eqn}}} \sum_{i=1}^{N_{\text{eqn}}} {\mathbf{R}_i}^2,
\end{aligned}
\end{equation}

where, the physics-based residual, $\mathbf{R}_i$ arises from \autoref{linear_Mat_GEN}, $i$ denotes the physics based residual $\textbf{R}$ obtained from \autoref{linear_Mat_res} after the finite volume based discretization of governing equations corresponding to the time instants ($t_i$) and parameter value $\Gamma_{\phi,i}$. $N_{\text {eqn}}$ is the number of points available for the computation of $(\text{L}_\text{Eqn})$. The total loss term can be written as:

\begin{equation}\label{eq:ltot1}
\begin{aligned}
\text{L} = {\lambda_{1}}\text{L}_{\text{eqn}} + {\lambda_{2}}\text{L}_{\text{Data}},
\end{aligned}
\end{equation}

where, $\lambda_1$ and $\lambda_2$ are the constant values multiplied with the $\text{L}_{\text{eqn}}$ and $\text{L}_{\text{Data}}$ respectively. Given $\mathbf{\Phi}_{\text {tot, pred }}=\left[\mathbf{\Phi}_{\text {pred}, 1}, \mathbf{\Phi}_{\text {pred}, 2}, \mathbf{\Phi}_{\text {pred}, 3 \ldots \ldots .} \ldots \mathbf{\Phi}_{\text{pred}, mN_t}\right]$, the derivative of ${\text{L}}_\text{eqn}$ concerning the hyper-parameters can be computed as follows:

\begin{equation}\label{eq:dldtheta}
\begin{aligned}
\partial \mathrm{L}_{\text {eqn }} / \partial \mathrm{\theta}=\left(\frac{1}{N_{\text {eqn }}} \sum_1^{N_{\text {eqn}}} 2 \textbf{R}_i \frac{\partial \textbf{R}_i}{\partial \theta}\right)=\left(\frac{1}{N_{\text {eqn }}} \sum_1^{N_{\text {eqn }}} 2 \textbf{R}_i \frac{\partial \textbf{R}_i}{\partial \mathbf{\Phi}_{\text{tot, pred}}} \frac{\partial \mathbf{\Phi}_{\text{tot, pred}}}{\partial \theta}\right).
\end{aligned}
\end{equation}

Now, if the physics-based residual, $\mathbf{R}_i$ is computed from an external solver, it is challenging to access the discretized form directly from the PINN solver, therefore we can not compute $\partial \mathrm{L}_{\text {eqn }} / \partial \theta$ from the computational graph using back-propagation. However if both the $\textbf{R}_i$ and $\frac{\partial \textbf{R}_i}{\partial \mathbf{\Phi}_{\text{tot, pred}}}$ are imported from the external solver, the modified physics-based loss therm, $\mathrm{L}_{\text{Dis}}$  that can be passed through the backpropagation is written as follows:
\begin{equation}\label{eq:LDIS}
\begin{aligned}
\mathrm{L}_{\text{Dis}}=\left(\frac{1}{N_{\text {eqn }}} \sum_1^{N_{\text {eqn }}}\left[2 \textbf{R}_i \frac{\partial \textbf{R}_i}{\partial \mathbf{\Phi}_{\text{tot, pred}}}\right]_{\operatorname{detached}} \mathbf{\Phi}_{\text{tot, pred}}\right).
\end{aligned}
\end{equation}

Since both the $\mathbf{R}_i$ and $\frac{\partial \textbf{R}_i}{\partial \mathbf{\Phi}_{\text{tot, pred}}}$ are imported from the external solver, therefore, detached from the computational graph, the derivative concerning the neural network parameter is only computed on $\mathbf{\Phi}_{\text{tot, pred}}$ outside the bracketed part $\left[
 \right]_{\text{detached}}$. Hence, the loss term that needs to be passed through the backpropagation is ($\text{L}= \lambda_1\mathrm{L}_{\text {Dis }}+\lambda_2\mathrm{L}_{\text {Data }}$). Finally, the first-order gradient-based optimization approach is considered in the current work, as mentioned earlier, to update the weight matrices and bias vectors based on $\partial \mathrm{L} / \partial \mathrm{\theta}$, following the relationship:

\begin{equation}\label{eq:pupdate}
\begin{aligned}
\theta_{new}= \theta_{old} - \alpha \frac{\partial \text{L}}{\partial \theta}, 
\end{aligned}
\end{equation}

where, $\theta_{new}$ and $\theta_{old}$ are the updated and old parameter values and $\alpha$ is the learning rate. It is important to note that the above formulation of the loss term \( \mathrm{L}_{\text{Dis}} \) is applicable under the assumption that first-order optimization methods are used. Now, if the residual $\mathbf{R}_i$ is computed from the external solver, Jacobian, $J = \frac{\partial \mathbf{R}_i}{\partial \Phi_{\text{tot, pred}}}$ needs to be updated at every epoch. However, $\frac{\partial \mathbf{R}}{\partial \Phi_{\text{tot, pred}}}$ is a sparse matrix, and the overall computational expense of PINN solver can increase several folds if the $J$ term is updated every epoch. It contains residual from the entire computational domain. Therefore, the computational expenditure can be reduced by updating the Jacobian term, $J$ only at a certain epoch interval. In the current work, the $J$ term is computed using finite-difference (FD), and no additional modification is required in the external solver, demonstrating the potential of our approach for seamless coupling of PINN with any other external forward solver, such as the OpenFOAM here. The algorithm of coupling PINN with an external solver is discussed in the authors' previous work \cite{halder2025physics}. In the following section, we demonstrate that even with access only to the system matrix $\mathbf{A}$ and the forcing term $\mathbf{b}$, as shown in \autoref{linear_Mat_res}, a neural network can still be coupled with an external CFD solver---even if the full discretized form cannot be included in the computational graph.

\section{Testcases} \label{sec:testcase}

The procedure for coupling an external CFD solver, where the governing equations are discretized using a finite-volume formulation with a neural network framework is demonstrated through two test cases: a full-order simulation of the Nonlinear Transport Equation model and a reduced-order model of viscous incompressible flow past a cylinder.
   
\subsection{Full Order System - Nonlinear Transport Equation Model}\label{Burgers_FOM}

 The scalar term $\phi$ is replaced with the vector $\mathbf{u}$ in \autoref{eq:gen_eqn} to obtain the nonlinear transport equation with the physical parameter, $v$. We apply finite-volume-based formulations for the governing equation. We consider a Eulerian frame on a space-time domain with $(\mathrm{x}, t) \in[0,L]^d \times[0,t]$ with $d = 2$. A square domain denoted by $\mathbf{x} \in[x_{1},x_{2}]^2$ is taken as an initial condition with uniform $\mathbf{u}$ value of $1$. The boundary values at $\mathbf{x} \in \partial[0,L]^2$ is taken as $0$. 

\begin{equation}\label{eq:Burgers}
\begin{aligned}
\begin{array}{ll}
\partial_t \mathbf{u}+\frac{1}{2} \nabla \cdot(\mathbf{u} \otimes \mathbf{u})=v \Delta \mathbf{u} & (\mathrm{x}, t) \in[0,L]^2 \times[0,t], \\
\mathbf{u}(\mathbf{x}, 0)=1.0 & \mathbf{x} \in[x_{1},x_{2}]^2, \\
\mathbf{u}(\mathbf{x}, t)=0 & (\mathbf{x}, t) \in \partial[0,L]^2 \times[0,t].
\end{array}
\end{aligned}
\end{equation}

In the external solver, the finite-volume approximation of \autoref{eq:Burgers} for a control volume centred at $p$ is as follows: 

\begin{equation}\label{eq:Burgers_fv}
\begin{aligned}
\int_{V_p} \frac{\partial \mathbf{u}}{\partial t} d V+\int_{V_p} \frac{1}{2} \nabla \cdot(\mathbf{u} \otimes \mathbf{u}) d V-\int_{V_p} \nabla \cdot\left(v \nabla \mathbf{u}\right) d V=0.
\end{aligned}
\end{equation}

 The temporal and spatial discretization is carried out in the external CFD solver. Similar to \autoref{eq_dis}, the discretized form of the spatial derivative terms in \autoref{eq:Burgers_fv} can be formulated as follows:

\begin{equation}\label{eq:Burgers_discretized}
\begin{aligned}
\frac{\partial}{\partial t} \int_{V_p} \mathbf{u} d V+\frac{1}{2}\sum_f A_f \cdot(\mathbf{u} \otimes \mathbf{u})_f-\sum_f A_f \cdot\left(v \nabla \mathbf{u}\right)_f=0.
\end{aligned}
\end{equation}

For the temporal derivative, we have considered an Euler time discretization scheme. The implicit formulation of the \autoref{eq:Burgers_discretized} can be written as follows:  

\begin{equation} \label{eq_dis1}
\begin{aligned}
{\mathbf{u}}_{p}^{(n+1)} = {\mathbf{u}}_{p}^{(n)} + dt(\frac{1}{2}\sum_f A_f \cdot( \mathbf{u}^{n+1}\mathbf{u}^{n+1})_f-\sum_f v A_f \cdot(\nabla \mathbf{u}^{n+1})_f,
\end{aligned}
\end{equation}

Although the convective term is nonlinear here, a linear system can now be formulated by linearizing the nonlinear convective term in the following way:

\begin{equation} \label{eq_dis2}
\begin{aligned}
{\mathbf{u}}_{p}^{(n+1)} = {\mathbf{u}}_{p}^{(n)} + dt(\sum_f A_f \cdot( \mathbf{u}^{n}\mathbf{u}^{n+1})_f-\sum_f v A_f \cdot(\nabla \mathbf{u}^{n+1})_f,
\end{aligned}
\end{equation}

However, in the convective term, both the $\mathbf{u}$ values can also be considered at $n^{\text{th}}$ time step. \autoref{eq_dis2} for all the cells in the computational domain can be written in Matrix format as shown in \autoref{sec:GE} as matrix $\mathbf{A}_{(n+1)}$ and
forcing vector $\mathbf{b}_{(n+1)}$ denotes  $\mathbf{A}$ and $\mathbf{b}$ associated with $\mathbf{R}_{(n+1)}$. $U^{n+1}$ is the vector consisting of the values of $\textbf{u}^{(n+1)}$ at all the computational cells. 

\begin{equation}\label{linear_Mat}
\mathbf{R}_{(n+1)} = \mathbf{A}_{(n+1)} {U^{n+1}}-\mathbf{b}_{(n+1)}.
\end{equation}

In section \autoref{NN-coupling}, we demonstrate that a PINN can be integrated with an external CFD solver even if it only has access to the residual $\mathbf{R}$ and not the full discretized form of the equations, such as \autoref{eq_dis1} and \autoref{eq_dis2} can not be included in the computational graph of the neural network. Now for the full order nonlinear transport Equation, we will formulate the $\text{L}_{\text{Dis}}$ shown in \autoref{eq:LDIS}, when $\mathbf{A}_{(n+1)}$ and vector $\mathbf{b}_{(n+1)}$ can be imported from the external solver. Time ($\mathbf{t}$) is an input to the neural network, which consists of $[{t}_1, {t}_2 ..... {t}_{N_t}]$. where $N_{t}$ is the total time instants. The output of the neural network associated with $N_{t}$ time instants is as follows:   

\begin{equation}\label{eq:Upred_vec}
\begin{aligned}
& \textbf{U}_{\text {pred }}=\left[U_{\text {pred}, 1}, U_{\text {pred}, 2}, U_{\text {pred},3} \ldots \ldots . U_{\text {pred},{N_t}}\right].
\end{aligned}
\end{equation} 
\\
The obtained $\textbf{U}_{\text {pred}}$ is passed to the external CFD solver to obtain the residual at $(n+1)^{\text{th}}$ time step as follows:

\begin{equation}\label{eq:Burgers_Linear_Mat}
\begin{aligned}
\mathbf{R}_{(n+1)} = \mathbf{A}_{(n+1)}U_{\text {pred},{(n+1)}}-\mathbf{b}_{(n+1)}.
\end{aligned}
\end{equation}

Now, we can write $\mathbf{R}_{\text{tot}}$ associated with numerical residual $\mathbf{R}_{i}$ for all the $N_{\text{eqn}}$ points, as $\mathbf{R}_{\text{tot}} = \left[\mathbf{R}_{2}, \mathbf{R}_{3}, \mathbf{R}_{N_\text{eqn}} \right]$. The physics based loss term $\text{L}_\text{eqn}$ can be written as follows:

\begin{equation}\label{eq:Leqn_NTE}
\begin{aligned}
\mathrm{L}_{\text{eqn}}=\frac{1}{N_{\text{eqn}}} \sum_{i=1}^{N_{\text{eqn}}} {\mathbf{R}_i}^2 = \frac{1}{N_{\text{eqn}}} \sum_{i=1}^{N_{\text{eqn}}} {({\mathbf{A}_{i}U_{\text {pred},{i}}-\mathbf{b}_{i}}})^2.
\end{aligned}
\end{equation}

If the governing equation is linear, ($\mathbf{A}, \mathbf{b}$ is constant, not dependent on $\mathbf{u}$), $\mathrm{L}_{\text{eqn}}$ can alone be considered as the physics-based residual of the neural network without considering any correction as in \autoref{eq:LDIS}. However, for the nonlinear transport equation, as considered in the current work, either $\mathbf{A}$ or $\mathbf{b}$ or both the terms could be dependent on the variables $\mathbf{u}$ and we have to obtain the modified physics-based loss term, $\mathrm{L}_{\text{Dis}}$, to be employed in the optimization machinery of the neural network as demonstrated in \autoref{eq:LDIS}. To obtain the Jacobian $J =\frac{\partial \mathbf{R}_{\text{tot}}}{\partial \mathbf{U}_{\text{pred}}}$ term for the computation of the $\mathrm{L}_{\text{Dis}}$ as follows:

\begin{equation}\label{eq:Burgers_grad1}
\begin{aligned}
\frac{\partial \mathbf{R}_{\text{tot}}}{\partial \mathbf{U}_{\text {pred}}}=\left[\begin{array}{cccc}
\frac{\partial \mathbf{R}_2}{\partial U_{\text {pred}, 1}} & \frac{\partial \mathbf{R}_2}{\partial U_{\text {pred}, 2}} & \ldots & \frac{\partial \mathbf{R}_2}{\partial U_{\text {pred}, N_{\text{eqn}}}} \\
\frac{\partial \mathbf{R}_3}{\partial U_{\text {pred}, 1}} & : & . . & \frac{\partial \mathbf{R}_3}{\partial U_{\text {pred}, N_{\text{eqn}}}} \\
: & : & . . & : \\
\frac{\partial \mathbf{R}_{N_\text{eqn}}}{\partial U_{\text {pred}, 1}} & : & . . & \frac{\partial \mathbf{R}_{N_\text{eqn}}}{\partial U_{\text {pred}, {N_\text{eqn}}}}
\end{array}\right].
\end{aligned}
\end{equation}

Now we can write, 

\begin{equation}\label{eq:LDIS_NTE}
\begin{aligned}
\mathrm{L}_{\text{Dis}}=\left(\frac{1}{N_{\text {eqn }}} \sum_1^{N_{\text {eqn }}}\left[2 \textbf{R}_i \frac{\partial \textbf{R}_i}{\partial \mathbf{U}_{\text{pred}}}\right]_{\operatorname{detached}} \mathbf{U}_{\text{pred}}\right).
\end{aligned}
\end{equation}

The Jacobian $J$ term, $\frac{\partial \mathbf{R}_{\text{tot}}}{\partial \mathbf{U}_{\text {pred}}}$ can further be simplified to a sparser matrix when the time and space discretization schemes are known. For example, using Euler time discretization, the residual $\mathbf{R}_{(n+1)}$ depends only on the $(n+1)^{\text{th}}$ and $n^{\text{th}}$ time steps for \autoref{eq_dis2}. This structure ensures that each row of the Jacobian contains only two non-zero entries, with all other elements being zero. Now for $\mathbf{R}_{(n+1)}$, while \autoref{eq_dis2} is considered for discretization, $\mathbf{A}_{(n+1)} = f(U_{\text {pred},n})$, $\mathbf{b}_{(n+1)} = f(U_{\text {pred},n})$. $\frac{\partial \mathbf{R}_{(n+1)}}{\partial U_{\text {pred}, n+1}}$ can be written as follows:

\begin{equation}\label{eq:grad_burgers_main}
\begin{aligned}
& \frac{\partial \mathbf{R}_{(n+1)}}{\partial U_{\text {pred}, (n+1)}}=\frac{\partial\left(\mathbf{A}_{(n+1)} U_{\text {pred}, (n+1)}+ \mathbf{b}_{(n+1)}\right)}{\partial U_{\text {pred}, (n+1)}}, \\
& = \mathbf{A}_{(n+1)},
\end{aligned}
\end{equation}

Therefore, the $J$ can be further simplified as follows:

\begin{equation}\label{eq:grad_burgers}
\begin{aligned}
\frac{\partial \mathbf{R}_\text{tot}}{\partial \mathbf{U}_{\text {pred }}}=\left[\begin{array}{cccc}
\frac{\partial \mathbf{R}_2}{\partial U_{\text {pred}, 1}} & \mathbf{A}_2 & \ldots . & 0 \\
0 & : & . . & 0 \\
: & : & . . & : \\
0 & : & \frac{\partial \mathbf{R}_{N_\text{eqn}}}{\partial \mathbf{U}_{\text {pred}, {N_\text{eqn}}-1}} & \mathbf{A}_{\text{Neqn}}
\end{array}\right],
\end{aligned}
\end{equation}

\begin{figure}[ht]
\centering
\begin{subfigure}[b]{1\textwidth}
\centering
% \hspace{-3cm}
{\label{fig:burgers2}\includegraphics[width=0.8\linewidth]{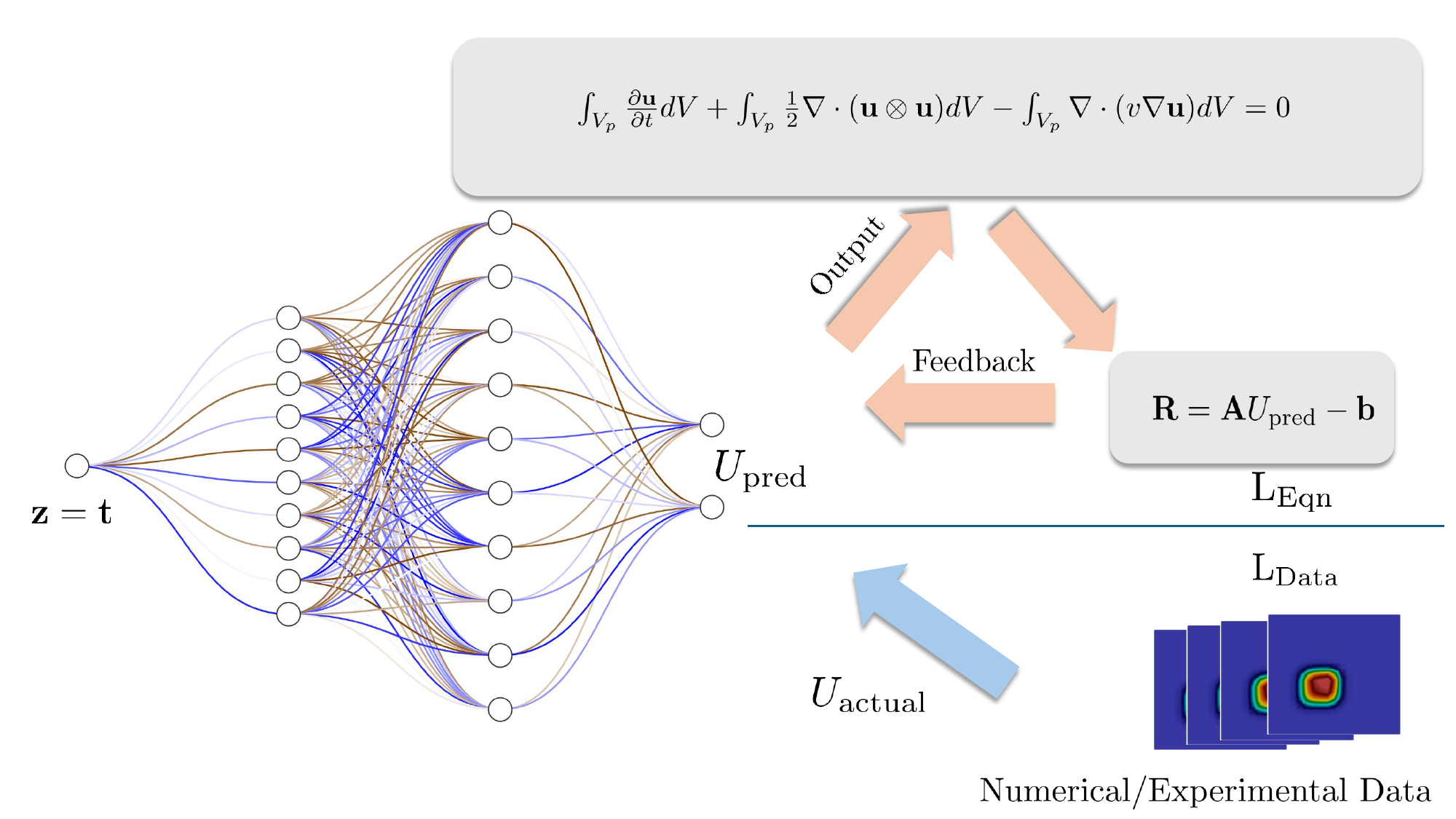}}
\end{subfigure}
\caption{Schematic of the discretized physics-based PINN and coupling procedure with any external CFD solver where, $\text{L}_\text{eqn}$ and $\text{L}_\text{Data}$ are physics-based and data-driven loss terms, respectively. }
\label{fig:ANN-DisPINN_Full}
\end{figure}

As different discretisation schemes are employed for the temporal and spatial derivatives, the forms of \autoref{eq:grad_burgers_main} and \autoref{eq:grad_burgers} will also be modified accordingly. Now, we demonstrate the algorithm for coupling the full order model from an External CFD solver with the PINN, where the neural network can only import $\mathbf{A}$ and $\mathbf{b}$. The Jacobian $J$ termed as $\frac{\partial \mathbf{R}_{\text{tot}}}{\partial \mathbf{U}_{\text {pred}}}$ is updated from the external CFD solver at epoch interval termed as $k_{int}$ which plays a crucial role deciding the accuracy of the predicted solution and training time.

\SetKwComment{Comment}{/* }{ */}

\begin{algorithm}[hbt!] 
\caption{DisPINN pseudo-code for coupling external full-order forward solver}\label{alg:algo-burgers}
\KwData{$U_\text{actual}, \textbf{z}, N_\text{Data}, N_\text{eqn}$}

$[{\boldsymbol{W}},{\boldsymbol{c}}] \gets \textbf{INIT} ([\boldsymbol{W},\boldsymbol{c}])$ \Comment*[r]{Initialize Weight and bias}
\While{$\text{epoch} \leq \text{Total Epoch No.}$}{
  
  $U_\text{pred} \gets \text{NN}(\boldsymbol{W},\boldsymbol{c}, \textbf{z})$ \Comment*[r]{prediction from chosen network $\text{NN}$}
  $\text{Ext. Solver} \gets U_\text{pred,i}$ \Comment*[r]{pass to external solver}
  $ \mathbf{R}_i, \mathbf{A}_i, \mathbf{b}_i \gets \text{Ext. Solver}$  \Comment*[r]{pass to PINN solver from external solver}
  $\mathrm{L}_{\text{eqn}} \gets \frac{1}{N_{\text {eqn }}}\sum_{i=1}^{N_{\text{eqn}}} \mathbf{R}_i^2  \text{or},\frac{1}{N_{\text{eqn}}} \sum_{i=1}^{N_{\text{eqn}}} {({\mathbf{A}_{i}U_{\text {pred},{i}}-\mathbf{b}_{i}}})^2$ \Comment*[r]{compute physics-based loss term}
  $\mathrm{L}_{\text {Data }} \gets \frac{1}{N_{\text {Data }}}\sum_{i=1}^{N_{\text {Data }}}\left(U_{\text {pred,i }}-U_{\text {actual,i }}\right)^2$ \Comment*[r]{compute data-based loss term}

  \eIf{$\text{rem}(epoch, k_{int}) = 0$}{  
    $\frac{\partial \mathbf{R}_\text{tot}}{\partial \textbf{U}_\text{pred}} \gets \text{Ext. Solver}$ \Comment*[r]{compute Jacobian from external solver}
    $J_\text{new,i} \gets \frac{\partial \mathbf{R}_i}{\partial \textbf{U}_\text{pred}}$\; 
    $\mathrm{L}_{\text{Dis}} \gets (\frac{1}{N_{\text {eqn}}} \sum_1^{N_{\text {eqn }}}2 \textbf{R}_i J_\text{new,i} \textbf{U}_\text{pred})$  \Comment*[r]{compute additional physics based loss term for backpropagation}
    }
    {$J_i \gets J_\text{new,i}$ \;
    $\mathrm{L}_{\text{Dis}} \gets (\frac{1}{N_{\text {eqn }}} \sum_1^{N_{\text {eqn }}}2 \textbf{R}_i J_i \textbf{U}_\text{pred})$ \;
  }

  $\mathrm{L} \gets (\mathrm{L}_{\text{Dis}} + \mathrm{L}_{\text {Data }})$ \Comment*[r]{compute total loss term for backpropagation}

  $\Delta {\boldsymbol{W}} \leftarrow-\alpha \mathrm{G}_{\mathrm{ADAM}}\left(\nabla_{{\boldsymbol{W}}} L\right), \Delta \boldsymbol{c} \leftarrow-\alpha \mathrm{G}_{\mathrm{ADAM}}\left(\nabla_{\boldsymbol{c}} L\right)$ \; \Comment*[r]{compute weight matrices and bias vector update}
  $\boldsymbol{W} \leftarrow {\boldsymbol{W}}+\Delta \boldsymbol{W}, \boldsymbol{c} \leftarrow \boldsymbol{c}+\Delta {\boldsymbol{c}}$ \;
}
\end{algorithm}

\subsection{Reduced Order System - Incompressible Navier Stokes Equation} \label{ROM-NS}

The neural network size often scales with the mesh size of the physical problem. However, for large mesh size values, the neural networks can have difficulties with scalability. Therefore, in this section, we will demonstrate the coupling of the reduced-order formulation of a full-order system with the physics-informed neural network framework. We have focused on a mathematical problem given by the unsteady, incompressible, parameterized Navier-Stokes equation. We consider an Eulerian frame in the space-time domain $Z=V_{p} \times[0, T] \subset \mathbb{R}^d \times \mathbb{R}^{+}$ with $d=2$. The vectorial field $\boldsymbol{u}: Z \rightarrow \mathbb{R}^d$ and pressure field $\boldsymbol{p}: Z \rightarrow \mathbb{R}$. The governing equation can be written as follows:

\begin{equation}\label{eq:NS-FOM}
\begin{aligned}
\begin{array}{ll}
\boldsymbol{u}_{\boldsymbol{t}}+\boldsymbol{\nabla} \cdot(\boldsymbol{u} \otimes \boldsymbol{u})-\boldsymbol{\nabla} \cdot 2 \nu \boldsymbol{\nabla}^{\boldsymbol{s}} \boldsymbol{u}=-\boldsymbol{\nabla} \boldsymbol{p} & \text { in } Z, \\
\boldsymbol{\nabla} \cdot \boldsymbol{u}=\mathbf{0} & \text { in } Z, \\
\boldsymbol{u}(t, x)=\boldsymbol{g}(\boldsymbol{x}) & \text { on } \Gamma_{\text{In}} \times[0, T], \\
\boldsymbol{u}(t, x)=\mathbf{0} & \text { on } \Gamma_0 \times[0, T], \\
(\nu \nabla \boldsymbol{u}-p \boldsymbol{I}) \boldsymbol{n}=\mathbf{0} & \text { on } \Gamma_{\text{Out}} \times[0, T], \\
\boldsymbol{u}(0, \boldsymbol{x})=\boldsymbol{k}(\boldsymbol{x}) & \text { in } T_0,
\end{array}
\end{aligned}
\end{equation}

where, $\Gamma=\Gamma_{\text{In}} \cup \Gamma_0 \cup \Gamma_{\text {Out }}$ is the boundary of $V_{p}$ and, is composed by three different parts $\Gamma_{\text{In}}, \Gamma_{\text {Out }} \text { and } \Gamma_0$ which are inlet, outlet and physical walls.  The kinematic viscosity ($\nu$) gives the parameter dependency. The function $\boldsymbol{g}(x)$ represents the boundary conditions for the non-homogeneous boundary and $\boldsymbol{k}(x)$ denotes the initial condition for the velocity at $t=0$, denoted by $T_0$. It is also assumed that the boundary conditions $f$ are not dependent on time. The finite volume formulation of \autoref{eq:NS-FOM} can be written for a volume $V_p$ centred at $p$ as follows:

\begin{equation}\label{eq:NS-FV1}
\begin{aligned}
\int_{V_p} \frac{\partial \boldsymbol{u}}{\partial t} \mathrm{~d} V+\int_{V_p} \boldsymbol{\nabla} \cdot(\boldsymbol{u} \otimes \boldsymbol{u}) \mathrm{d} V 
& -\int_{V_p} \boldsymbol{\nabla} \cdot 2 \nu \boldsymbol{\nabla}^s \boldsymbol{u} \mathrm{~d} V = -\int_{V_p} \boldsymbol{\nabla} \boldsymbol{p} \mathrm{~d} V,
\end{aligned}
\end{equation}
\begin{equation}\label{eq:NS-FV2}
\begin{aligned}
\int_{V_p} \boldsymbol{\nabla} \cdot \boldsymbol{u}\mathrm{~d} V=0.
\end{aligned}
\end{equation}

Within the finite volume discretization, all the divergence terms are re-written in terms of fluxes over boundaries of each finite volume, making use of the Gauss's theorem:

\begin{equation}\label{eq:NS-FV3}
\begin{aligned}
\frac{\partial}{\partial t} \int_{V_p} \boldsymbol{u} d V+\sum_f A_f \cdot(\boldsymbol{u} \otimes \boldsymbol{u})_f-\sum_f A_f \cdot(2\nu \boldsymbol{\nabla}^s \boldsymbol{u})_f= - \sum_f A_f  \cdot \boldsymbol{p}_f.
\end{aligned}
\end{equation}

The term originated from the divergence of the velocity can be discretized as follows:

\begin{equation}\label{eq:continuity-FV}
\begin{aligned}
\sum_f A_f  \cdot \boldsymbol{u}_f = 0.
\end{aligned}
\end{equation}

In \autoref{eq:NS-FV3}, all the coefficients that multiply the acceleration term $\frac{\partial}{\partial t} \int_{V_p} \boldsymbol{u} d V$ can be recast in a matrix form $\boldsymbol{M}$. Similarly, the second term, $\sum_f A_f \cdot(\boldsymbol{u} \otimes \boldsymbol{u})_f$ (convective term) produces the matrix $\boldsymbol{C}$, and the diffusive term $\sum_f A_f \cdot(2\nu \boldsymbol{\nabla}^s \boldsymbol{u})_f$ generates the matrix $\boldsymbol{D}$. The pressure term $\sum_f A_f  \cdot \boldsymbol{p}_f$ gives rise to the matrix $\boldsymbol{B}$. The coefficients from \autoref{eq:continuity-FV} assemble the matrix $\boldsymbol{P}$. For each finite volume, the interpolation coefficients obtained during the discretization process are used to form an algebraic form of the equation as follows:

\begin{equation}\label{eq:fom-algebric}
\begin{aligned}
\boldsymbol{M \dot { u }}+\boldsymbol{C}(\boldsymbol{u}) \boldsymbol{u}+\nu \boldsymbol{D} \boldsymbol{u}+\boldsymbol{B} \boldsymbol{p} & = 0, \\
\boldsymbol{P} \boldsymbol{u} & =0.
\end{aligned}
\end{equation}

Now, the main idea behind the reduced-order model is to find a spatial basis $\phi(x)$ that spans over a subspace $\mathcal{S}$ for a state vector (i.e., velocity and pressure, etc.) $\boldsymbol{\zeta}(\boldsymbol{x}, \boldsymbol{\mu}, t) \approx \boldsymbol{\zeta}_s=\sum_{i=1}^{N_{\zeta}^s} \alpha_i(t, \boldsymbol{\mu}) \boldsymbol{\phi}_i(\boldsymbol{x})$ where, the reduced field is $\boldsymbol{\zeta}_s$. $\alpha_i(t, \boldsymbol{\nu})$ are the temporal coefficients which depends on the parameter vector $\nu$. $N_{\zeta}^s$ is the cardinality of the reduced space for the variable $\zeta$. Several techniques are available in the literature for constructing the reduced basis, including the Proper Orthogonal Decomposition (POD), Proper Generalised Decomposition (PGD), and the Reduced Basis (RB) method based on a greedy sampling strategy. Details of the different techniques can be found in \cite{rozza2008reduced,chinesta2015model,kalashnikova2010stability,dumon2011proper}. In this work, we employ the POD method applied to the snapshots that account for both time and parameter dependence. Now, we briefly recall the standard POD procedure.

\subsubsection{Proper Orthogonal Decomposition and Galerkin Projection}

In the POD approach, First, the state vector solutions are gathered using a high-fidelity solver. For example, in the current case scenario, the discrete time instances are denoted as $t^j \in \{t^1, \ldots, t^{N_t}\} \subset [0, T]$, associated with a given parameter value $\nu^j$ (where $j$ varying from $0$ to $m$). The total number of snapshots is thus $ N_S = m \times N_t$. The snapshot matrices $\mathcal{Z}_{\boldsymbol{\zeta}}$ and $\mathcal{Z}_{\boldsymbol{p}}$ for velocity and pressure, respectively, are formed by collecting $N_s$ full-order solutions, where $N_u^h$ and $N_p^h$ denote the degrees of freedom associated with the velocity and pressure fields.

\begin{equation}
\begin{gathered}
\mathcal{Z}_{\boldsymbol{u}} = \left[\boldsymbol{u}(\nu^1, t^1), \ldots, \boldsymbol{u}(\nu^m, t^{N_t})\right] \in \mathbb{R}^{N_u^h \times N_s}, \\
\mathcal{Z}_{\boldsymbol{p}} = \left[\boldsymbol{p}(\nu^1, t^1), \ldots, \boldsymbol{p}(\nu^m, t^{N_t})\right] \in \mathbb{R}^{N_p^h \times N_s}.
\end{gathered}
\end{equation}

Once the snapshot matrices are assembled, the reduced-order model (ROM) can be constructed to approximate solutions for new parameters and time instances. Given a set of vector-valued velocity functions $\boldsymbol{u}(t): Z \rightarrow \mathbb{R}^d$, represented by a collection of $N_s$ realizations $\boldsymbol{u}_1, \ldots, \boldsymbol{u}_{N_s}$ (also referred to as \emph{snapshots}). Now, the Proper Orthogonal Decomposition (POD) aims to identify, for a prescribed reduced dimension $N_{\text{POD}} \leq N_s$, an optimal set of basis functions $\{\varphi_1, \ldots, \varphi_{N_{\text{POD}}}\}$ and corresponding coefficients $\{a_i^k\}$ that minimize the projection error:

\begin{equation} \label{POD-monimization}
E_{N_{\text{POD}}} = \sum_{i=1}^{N_s} \left\| \boldsymbol{u}_i - \sum_{k=1}^{N_{\text{POD}}} a_i^k \varphi_k \right\|,
\end{equation}

subjected to the orthonormality condition:

\begin{equation}
\langle \varphi_i, \varphi_j \rangle_{L^2(V_p)} = \delta_{ij}, \quad \forall i,j=1,\ldots, N_s.
\end{equation}

The minimization problem stated in \autoref{POD-monimization} is equivalent to solving the following eigenvalue problem:

\begin{equation}
\boldsymbol{\mathcal{C}}^u \boldsymbol{G}^u = \boldsymbol{G}^u \boldsymbol{\lambda}^u,
\end{equation}

where $\boldsymbol{\mathcal{C}}^u$ is the correlation matrix obtained from $\mathcal{Z}_{u}$ and defined by:

\begin{equation}
\mathcal{C}^u_{ij} = \langle \boldsymbol{u}_i, \boldsymbol{u}_j \rangle_{L^2(V_p)}, \quad i,j=1,\ldots,N_s.
\end{equation}

$\boldsymbol{G}^{\boldsymbol{u}}$ is a matrix of eigenvectors and $\boldsymbol{\lambda}^{\boldsymbol{u}}$ is a matrix of eigenvalues.The POD basis functions are then constructed as:

\begin{equation}
\varphi_i = \frac{1}{N_s \lambda_i^u} \sum_{j=1}^{N_s} \boldsymbol{u}_j G_{ij}^u.
\end{equation}

Following this procedure, we compute separate POD spaces for velocity and pressure:

\begin{equation}
\begin{aligned}
L_u &= \left[\varphi_1, \ldots, \varphi_{N_u^r} \right] \in \mathbb{R}^{N_u^h \times N_u^r}, \\
L_p &= \left[\chi_1, \ldots, \chi_{N_p^r} \right] \in \mathbb{R}^{N_p^h \times N_p^r},
\end{aligned}
\end{equation}

where $N_u^r, N_p^r < N_s$ are selected based on the decay of the eigenvalues $\boldsymbol{\lambda}^u$ and $\boldsymbol{\lambda}^p$. Once the POD basis is defined, the velocity and pressure fields are approximated as:

\begin{equation}
\begin{aligned}
\boldsymbol{u}^r(\boldsymbol{x}, t, \mu) &\approx \sum_{i=1}^{N_u^r} a_i(t, \mu) \varphi_i(\boldsymbol{x}), \\
\boldsymbol{p}^r(\boldsymbol{x}, t, \mu) &\approx \sum_{i=1}^{N_p^r} b_i(t, \mu) \chi_i(\boldsymbol{x}),
\end{aligned}
\end{equation}

where, the coefficients $a_i$ and $b_i$ depend on time and the parameter $\mu$, while the basis functions depend only on the spatial variable. Applying a Galerkin projection of the full-order governing equations \autoref{eq:NS-FV1} and \autoref{eq:NS-FV2} onto the POD space leads to the reduced-order system of ordinary differential equations. 

\begin{equation}\label{eq:red_system1}
\begin{aligned}
\boldsymbol{M}_r \dot{a} - \nu \boldsymbol{D}_r a + \boldsymbol{C}_r(a) a + \boldsymbol{B}_r b &= 0, \\
\end{aligned}
\end{equation}

and,
\begin{equation}\label{eq:red_system2}
\begin{aligned}
\boldsymbol{P}_r a &= 0,
\end{aligned}
\end{equation}

where the reduced matrices are defined as:

\begin{equation}\label{Matrices}
\begin{aligned}
\boldsymbol{M}_{r_{ij}} &= \langle \varphi_i, \varphi_j \rangle_{L^2(V_p)}, \\
\boldsymbol{A}_{r_{ij}} &= \langle \varphi_i, \nabla \cdot 2 \nabla^s \varphi_j \rangle_{L^2(V_p)}, \\
\boldsymbol{B}_{r_{ij}} &= \langle \varphi_i, \nabla \chi_j \rangle_{L^2(V_p)}, \\
\boldsymbol{P}_{r_{ij}} &= \langle \chi_i, \nabla \cdot \varphi_j \rangle_{L^2(V_p)}.
\end{aligned}
\end{equation}

Once we define the reduced basis spaces $L_u$ and $L_p$ through the basis functions $\varphi_i$ and $\chi_i$, all the reduced matrices of equation \autoref{Matrices} can be computed from the offline stage without difficulties, except for the reduced order matrix $\boldsymbol{C}_r(a)$, originating from the nonlinear convective term. The strategy employed here involves storing a third-order tensor $\mathbf{C}_r$ whose entries are given by: 

\begin{equation}
{\mathrm{C}}_{r_{ijk}} = \langle \varphi_i, \nabla \cdot (\varphi_j \otimes \varphi_k) \rangle_{L^2(V_p)}.
\end{equation}

Now, during the online stage, at each fixed point iteration of the solution procedure, each entry of the contribution to the reduced residual given by the convective term $\mathcal{R}_c^r=C_{\boldsymbol{r}}({a}) {a},$ can be computed with:

\begin{equation}
    \boldsymbol{\mathcal { R }}_{c_i}^r=\left(\boldsymbol{C}_{\boldsymbol{r}}({a}) {a}\right)_i={a}^T \mathbf{C}_{r_i . .} {a} .
\end{equation}

In addition, different strategies can be employed to enrich the velocity space such that the inf-sup condition is met, such as by suprimizer enrichment, which is explained in \cite{STABILE2018273,ballarin2015supremizer}.

\begin{remark}
It is important to remark that despite the reduction of computational time, such Galerkin projection-based ROM may not be applied to quicker prediction of solution field at a new parameter or time-instants (such as in real-time), which is possible using physics-informed neural network. Therefore, in the following section, we will demonstrate how to couple an external reduced-order CFD solver with neural network.
\end{remark}

\subsubsection{Coupling PINN with Reduced Order System} \label{PINN-ROM}

This section will discuss how to couple the reduced order system as described in \autoref{eq:red_system1} and \autoref{eq:red_system2} as a physics-based loss term in the neural network. \autoref{eq:red_system1} can be formulated in the following manner: 

\begin{equation}\label{eq:matrices1}
\begin{aligned}
\dot{a} = \mathcal{X}(a,b), 
\end{aligned}
\end{equation}
\\
where, the spatial residual $\mathcal{X}(a,b)$ can be written as the following way: 

\begin{equation}\label{eq:spatial_res}
\begin{aligned}
\mathcal{X}(a,b) = \boldsymbol{M}_r^{-1}(\nu \boldsymbol{D}_r a-\boldsymbol{C}_r(a) a-\boldsymbol{B}_r b). 
\end{aligned}
\end{equation}
\\
The physics-based residual arising from the reduced order system of equations \autoref{eq:red_system1} and \autoref{eq:red_system2} can now be written as follows: 

\begin{equation}\label{eq:matrices2}
\begin{aligned}
\mathbf{R}_{\text{red1}} &= \dot{a} - \mathcal{X}(a,b), \\
\mathbf{R}_{\text{red2}} &= \boldsymbol{P}_r a. 
\end{aligned}
\end{equation}
\\
Similar to the full-order model presented in \autoref{Burgers_FOM}, the PINN utilises the residual terms $\mathbf{R}_{\text{red1}}$ and $\mathbf{R}_{\text{red2}}$ corresponding to \autoref{eq:red_system1} and \autoref{eq:red_system2}, respectively, as physics-based loss terms which are imported from an external CFD solver. It is assumed that the neural network can not include the algebraic forms of these residuals within its computational graph. As mentioned earlier, in the present reduced-order model for viscous incompressible flow, the input vector $\mathbf{z}$ represents an unsteady problem defined over $N_t$ time instances, i.e., $\mathbf{t} = [ t^1, t^2, \ldots, t^{N_t} ]$, and $m$ parameter values of the viscosity $\nu$, i.e., $\boldsymbol{\nu} = [ \nu^1, \nu^2, \ldots, \nu^m ]$. The output $\mathbf{Q}$ consists of the reduced-order variables $a$ and $b$, i.e., $\mathbf{Q} = [a, b]$. Therefore, the dimensions of the input and output are $\mathbf{z} \in \mathbb{R}^{mN_t \times 2}$ and $\mathbf{Q} \in \mathbb{R}^{mN_t \times (N_u^r+N_p^r)}$, respectively. The $i^{\text{th}}$ row of the input $[t^{i}, \nu^{i}]$ corresponds to the output $\text{Q}_{\text{pred},i}$ and associated residuals obtained from external reduced order solver are $\mathbf{R}_{\text{red1},i}$ and $\mathbf{R}_{\text{red2},i}$. The modified loss term $\text{L}_{\text{Dis}}$, introduced earlier in \autoref{eq:LDIS}, can now be formulated for the present reduced-order PINN system as follows, where $\mathbf{Q}_{\text{pred}}$ denotes the neural network prediction of $\mathbf{Q}$:

\begin{equation}\label{eq:LDIS_ROM1}
\begin{aligned}
\mathrm{L}_{\text{Dis,Mom}} = \left(\frac{1}{N_{\text{eqn}}} \sum_1^{N_{\text{eqn}}} \left[ 2 \mathbf{R}_{\text{red1},i} \frac{\partial \mathbf{R}_{\text{red1},i}}{\partial \mathbf{Q}_{\text{pred}}} \right]_{\text{detached}} \mathbf{Q}_{\text{pred}} \right),
\end{aligned}
\end{equation}

\begin{equation}\label{eq:LDIS_ROM2}
\begin{aligned}
\mathrm{L}_{\text{Dis,Pressure}} =  \left(\frac{1}{N_{\text{eqn}}} \sum_1^{N_{\text{eqn}}} \left[ 2 \mathbf{R}_{\text{red2},i} \frac{\partial \mathbf{R}_{\text{red2},i}}{\partial \mathbf{Q}_{\text{pred}}} \right]_{\text{detached}} \mathbf{Q}_{\text{pred}} \right).
\end{aligned}
\end{equation}

In this work, the time derivative term $\dot{a}$ in $\mathbf{R}_{\text{red1}}$ is computed using automatic differentiation. However, the spatial residual terms $\mathcal{X}(a, b)$ are imported from the external solver. The rationale behind such a choice is that when finer time steps are used, the input size can increase if all time instants are retained to match the time discretization of the external reduced-order solver. Now, we try to derive the formulation of $\mathrm{L}_{\text{Dis}}$ when the temporal derivative ($\dot{a}$) is computed using automatic differentiation in \autoref{eq:matrices1} and spatial residual $\mathcal{X}(a,b)$ is obtained from the external reduced order system. Now \autoref{eq:LDIS_ROM1} and \autoref{eq:LDIS_ROM2} can be further modified. We start with, 

\begin{equation}\label{eq:leqn_mom}
\begin{aligned}
\mathrm{L}_{\text{eqn, mom}}=\frac{1}{N_{\text{eqn}}} \sum_{i=1}^{N_{\text{eqn}}} {\mathbf{R}_\text{red1,i}}^2,
\end{aligned}
\end{equation}

\begin{equation}\label{eq:leqn_pressure}
\begin{aligned}
\mathrm{L}_{\text{eqn, pressure}}=\frac{1}{N_{\text{eqn}}} \sum_{i=1}^{N_{\text{eqn}}} {\mathbf{R}_\text{red2,i}}^2.
\end{aligned}
\end{equation}

Now, the derivative of $\mathrm{L}_{\text{eqn,mom}}$ concerning the hyper-parameters can be computed as follows:

\begin{equation}
\label{eq:grad_momentum_loss}
\begin{aligned}
\frac{\partial \mathrm{L}_{\text{eqn,mom}}}{\partial \theta}
&= \frac{1}{N_{\text{eqn}}} 
\sum_{i=1}^{N_{\text{eqn}}} 
  2 \, \mathbf{R}_{\text{red1},i} 
  \frac{\partial \mathbf{R}_{\text{red1},i}}{\partial \theta}, \\
&= \frac{1}{N_{\text{eqn}}} 
\sum_{i=1}^{N_{\text{eqn}}} 
  2 \, \mathbf{R}_{\text{red1},i} 
  \frac{\partial \left(\dot{a} - \mathcal{X}(a,b)\right)}{\partial \theta},\\
&= \frac{1}{N_{\text{eqn}}} 
\sum_{i=1}^{N_{\text{eqn}}} (
  2 \, \mathbf{R}_{\text{red1},i} 
  \frac{\partial \left(\dot{a}\right)}{\partial \theta} - 2 \, \mathbf{R}_{\text{red1},i} 
  \frac{\partial \mathcal{X}(a,b)}{\partial \theta}).
\end{aligned}
\end{equation}

Since $\dot{a}$ is computed using automatic differentiation and hence included in the computational graph, \autoref{eq:grad_momentum_loss} can be further written as follows:

\begin{equation}\label{eq:LDIS_Burgers}
\begin{aligned}
\mathrm{L}_{\text{Dis, Mom}}=\frac{1}{N_{\text {eqn }}} \sum_1^{N_{\text {eqn }}}\left(\left[2 \textbf{R}_i\right]_{\operatorname{detached}} \textbf{R}_i\right) - (\left[2 \textbf{R}_i \frac{\partial \mathcal{X}(a,b)}{\partial \mathbf{Q}_{\text{pred}}}\right]_{\operatorname{detached}} \mathbf{Q}_{\text{pred}}).
\end{aligned}
\end{equation}

Similarly, the derivative of $\mathrm{L}_{\text{eqn, pressure}}$ with respect to the neural network parameter, $\theta$ can be written as $\frac{\partial \mathrm{L}_{\text{eqn, pressure}}}{\partial \theta}
=\frac{1}{N_{\text{eqn}}} 
\sum_{i=1}^{N_{\text{eqn}}} 
  2 \, \mathbf{R}_{\text{red2},i} 
  \frac{\partial \mathbf{R}_{\text{red2},i}}{\partial \theta}$ while the corresponding $\mathrm{L}_{\text{Dis, pressure}}$ can be written as \autoref{eq:LDIS_ROM2}. Therefore, the total loss term can be computed as follows: 

  \begin{equation}\label{eq:ltot2}
\begin{aligned}
\text{L} = {\lambda_{1}}\text{L}_{\text{Dis, Mom}} + {\lambda_{2}}\text{L}_{\text{Dis, pressure}} + {\lambda_{3}}\text{L}_{\text{Data}},
\end{aligned}
\end{equation}

where, ${\lambda_{1}}$, ${\lambda_{2}}$ and ${\lambda_{3}}$ are the constatnts associated with different loss terms.

\begin{figure}[ht]
\centering
\begin{subfigure}[b]{1\textwidth}
\centering
% \hspace{-3cm}
{\label{fig:burgers3}\includegraphics[width=0.8\linewidth]{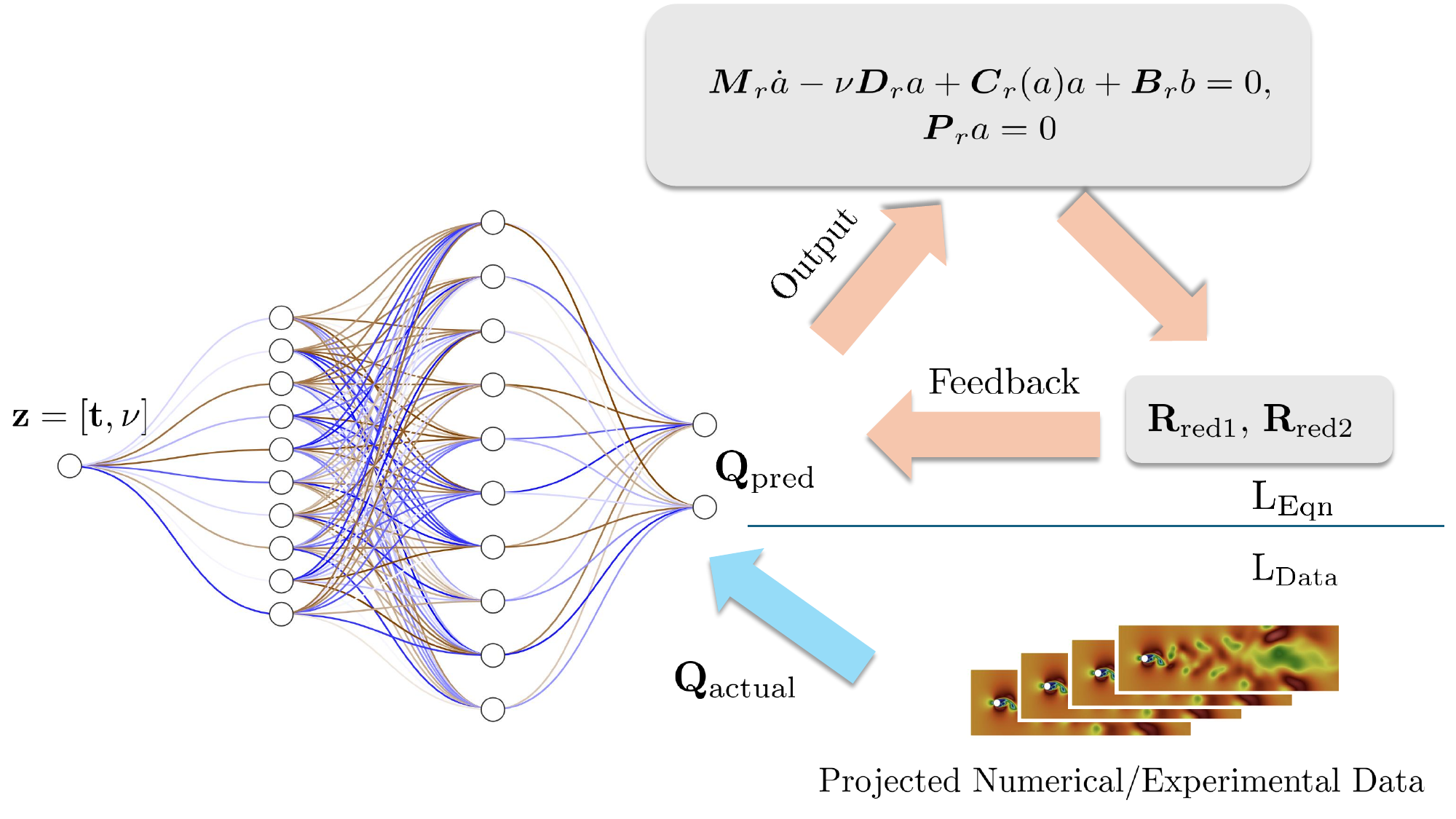}}
\end{subfigure}
\caption{Schematic of the discretized physics-based PINN and coupling procedure with any external reduced order formulation where $\text{L}_\text{eqn}$ and $\text{L}_\text{Data}$ are physics-based and data-driven loss terms, respectively. }
\label{fig:ANN-DisPINN}
\end{figure}

\SetKwComment{Comment}{/* }{ */}

\begin{algorithm}[hbt!] 
\caption{DisPINN pseudo-code for coupling external reduced order solver}\label{alg:algo}
\KwData{$\mathbf{Q}_\text{actual}, \textbf{z}, N_\text{Data}, N_\text{eqn}$}

$[\boldsymbol{W},\boldsymbol{c}] \gets \textbf{INIT} ([\boldsymbol{W},\boldsymbol{c}])$ \Comment*[r]{Initialize Weight and bias}
\While{$\text{epoch} \leq \text{Total Epoch No.}$}{
  
  $\mathbf{Q}_\text{pred} \gets \text{NN}(\boldsymbol{W},\boldsymbol{c},f, \textbf{z})$ \Comment*[r]{prediction from chosen network $\text{NN}$}
  $\text{Ext. Solver} \gets \text{Q}_\text{pred,i}$ \Comment*[r]{pass to external solver}
  $ \mathbf{R}_{\text{red1},i}, \mathbf{R}_{\text{red2},i} \gets \text{Ext. Solver}$  \Comment*[r]{pass to PINN solver from external solver}
  $\mathrm{L}_{\text{eqn}} \gets (\frac{1}{N_{\text {eqn }}}\sum_{i=1}^{N_{\text{eqn}}} {\mathbf{R}_{\text{red1},i}}^2 + \frac{1}{N_{\text {eqn }}}\sum_{i=1}^{N_{\text{eqn}}} {\mathbf{R}_{\text{red2},i}}^2)$ \Comment*[r]{compute physics-based loss term}
  $\mathrm{L}_{\text {Data }} \gets \frac{1}{N_{\text {Data }}}\sum_{i=1}^{N_{\text {Data }}}\left(\text{Q}_{\text {pred,i }}- \text{Q}_{\text {actual,i }}\right)^2$ \Comment*[r]{compute data-based loss term}

  \eIf{$\text{rem}(epoch, k_{int}) = 0$}{  
    $\frac{\partial \mathbf{R}_{\text{red1},i}}{\partial \mathbf{Q}_{\text{pred}}}, \frac{\partial \mathbf{R}_{\text{red2},i}}{\partial \mathbf{Q}_{\text{pred}}} \gets \text{Ext. Solver}$ \Comment*[r]{compute Jacobian from external solver}
    $J_\text{new1,i} \gets \frac{\partial \mathbf{R}_{\text{red1},i}}{\partial \mathbf{Q}_{\text{pred}}}$\; 
 
    $J_\text{new2,i} \gets \frac{\partial \mathbf{R}_{\text{red2},i}}{\partial \mathbf{Q}_{\text{pred}}}$\; 
    
    $\mathrm{L}_{\text{Dis}} \gets (\frac{1}{N_{\text {eqn }}} \sum_1^{N_{\text {eqn }}}2 \mathbf{R}_{\text{red1},i} J_\text{new,i} \mathbf{Q}_\text{pred}) + (\frac{1}{N_{\text {eqn }}} \sum_1^{N_{\text {eqn }}}2 \mathbf{R}_{\text{red2},i} J_\text{new,i} \mathbf{Q}_\text{pred})$  \Comment*[r]{compute additional physics based loss term for backpropagation}
    }
    {$J_{1,i} \gets J_\text{new1,i}$ \;
    $J_{2,i} \gets J_\text{new2,i}$ \;
    $\mathrm{L}_{\text{Dis}} \gets (\frac{1}{N_{\text {eqn }}} \sum_1^{N_{\text {eqn }}}2 \mathbf{R}_{\text{red1},i} J_{1,i} \mathbf{Q}_\text{pred} + \frac{1}{N_{\text {eqn }}} \sum_1^{N_{\text {eqn }}}2 \mathbf{R}_{\text{red2},i} J_{2,i} \mathbf{Q}_\text{pred})$ \;
  }

  $\mathrm{L} \gets (\mathrm{L}_{\text{Dis}} + \mathrm{L}_{\text {Data }})$ \Comment*[r]{compute total loss term for backpropagation}

  $\Delta {\boldsymbol{W}} \leftarrow-\alpha \mathrm{G}_{\mathrm{ADAM}}\left(\nabla_{{\boldsymbol{W}}} L\right), \Delta {\boldsymbol{c}} \leftarrow-\alpha \mathrm{G}_{\mathrm{ADAM}}\left(\nabla_{{\boldsymbol{b}}} L\right)$ \; \Comment*[r]{compute weight matrices and bias vector update}
  ${W} \leftarrow {W}+\Delta {W}, {\boldsymbol{c}} \leftarrow {\boldsymbol{c}}+\Delta {\boldsymbol{c}}$ \;
}
\end{algorithm}

\section{\bf{Numerical Experiments}}\label{sec:PINN}

This section demonstrates the numerical experiments on the test cases mentioned in \autoref{sec:testcase}. The first test case is a nonlinear transport equation, where the finite-volume formulation of the full-order governing equation is coupled with the PINN solver. Here, time instants are considered the only input in the neural network. As a second test case, viscous incompressible flow past a cylinder is considered. A reduced-order framework is first derived from the finite-volume-based formulation of the Navier-Stokes equation, followed by coupling with the physics-informed neural network. Here, the input consists of both the physical time and kinematic viscosity. We note that even when only the residuals from the governing equations in the external solver are available, or when limited knowledge of the discretized system is available to the PINN solver, the proposed approach still enables effective neural network training through the introduction of a modified physics-based loss term.

\subsection{Non-linear Transport Equation Model (NTE)}\label{sec:burgers}

\begin{figure}[ht]
\centering
\begin{subfigure}[b]{1\textwidth}
\centering
% \hspace{-3cm}
{\label{fig:burgers1}\includegraphics[width=0.5\linewidth]{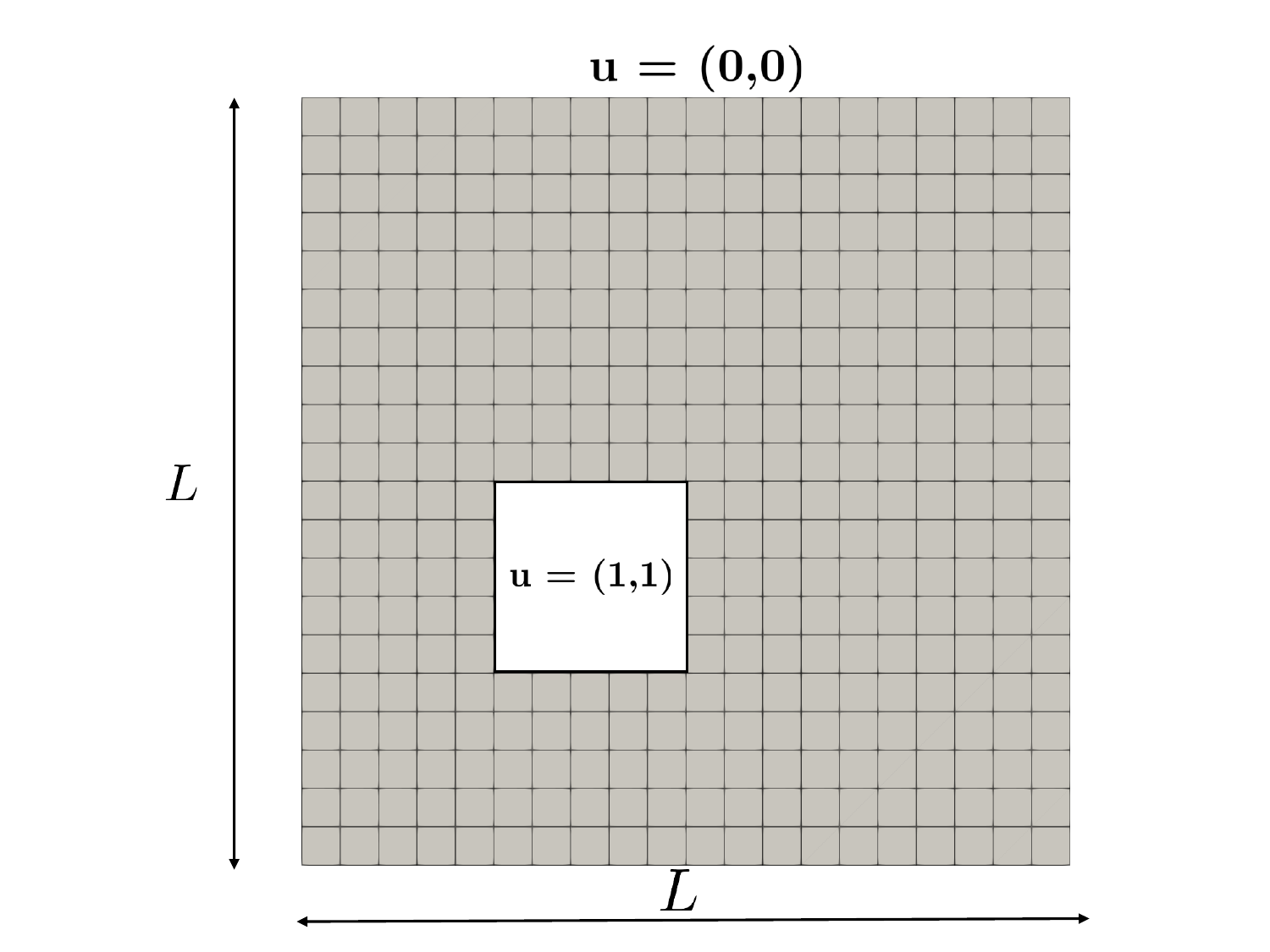}}
\end{subfigure}
\caption{Mesh topology, boundary and initial values.}
\label{fig:Burgers-Mesh}
\end{figure}

The benchmark full-order simulation based on the finite volume-based formulation of the nonlinear transport equation, as indicated in \autoref{eq:Burgers}, is carried out on a two-dimensional square domain of length ${L} = 1$ m. We consider a Eulerian frame on a space-time domain with $(\mathrm{x}, t) \in[0,1]^d \times[0,0.35]$ with $\text{d} = 2$. A uniform square domain denoted by $\mathbf{x} \in[0.25,0.5]^2$ with a uniformly distributed $\mathbf{u}$ value of $1$ (in both $x$ and $y$ directions) is considered as an initial condition. The boundary values at $\mathbf{x} \in \partial[0,1]^2$ is taken as $0$. In this test case, the temporal derivative is discretized using a first-order implicit Euler scheme, the convective term uses a second-order linear upwind scheme, and the diffusion term is discretized using a central differencing scheme with non-orthogonal correction (linear corrected). The computational domain is divided into $21\times21$ equal mesh size. The time step is kept constant and equal to $\text{t} = 0.001$ s, and the simulation is run till $\text{T} = 0.35$ s. The time-steps are sufficiently small to meet the CFL criteria at all parts of the computational domain. The value of the diffusion coefficient $v$, mentioned in \autoref{eq:Burgers}, is taken as $0.01$. The mesh topology and boundary conditions are shown in \autoref{fig:Burgers-Mesh}. This benchmark CFD-based solution, collected at sparse temporal locations, will be utilized to generate the data-driven loss term in neural network training and to compute the absolute error distribution of the predicted solution field from different data-driven and physics-based networks. However, the sparse solution field can also be obtained from the experimental measurements.

Our current objective is to couple the finite volume-based formulation of Eqn. \ref{eq:Burgers} in the external solver with the physics-informed neural network framework to obtain a functional relationship between the solution field $\bf{u}$ as output and time instants $\mathbf{t}$ as input. Once the network is trained, the solution can be obtained in real-time at any time instant. Since the PINN and CFD solver belong to a separate environment, the discretized form of the governing equation \autoref{eq_dis2} can not be included in the computational graph of the PINN solver. Deepinder et al.\cite{aulakh2022generalized} proposed the introduction of the finite-volume-based discrete loss term in the PINN framework. However, the proposed methods were limited to only linear governing equations.  In the present scenario,  we have assumed that the external solver can export the numerical residual obtained from the finite volume discretization of the governing \autoref{eq_dis2}, the system matrix $\mathbf{A}$ and vector $\mathbf{b}$ described in \autoref{linear_Mat} to the PINN solver.  

\begin{figure}[ht!]
\centering
\begin{subfigure}[b]{1\textwidth}
\centering
% \hspace{-3cm}
\subfloat[Benchmark]{\label{fig:Benchmark0}\includegraphics[width=1.0\linewidth]{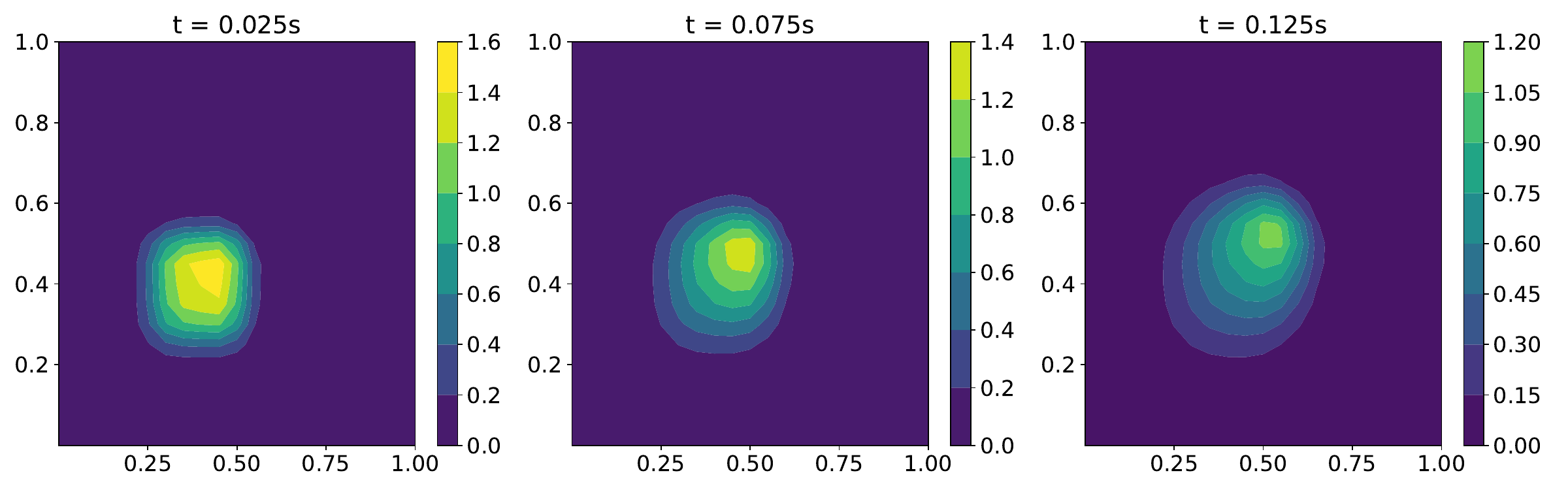}}
\end{subfigure}
\begin{subfigure}[b]{1\textwidth}
%\centering
% \hspace{-3cm}
\subfloat[Prediction from ANN-DisPINN]{\label{fig:Prediction0}\includegraphics[width=1.0\linewidth]{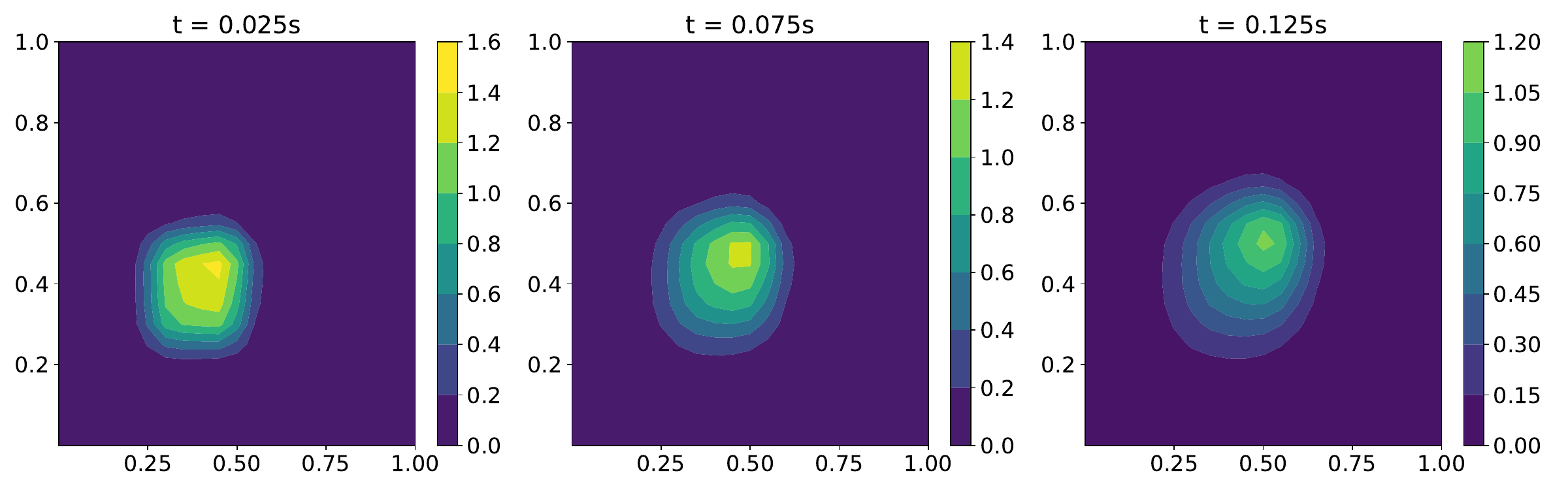}}
\end{subfigure}

\begin{subfigure}[b]{1\textwidth}
%\centering
% \hspace{-3cm}
\subfloat[Absolute Error]{\label{fig:Error0}\includegraphics[width=1.0\linewidth]{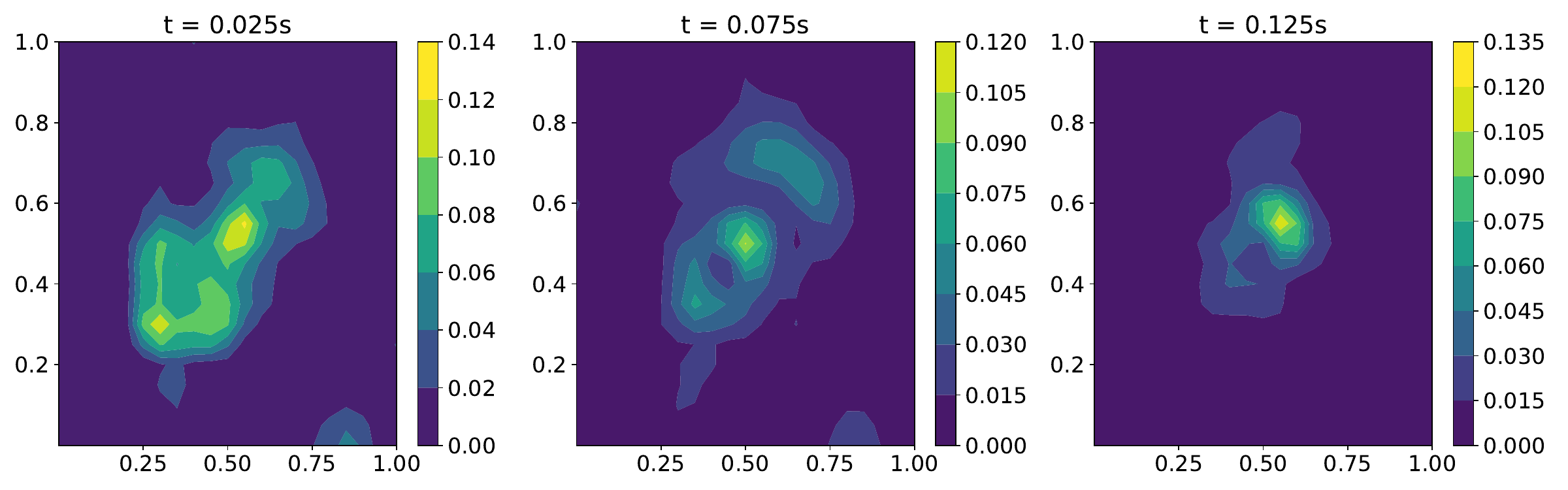}}
\end{subfigure}

\medskip
\caption{Comparison of the benchmark CFD solution and predicted results obtained from the proposed DisPINN. The solution snapshot at the initial time (t = 0) is only considered for the data-driven loss term. The physics-based residual, $\mathbf{R}$ obtained from \autoref{eq:Burgers_discretized}, coefficient matrix $\mathbf{A}$ and $\mathbf{b}$ are imported from the external CFD solver at a time interval of $0.01$ s. The figure shows (a) Benchmark CFD-based Results, (b) Predicted Results from DisPINN, and (c) Absolute Error between the benchmark and predicted results.}
\label{fig:comaparison_int}
\end{figure}

In the neural network, $4$ layers with $124, 100, 80$ and $64$ neurons are considered, respectively. $\text{softplus}$ activation function is considered and the learning rate is $0.006$. The total number of the $\text{epoch}$ (number of iterations considered in the optimization) is $3999$. The algorithms proposed in \autoref{alg:algo-burgers} are applied to obtain the modified physics-based loss term $\text{L}_{\text{Dis}}$, which will be utilized in the current DisPINN training, despite limited access to the discretized form of the governing equation in the CFD solver. Now, in the neural network platform, we consider numerical data from previously generated snapshots of the benchmark solution for the data-driven loss term $\text{L}_\text{data}$ at several $N_\text{data}$ locations. First, we have considered the time instants $\text{t} = 0$ to $\text{t} = 0.2$ s for training the network.  \autoref{fig:comaparison_int} compares the prediction from our proposed DisPINN with the benchmark solution. The solution snapshot at the initial time (t = 0) is only considered for the data-driven loss term. On the contrary, the physics-based residual, $\mathbf{R}$ obtained from \autoref{eq:Burgers_discretized}, coefficient matrix $\mathbf{A}$ and $\mathbf{b}$ are imported from the external CFD solver at a time interval of $0.01$ s for the computation of the physics-based loss term. The predicted results match significantly well with the benchmark results at time instants of $\text{t} = 0.025$, $\text{t} = 0.075$ and $\text{t} = 0.125$, as shown in \autoref{fig:Benchmark0} and \autoref{fig:Prediction0}. The maximum absolute error incurred at those time steps is $0.13$, $0.11$ and $0.12$ respectively, as shown in \autoref{fig:Error0}. The absolute error diminishes as the numerical simulation data in the data-driven loss term increases.

\begin{figure}[ht]
\centering
\begin{subfigure}[b]{1\textwidth}
\centering
% \hspace{-3cm}
\subfloat[Data+Physics from External CFD Solver]{\label{fig:Burgers_Physics_5}\includegraphics[width=1.0\linewidth]{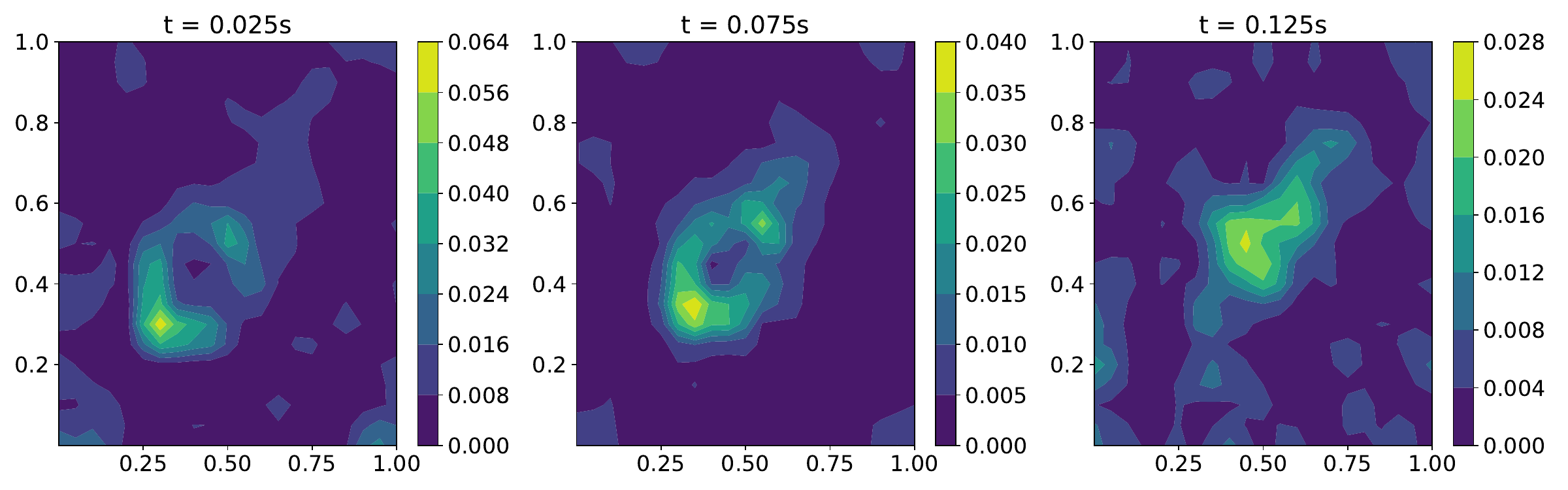}}
\end{subfigure}
\begin{subfigure}[b]{1\textwidth}
%\centering
% \hspace{-3cm}
\subfloat[Only Data]{\label{fig:Burgers_Data_5}\includegraphics[width=1.0\linewidth]{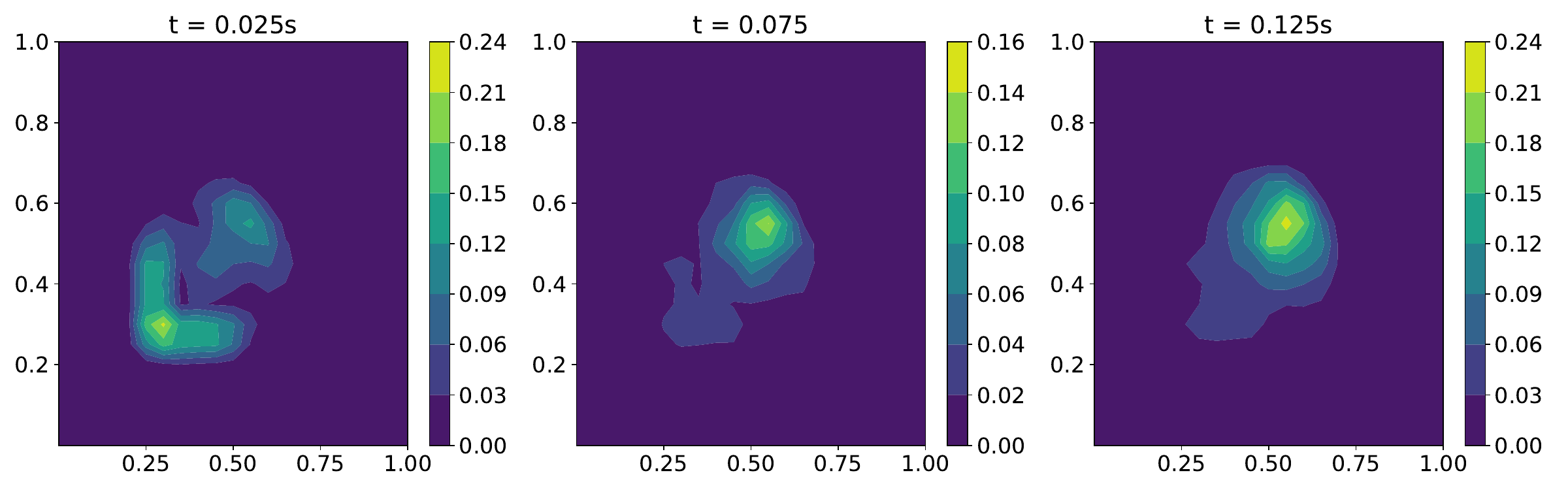}}
\end{subfigure}
\medskip
\caption{Absolute error with respect to the benchmark data when (a) both data-driven and physics-driven loss term is considered, (b) only the data-driven loss term is considered. For the numerical experiments, $3$ intermediate time-instants where numerical data ($0$, $0.05$, $0.1$ s) are considered for the data-driven loss term.  The physics-based residual is imported from the external CFD solver at a time interval of $0.01$ s between time instants of $0$ s to $0.2$ s.}
\label{fig:Comparison_Burgers}
\end{figure}

The prediction error is improved when more numerical data is added to the physics-informed neural network (as a data-driven loss term) from the previously computed snapshots. We have considered $3$ intermediate time-instants where numerical data ($0$, $0.05$, $0.1$ s) is considered for the data-driven loss term in the physics-informed neural network and we want to carry out the same numerical experiments, reconstruction of the solution field at time instants $\text{t} = 0.025$ s, $\text{t} = 0.075$ s and $\text{t} = 0.125$ s. The physics-based loss terms are considered at every $0.01$ s (starting from $0s$ to $0.2s$) for all the spatial grid locations. \autoref{fig:Comparison_Burgers} presents a comparison of the absolute error distribution between the benchmark solution and the predictions obtained from two neural network models: one trained using only a data-driven loss term, and the other trained with a combination of physics-based and data-driven loss terms. The maximum absolute error associated with the solution fields at time instants $\text{t} = 0.025$, $\text{t} = 0.075$ and $\text{t} = 0.125$ is $0.063$, $0.04$ and $0.03$ when both the data-driven and physics-driven loss terms are considered as shown in \autoref{fig:Burgers_Physics_5}. It is important to notice that the maximum absolute error is drastically reduced as compared to the case where only the snapshot at the initial condition ($t = 0$) is considered for the data-driven loss terms alongside the physics-based loss terms. \autoref{fig:Burgers_Data_5} shows the absolute error distribution when only data-driven residual is considered without the physics information, and the maximum absolute errors incurred in the solution fields are $0.21$, $0.14$ and $0.22$, respectively corresponding to the above-mentioned time instants.

\begin{figure}[ht]
\centering
    \begin{subfigure}[b]{0.49\textwidth}
        \centering
        \includegraphics[width=\linewidth]{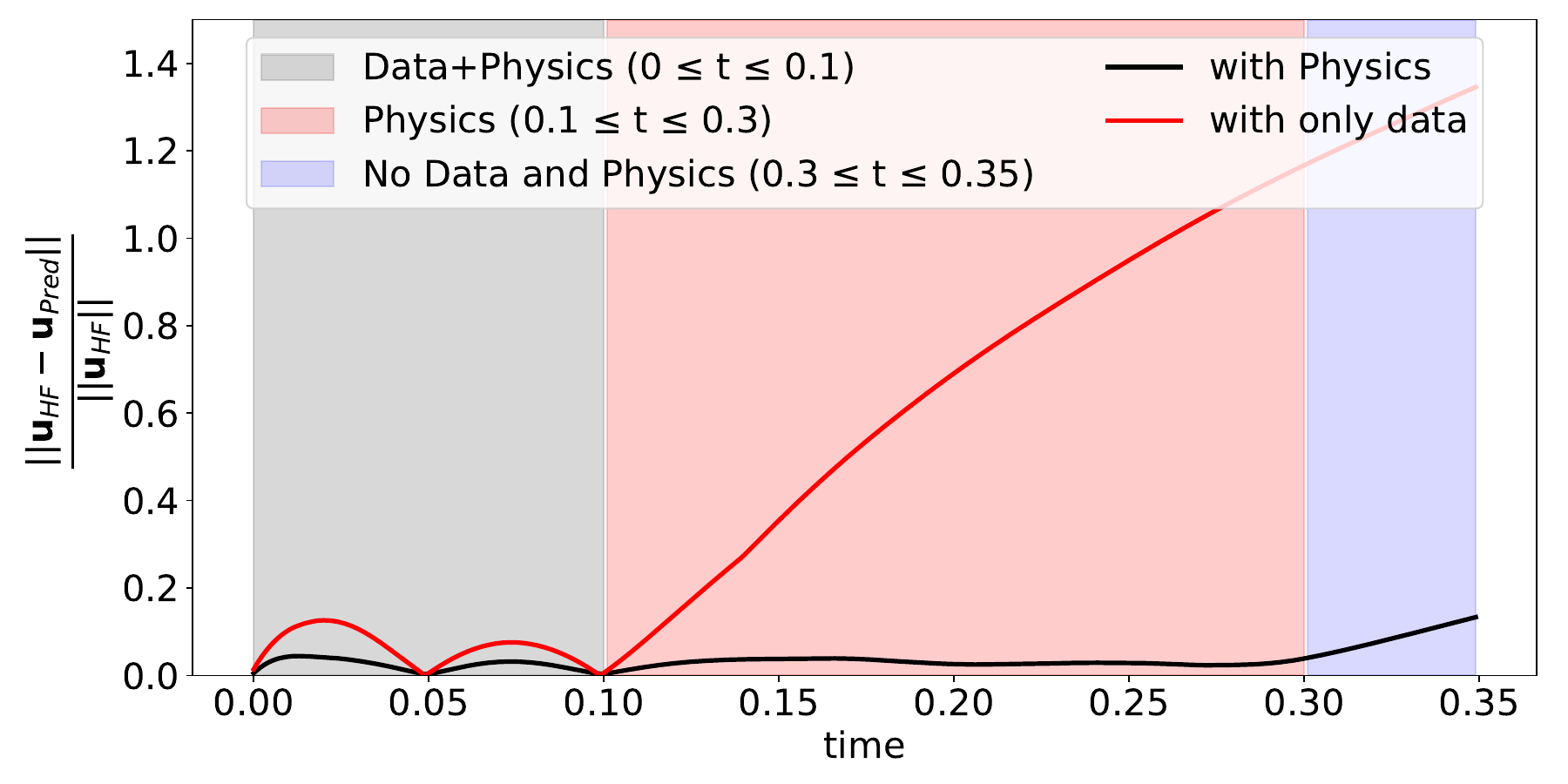}
        \caption{Loss with correction term}
        \label{fig:Error_Burgers_Full}
    \end{subfigure}
    \hfill
    \begin{subfigure}[b]{0.49\textwidth}
        \centering
        \includegraphics[width=\linewidth]{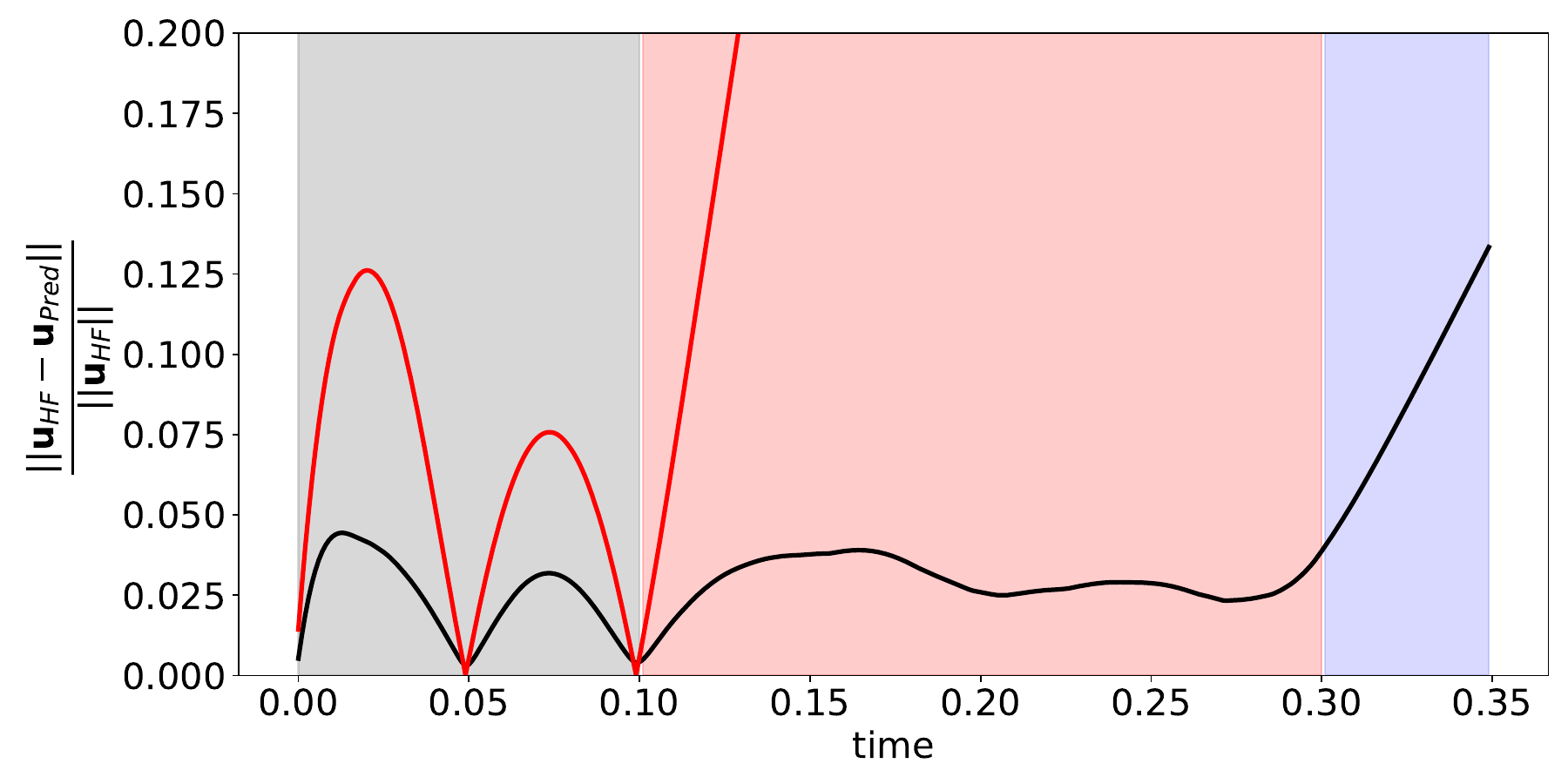}
        \caption{Loss without correction term}
        \label{fig:Error_Burgers_zoomed}
    \end{subfigure}
\medskip
\caption{Relative L2 error over time and two modeling strategies are evaluated: \textbf{Case I} - black line, which incorporates discretized physics via a physics-based loss term alongside a data-driven loss term, and \textbf{Case II} -   red line, which relies solely on the data-driven loss term. (a) Error distribution at three regimes (b) zoomed view close to maximum L2 error variation.}
\label{fig:loss_wwo_cor_data_physics}
\end{figure}

Now, we assess the prediction capability at a future time instant (which means the particular time instant is outside the training regime where the data-driven and physics-based residuals are considered). \autoref{fig:comaparison_0325_data_physics} shows absolute error distributions at $t = 0.325$ s while considering the combination of data and physics in the neural network and solely data-driven network. The data-driven (considered at $t = 0$, $t = 0.05$ and $t = 0.1$ s) loss terms are considered at the same time-instants as mentioned earlier. Whereas, the physics-based loss terms are considered starting from $0$ s to $0.3$ s in a time step interval of $0.01$ s. The maximum absolute errors incurred are $1.15$ and $0.075$, respectively, when only data-driven neural network-based prediction is considered and when physics-based residual is considered from the external solver alongside the data-driven loss terms. \autoref{fig:loss_wwo_cor_data_physics} shows the relative error varying over the time-instants starting from $0$ to $0.35$ s while the relative errors are computed as follows: 

\begin{figure}[ht]
\centering
\begin{subfigure}[b]{1\textwidth}
\centering
% \hspace{-3cm}
\subfloat[Prediction at $t = 0.325$ s]{\label{fig:Benchmark_pred}\includegraphics[width=0.75\linewidth]{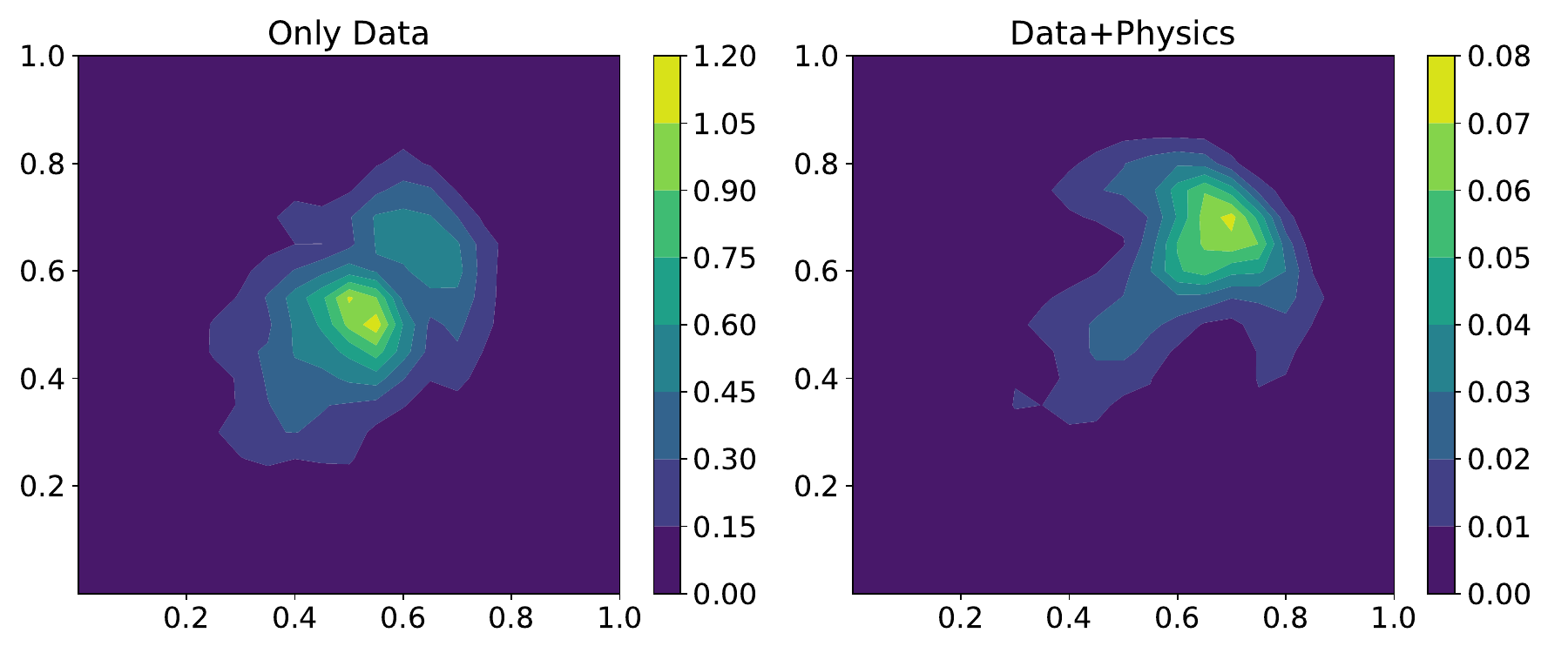}}
\end{subfigure}
\medskip
\caption{Absolute error at a future time instant of $0.325$ s with respect to the benchmark data when (a) only the data-driven loss term is considered, (b) both data-driven and physics-driven loss term is considered. For the numerical experiments, $3$ intermediate time-instants where numerical data ($0$, $0.05$, $0.1$ s) are considered for the data-driven loss term.  The physics-based residual is imported from the external CFD solver at a time interval of $0.01$ s between time instants of $0$ s to $0.3$ s.}
\label{fig:comaparison_0325_data_physics}
\end{figure}

%Now, we will discuss the effect of the additional correction we have introduced in the 

\begin{equation}\label{eq:L2error}
\text{Relative } L^2 \text{ error} = \frac{\| \mathbf{u}_{\text{HF}} - \mathbf{u}_{\text{pred}} \|_2}{\| \mathbf{u} \|_2}
= \frac{\left( \sum_{i=1}^{N} (\mathbf{u}_{\text{HF},i} - \mathbf{u}_{\text{pred},i})^2 \right)^{1/2}}{\left( \sum_{i=1}^{N} \mathbf{u}_{\text{HF},i}^2\right)^{1/2}},
\end{equation}

Where \( \mathbf{u}_{\text{HF}} \) denotes the high-fidelity prediction obtained using CFD, and \( \mathbf{u}_{\text{pred}} \) represents the prediction from the neural network. The total simulation time from \( 0 \) to \( 0.35\,\text{s} \) is divided into three regimes, and two modeling strategies are evaluated: \textbf{Case I}, which incorporates discretized physics via a physics-based loss term alongside data-driven loss term, and \textbf{Case II}, which relies solely on the data-driven loss term. In the \textbf{first regime} (\( 0 \) to \( 0.1\,\text{s} \)), solution snapshots are included at \( 0\,\text{s} \), \( 0.05\,\text{s} \), and \( 0.1\,\text{s} \) for computing the data-driven loss. Additionally, a physics-based loss term is constructed using residuals from an external solver, evaluated every \( 0.01\,\text{s} \). In this regime, the relative error for Case II is higher than that for Case I. In the \textbf{second regime}, no solution snapshots are used for the data-driven loss term. As a result, Case II exhibits a significantly higher relative error with a monotonically increasing trend. In contrast, Case I leverages the physics-based residuals computed at \( 0.01\,\text{s} \) intervals, resulting in improved accuracy and stability. The \textbf{third regime} lies outside the training window. Here, predictions rely solely on the physics-based loss term (up to \( 0.3\,\text{s} \)) and the data-driven loss term (up to \( 0.1\,\text{s} \)), based on extrapolated temporal information. Although both Case I and Case II show a monotonic increase in relative error, the error in Case I remains significantly lower than that in Case II.

\begin{figure}[ht]
\centering
\begin{subfigure}[b]{0.8\textwidth} % Reduced width for the top figure
\centering
\subfloat[Absolute prediction error when physics based loss term considered with ($\text{L}_\text{Dis}$) termed as (with correction) and with ($\text{L}_\text{eqn}$) termed as (without correction) ]{\label{fig:error_corrr}\includegraphics[width=1\linewidth]{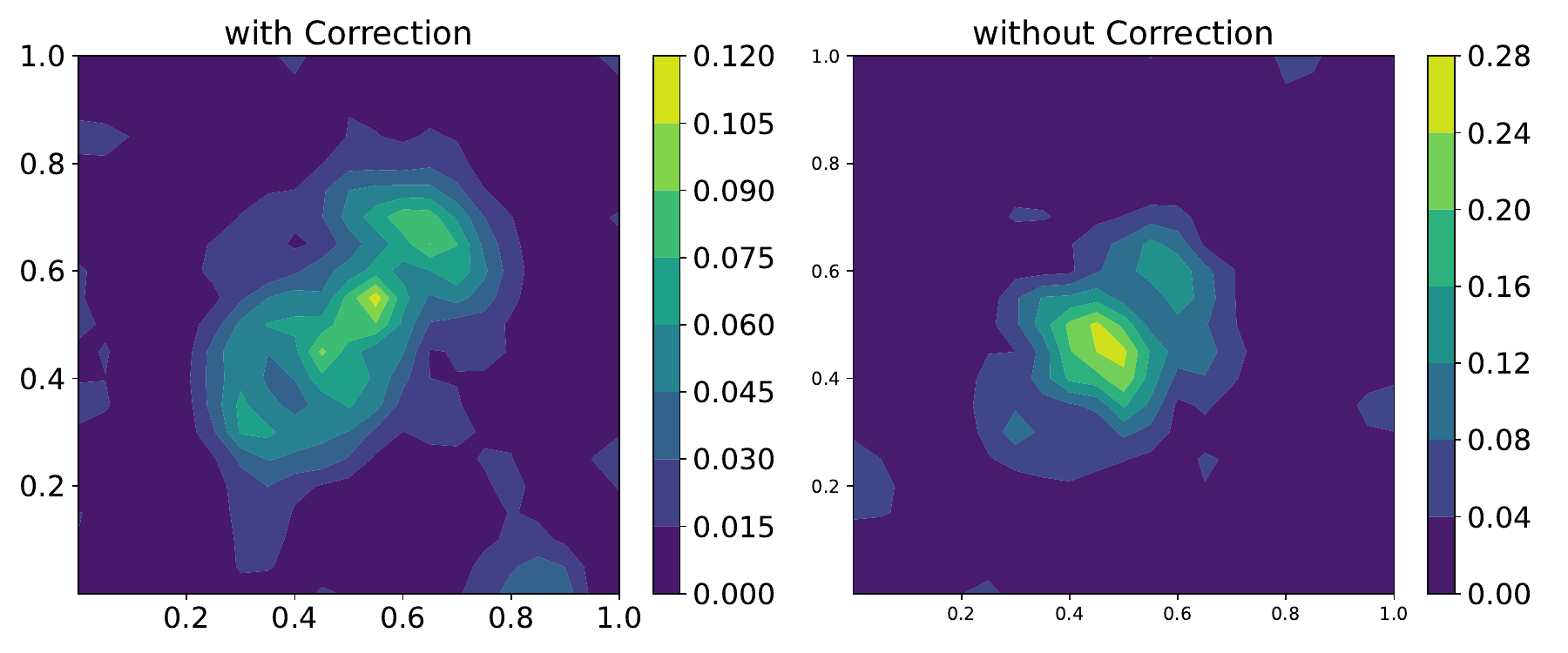}} % Reduced image width
\end{subfigure}
\medskip
\begin{subfigure}[b]{0.55\textwidth} % Further reduced width for bottom figures
    \centering
    \includegraphics[width=1\linewidth]{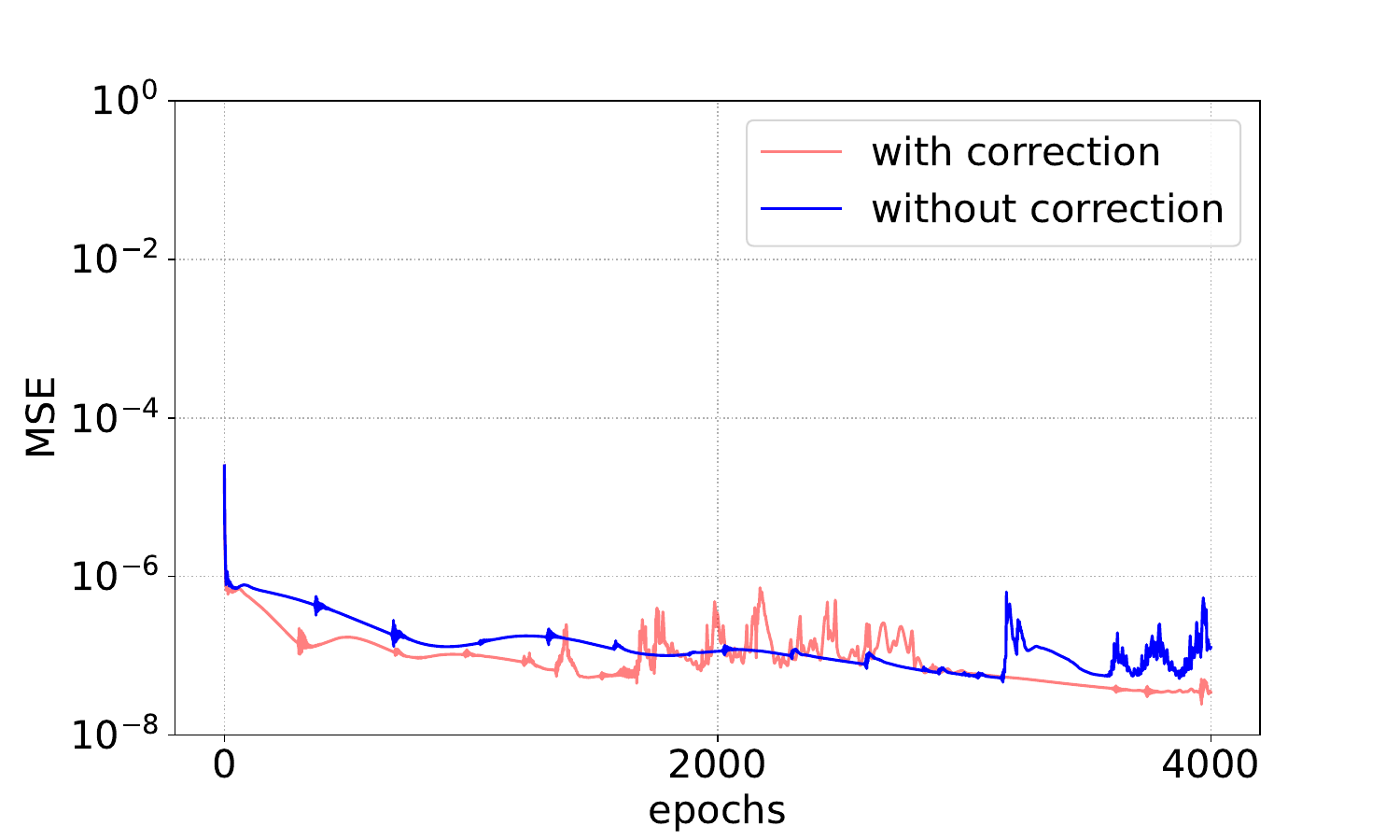} % Reduced image width
    \caption{MSE loss decay with number of epoch with and without the correction term - $\text{L}_\text{Dis}$ termed as (with correction) and with $\text{L}_\text{eqn}$ termed as (without correction).}
    \label{fig:loss_w_corr}
\end{subfigure}
\medskip
\caption{Examining the effect of introducing an additional correction in the computation of the physics-based loss term for the discretized PINN. The solution snapshot at the initial time (t = 0) is only considered for the data-driven loss term. The physics-based residuals are imported from the external CFD solver at a time interval of $0.01$ s.}
\label{fig:loss_wwo_cor1}
\end{figure}

We now examine the effect of introducing an additional correction in the computation of the physics-based loss term for the discretized PINN. While carrying out this experiment, the solution snapshot at the initial time (t = 0) is only considered for the data-driven loss term. The physics-based residuals are imported from the external CFD solver at a time interval of $0.01$ s.  In particular, we consider $\text{L}_\text{dis}$, as defined in \autoref{eq:LDIS_NTE}, instead of the standard residual loss $\text{L}_\text{eqn}$ shown in \autoref{eq:Leqn_NTE}, to better capture the nonlinearity present in the convective term. The rationale behind this correction lies in the treatment of the system matrix $\mathbf{A}$ and the vector $\mathbf{b}$ during backpropagation. In the formulation of $\text{L}_\text{eqn}$, both $\mathbf{A}$ and $\mathbf{b}$ are assumed to be constant. However, since  $\mathbf{A}$ or $\mathbf{b}$ or both can depend on the solution variable $\mathbf{u}$, this assumption introduces inaccuracies in the gradient calculation. The correction term in $\text{L}_\text{dis}$ addresses this issue by explicitly accounting for the dependency of $\mathbf{A}$ and $\mathbf{b}$ on $\mathbf{u}$. We compare the absolute error of the solution field with and without the additional correction term in \autoref{fig:error_corrr}. In this numerical experiment, we have considered the initial condition as the data-driven loss term, and the physics-driven loss term is computed from time instant $0$ s to $0.2$ s in steps of $0.01$ s. The maximum absolute error associated with the prediction with and without the consideration of the correction term is $0.11$ and $0.25$, which demonstrates the effectiveness of the additional correction introduced in the physics-based loss computation. In addition, we also compare the decay of the physics-driven and data-driven $\text{MSE}$ loss terms with and without the correction terms in \autoref{fig:loss_w_corr}. The physics-based loss term drops more in the case of the one with the correction term $(3.5\times10^{-8})$ as compared to the one without correction $(1.3\times10^{-7})$. The spectral properties of the coefficient matrix $\mathbf{A}$ are expected to influence the effectiveness of the optimization procedure employed to minimize the residuals obtained from \autoref{eq:LDIS_NTE} and \autoref{eq:Leqn_NTE}, which will be thoroughly investigated in the authors' subsequent works on this topic. Readers are also referred to the discussion on the decay of the MSE loss term associated with the physics-based loss term in the next section, within the context of the reduced-order system.

\subsection{Reduced-Order Unsteady Incompressible Flow Problem}\label{sec:incomp}

\begin{figure}[ht]
           \centering

        % First Row
        \begin{subfigure}[b]{1\textwidth}
            \centering
            \subfloat[Modes = $1$]{\label{fig:umodes1}\includegraphics[width=.25\linewidth]{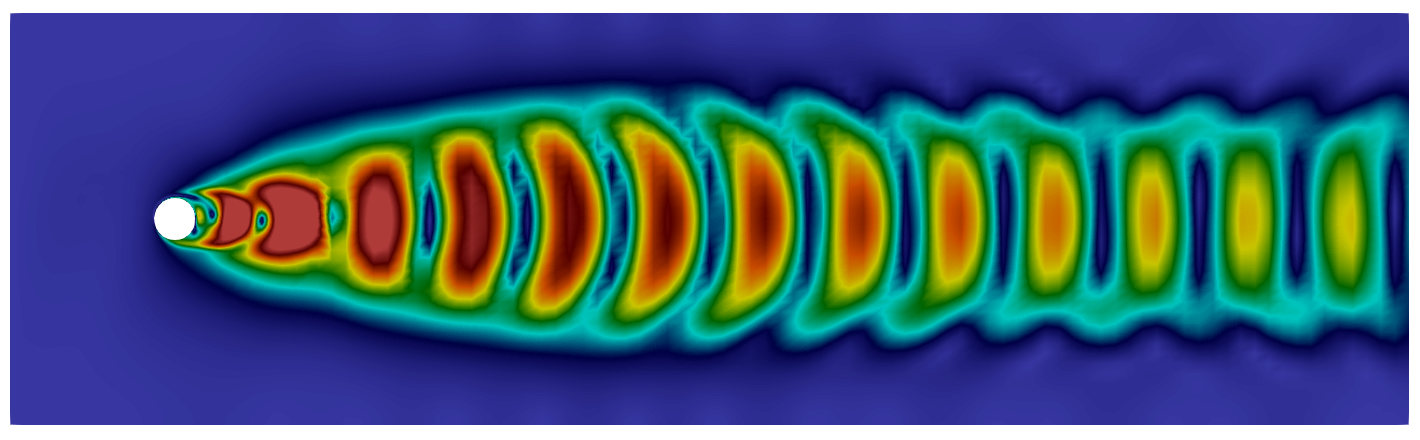}}
            \subfloat[Modes = $2$]{\label{fig:umodes2}\includegraphics[width=.25\linewidth]{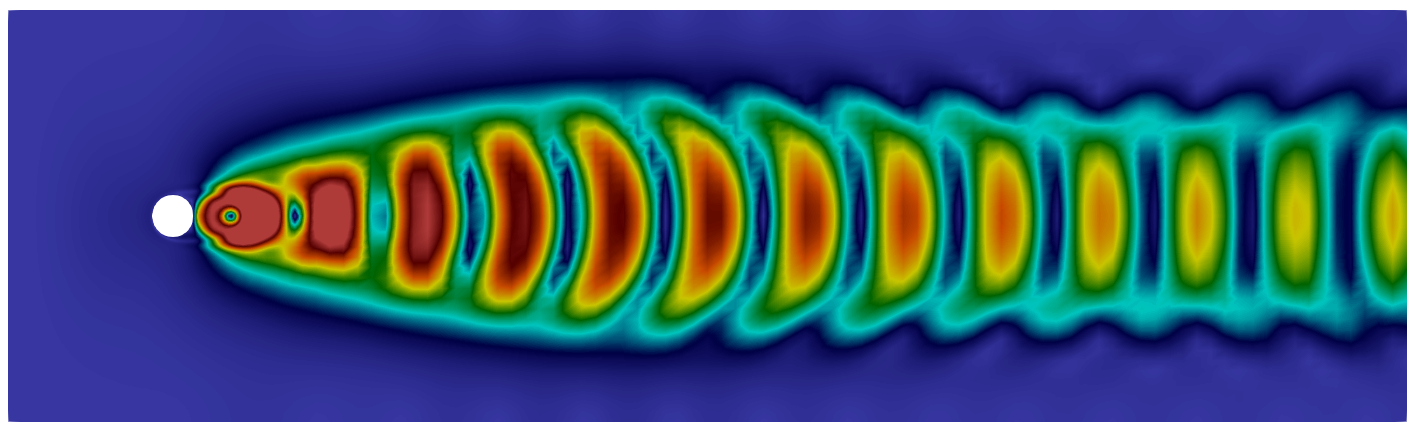}}
            \subfloat[Modes = $3$]{\label{fig:umodes3}\includegraphics[width=.25\linewidth]{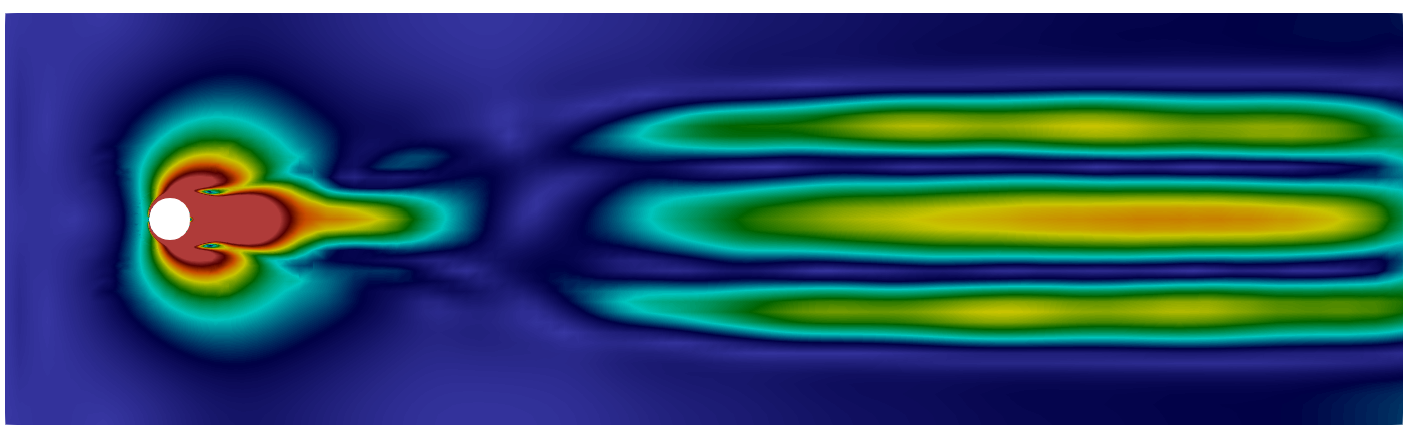}}
        \end{subfigure}

        \vspace{0.3cm} % Space between rows

        % Second Row
        \begin{subfigure}[b]{1\textwidth}
            \centering
            \subfloat[Modes = $4$]{\label{fig:umodes4}\includegraphics[width=.25\linewidth]{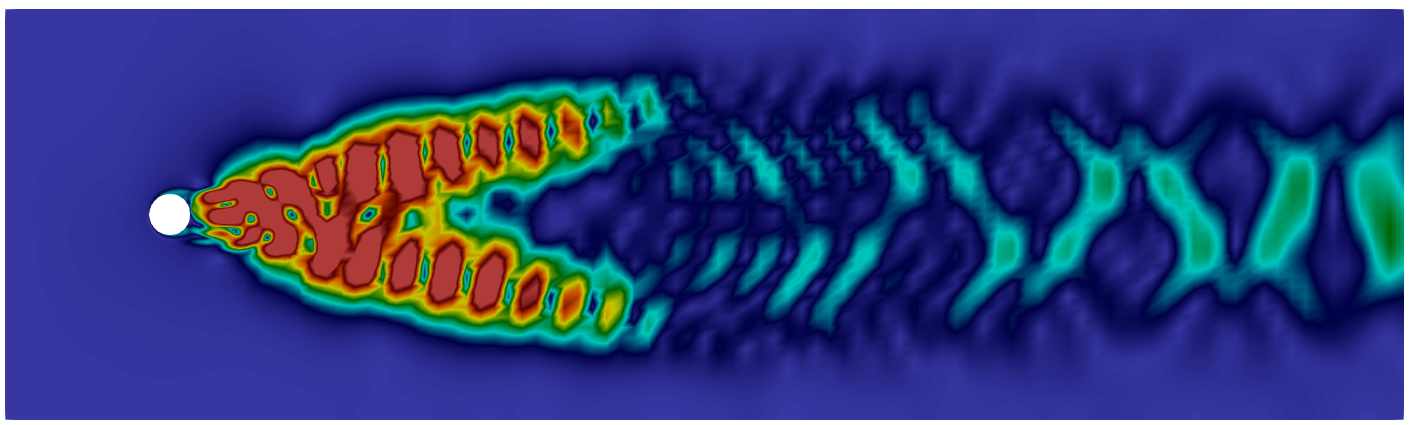}}
            \subfloat[Modes = $5$]{\label{fig:umodes5}\includegraphics[width=.25\linewidth]{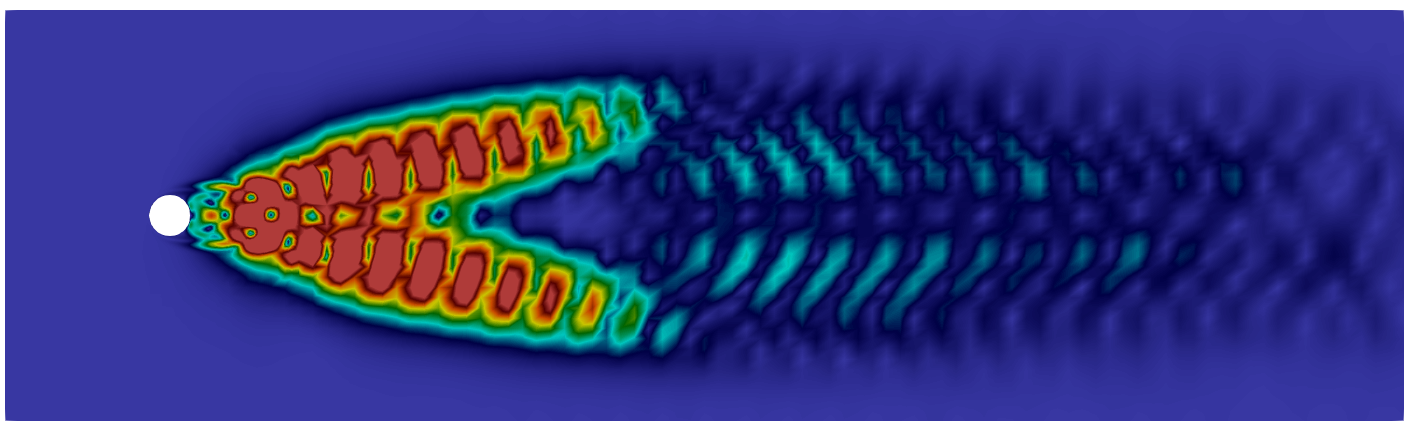}}
            \subfloat[Modes = $6$]{\label{fig:umodes6}\includegraphics[width=.25\linewidth]{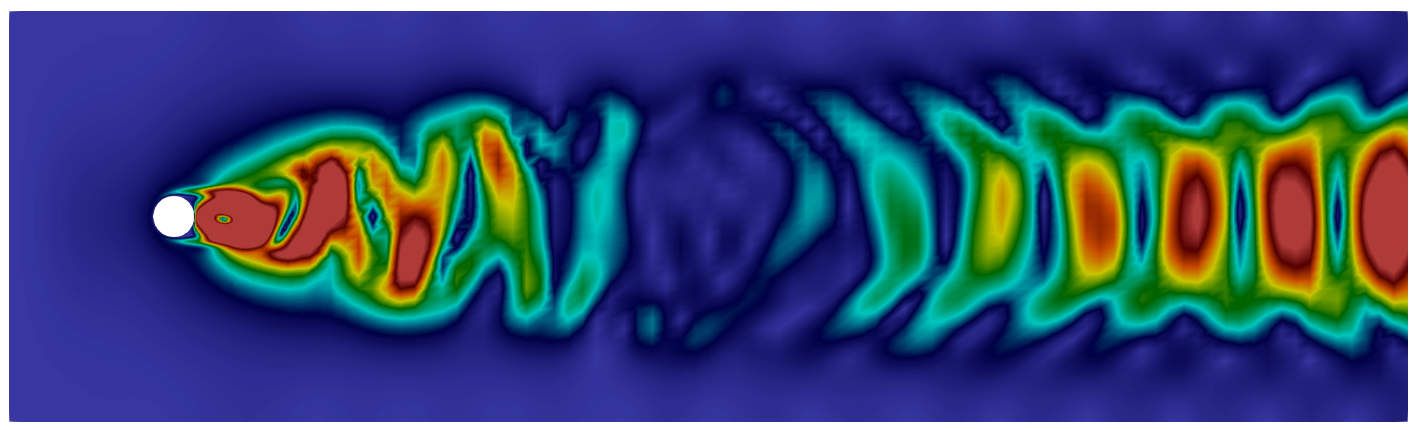}}
        \end{subfigure}

    \begin{subfigure}[b]{1\textwidth}
            \centering
{\label{fig:uscale1}\includegraphics[width=0.3\linewidth]{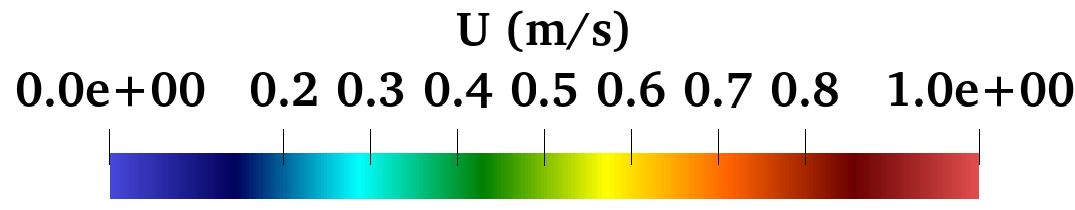}}          
        \end{subfigure}

    \medskip
    \caption{First $6$ velocity modes for the reduced order system of the incompressible flow past a cylinder.}
    \label{fig:vmodes}
\end{figure}

\begin{figure}[ht]
           \centering

        % First Row
        \begin{subfigure}[b]{1\textwidth}
            \centering
            \subfloat[Modes = $1$]{\label{fig:pmodes1}\includegraphics[width=.25\linewidth]{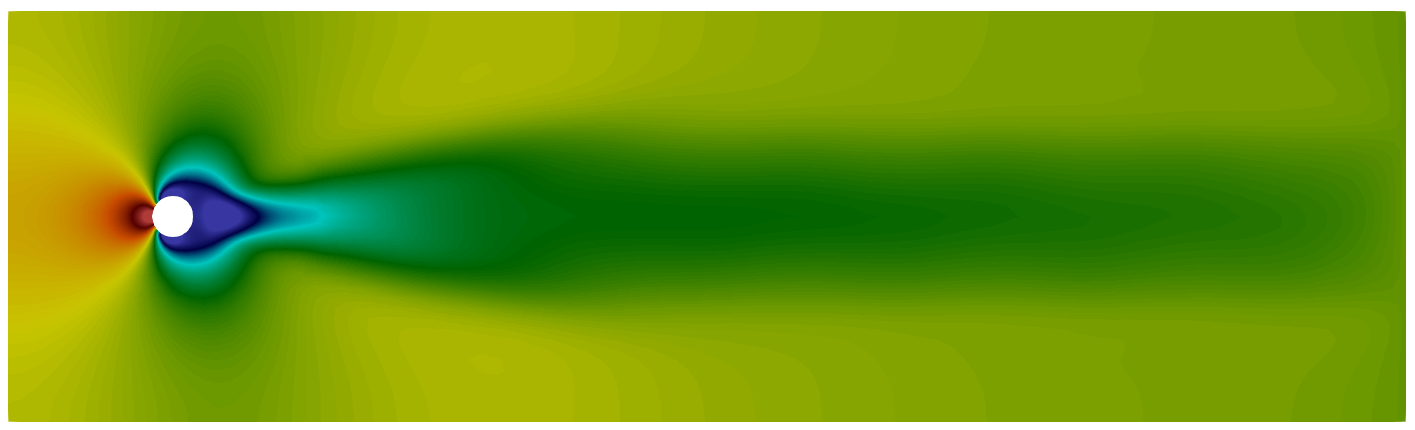}}
            \subfloat[Modes = $2$]{\label{fig:pmodes2}\includegraphics[width=.25\linewidth]{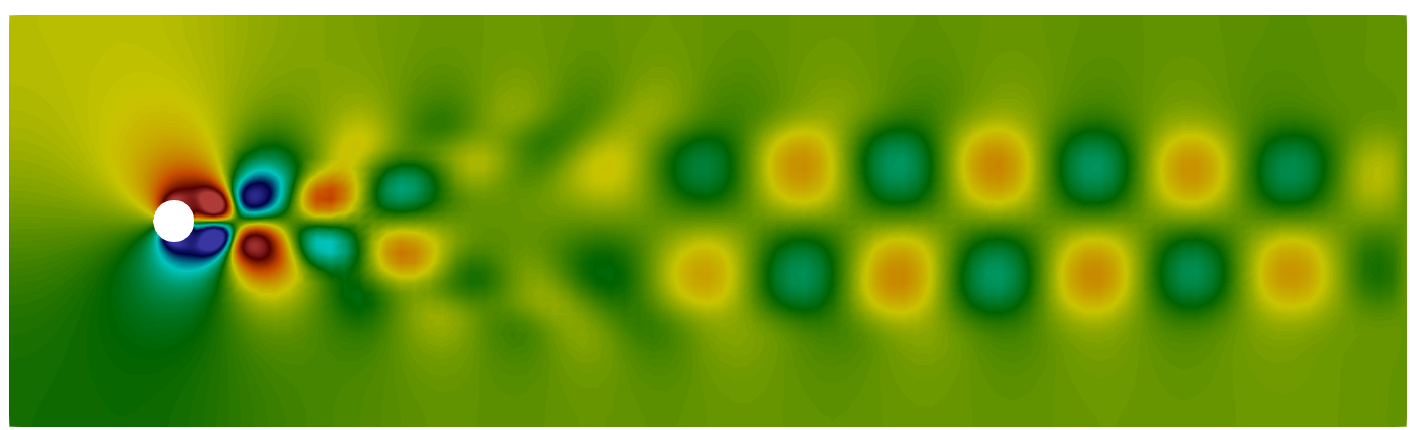}}
            \subfloat[Modes = $3$]{\label{fig:pmodes3}\includegraphics[width=.25\linewidth]{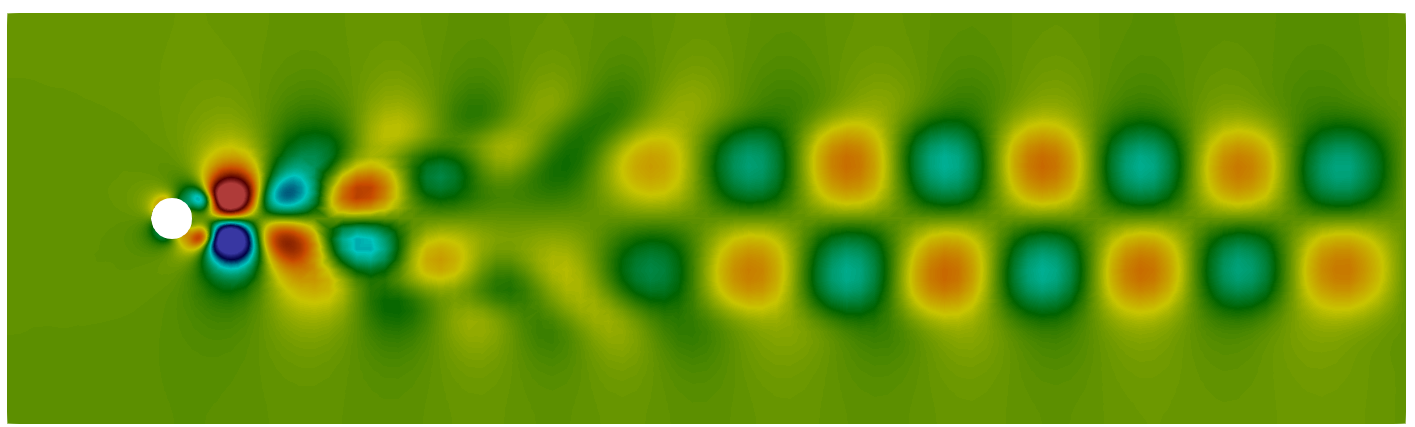}}
        \end{subfigure}

        \vspace{0.3cm} % Space between rows

        % Second Row
        \begin{subfigure}[b]{1\textwidth}
            \centering
            \subfloat[Modes = $4$]{\label{fig:pmodes4}\includegraphics[width=.25\linewidth]{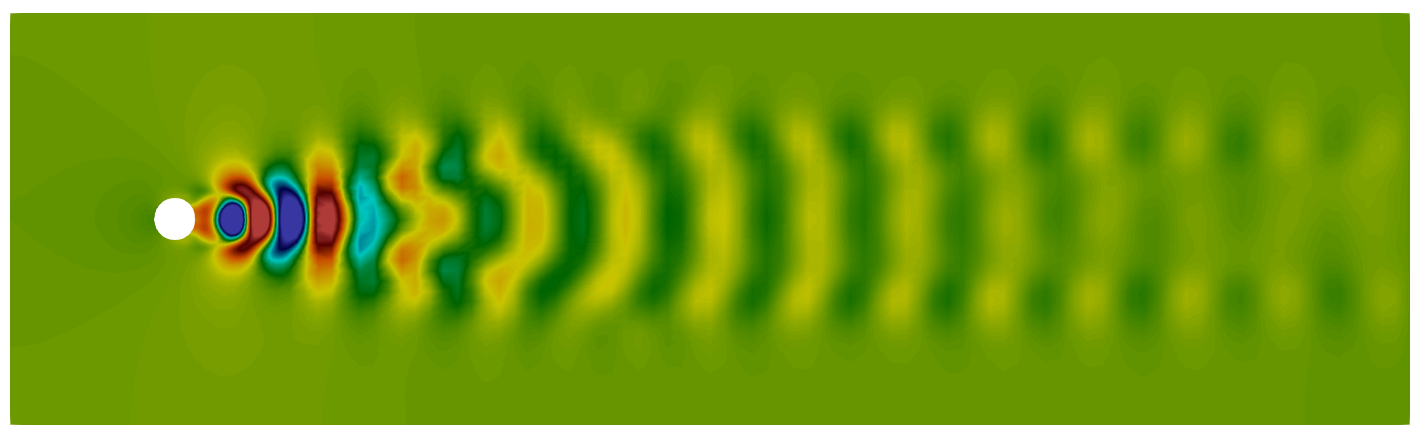}}
            \subfloat[Modes = $5$]{\label{fig:pmodes5}\includegraphics[width=.25\linewidth]{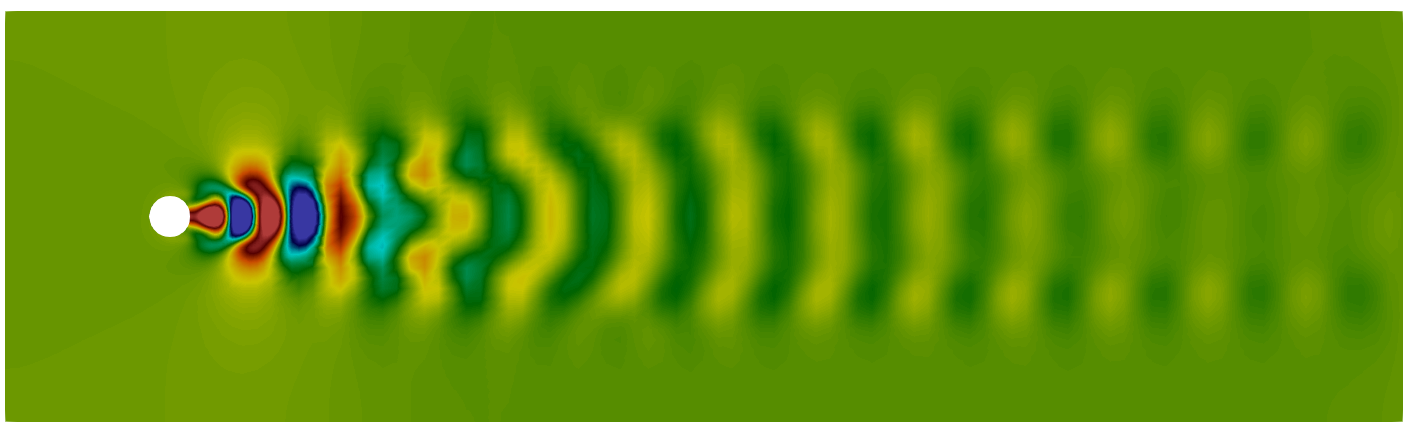}}
            \subfloat[Modes = $6$]{\label{fig:pmodes6}\includegraphics[width=.25\linewidth]{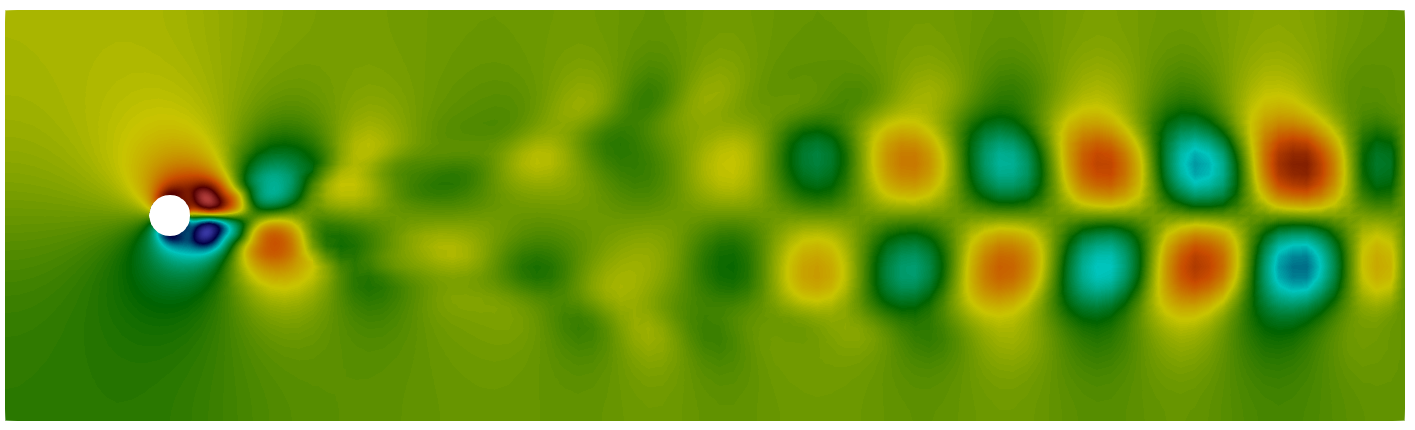}}
        \end{subfigure}

    \begin{subfigure}[b]{1\textwidth}
            \centering
{\label{fig:uscale2}\includegraphics[width=0.3\linewidth]{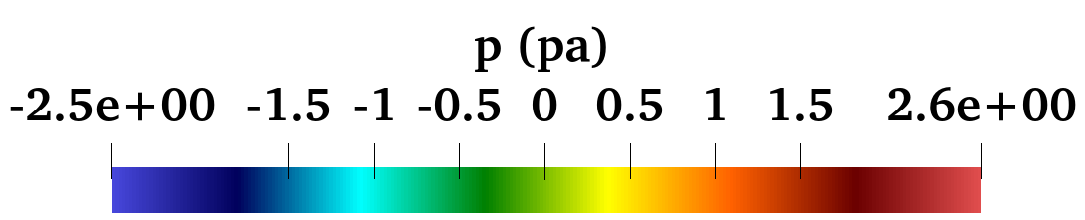}}          
        \end{subfigure}

    \medskip
    \caption{First $6$ pressure modes for the reduced order system of incompressible flow past a cylinder}
    \label{fig:pmodes}
\end{figure}

This example aims to apply the proposed methodology to a more complex flowfield and mesh structure, a benchmark flow over a circular cylinder. To reduce the degrees of freedom associated with the computational domain, we introduce a reduced-order model, which is then coupled with a physics-informed neural network. This work introduces physical parametrization due to parametrized physical viscosity $\nu$. First, we elaborate on the computational model considered for the full-order system. The number of computational cells is $11644$. The time step is considered $0.01$, sufficiently small to meet the CFL condition in all parts of the domain. Now, to train the reduced-order model, two kinematic viscosities $\nu$ of $0.005$ and $0.006$ $\text{m}^2/\text{s}$ are considered. In these experiments, the simulation is run for each kinematic viscosity up to $60$ $\text{s}$.  We have considered only the last $10$ s of the flow physics to generate the POD basis. Within these time windows, the snapshots are collected at every $0.01$ s, resulting in $1000$ snapshots for every $\nu$. Therefore, the total number of snapshots for the training of the neural network is $2000$. Now, the snapshots of two full-order simulations are used for the POD basis generation. $20$ modes are considered for velocity and $10$ modes for the pressure and $10$ supremizers modes for the generation of the reduced order system as shown in \autoref{eq:red_system1} and \autoref{eq:red_system2}. \autoref{fig:vmodes} and \autoref{fig:pmodes} show the first $6$ velocity and pressure modes, respectively. It is worth noting that the number of modes used for the velocity and pressure fields significantly affects both the accuracy and stability of the reduced-order solution. Increasing the number of modes generally improves the accuracy of the reduced-order simulation; however, it also increases the dimensionality of the neural network, which may lead to poor training performance. Therefore, the number of modes considered here is chosen to strike a balance between computational efficiency and prediction accuracy. Now, the obtained reduced-order system is utilized as the physics-based residual for the neural network, as shown in \autoref{eq:matrices1}. 

\begin{figure}[ht]
\centering
    \begin{subfigure}[b]{0.49\textwidth}
        \centering
        \includegraphics[width=0.85\linewidth]{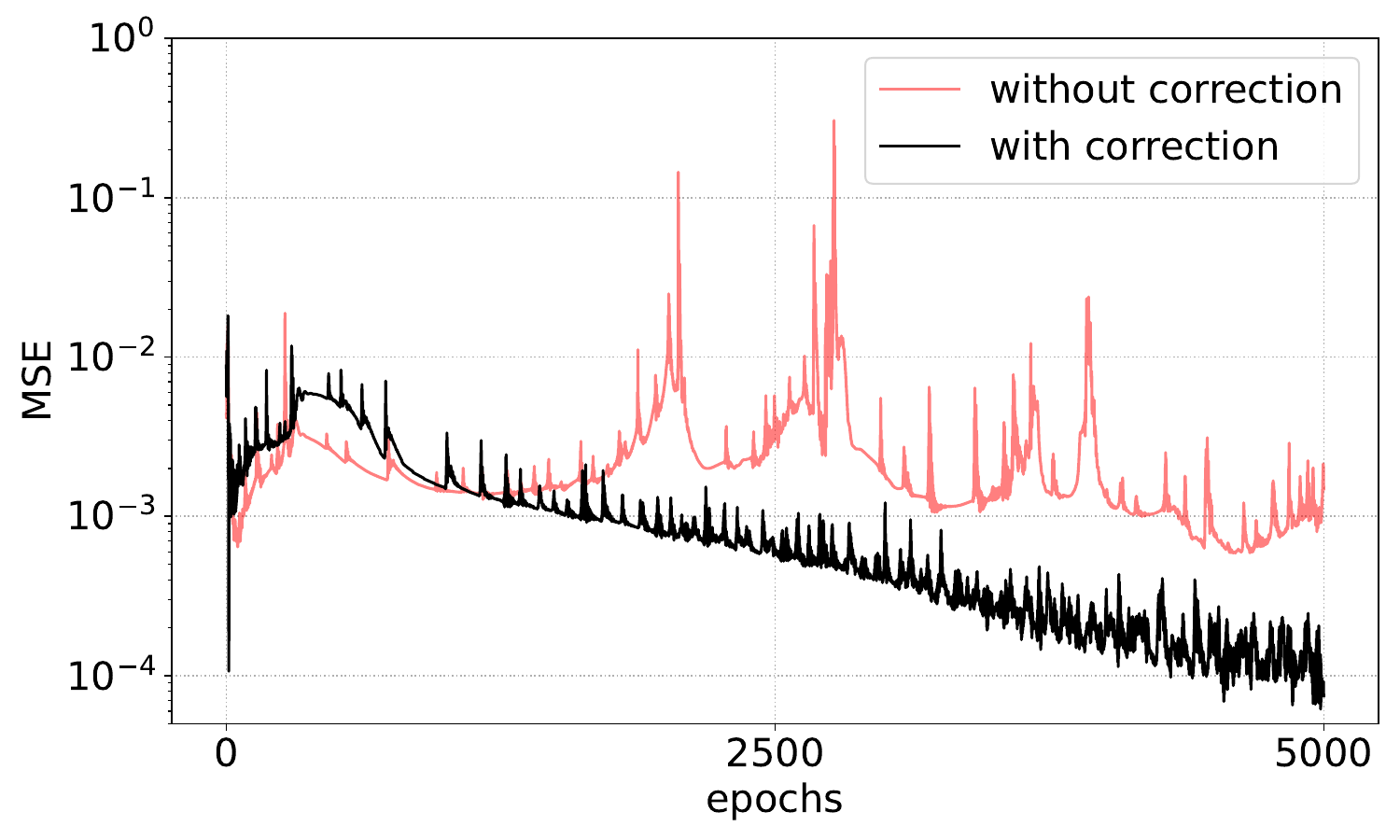}
        \caption{Minimization of MSE loss term based on ${\mathbf{R}_\text{red1}}$}
        \label{fig:MSE_Momentum}
    \end{subfigure}
    \hfill
    \begin{subfigure}[b]{0.49\textwidth}
        \centering
        \includegraphics[width=0.85\linewidth]{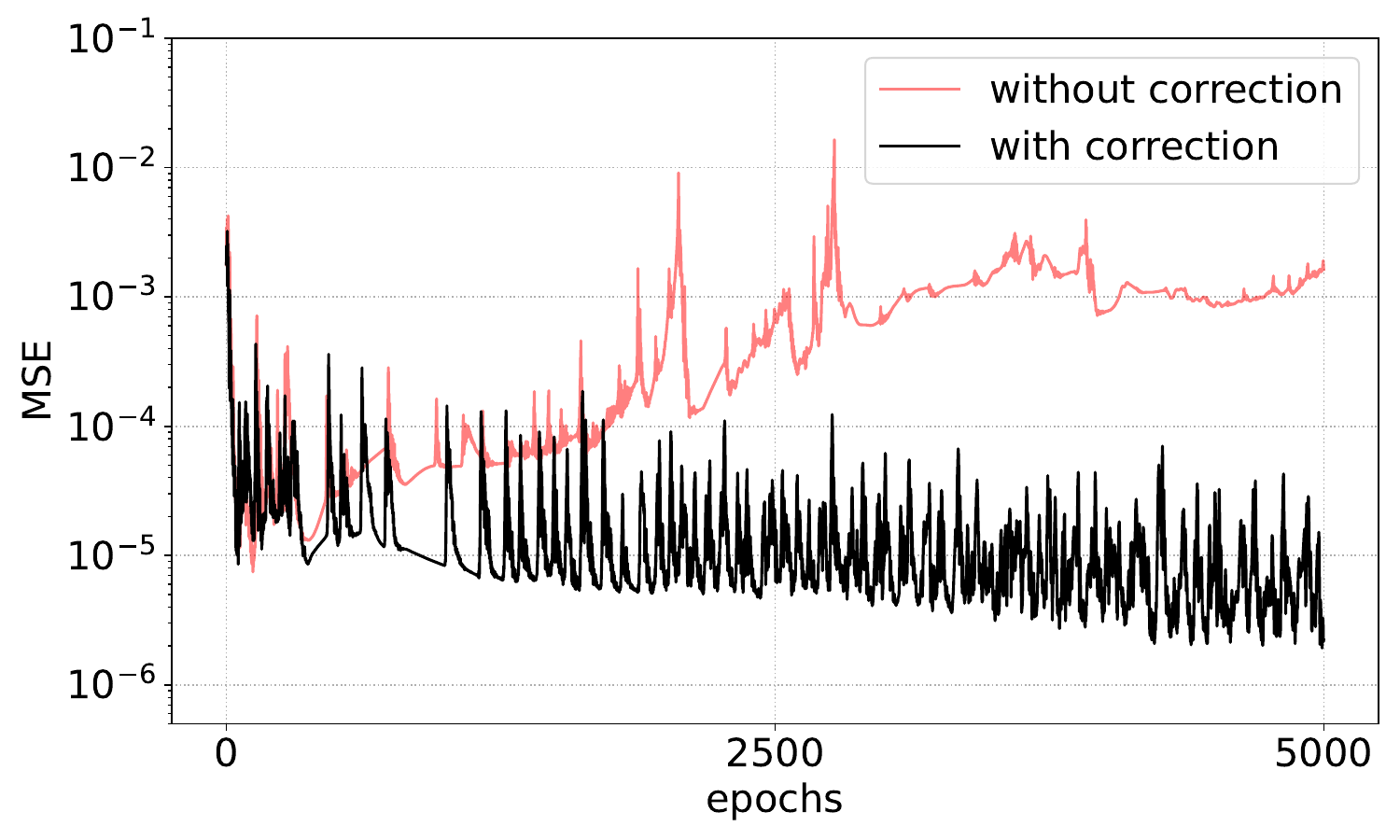}
        \caption{Minimization of loss term based on $\mathbf{R}_\text{red2}$}
        \label{fig:MSE_Preesure}
    \end{subfigure}
\medskip
\caption{Comparison of the performance of the loss terms $\mathrm{L}_{\text{Dis}}$ (denoted as with correction) and $\mathrm{L}_{\text{eqn}}$ (denoted as without correction) associated with the reduced formulation of (a) momentum and (b) pressure equation.}
\label{fig:loss_wwo_cor}
\end{figure}

For the development of the neural network, $4$ layers with $124, 100, 80$ and $64$ neurons are considered, respectively. $\text{softplus}$ activation function is considered and the learning rate is $0.001$. The total number of the $\text{epoch}$ is $5000$. In this particular example, our objective is to develop a discretized PINN when the neural network can not include the \autoref{eq:red_system1} and \autoref{eq:red_system2} in its computational graph. The input to the neural network consists of time and the viscosity parameter $\nu$. The output corresponds to the POD coefficients of the reduced-order system. The time input ranges from $50$ s to $60$ s with a step size of $0.01$ s, while two different values of the physical parameter $\nu$ ($0.005$ and $0.006$ $\text{m}^2/\text{s}$) are considered, resulting in an input size of $2000 \times 2$. The output of the neural network comprises $20$ velocity modes, $10$ pressure modes, $10$ supremizer modes, and $1$ lift mode, yielding an output size of $2000 \times 41$.  As discussed in \autoref{NN-coupling}, the primary objective of this work is to reconstruct a spatio-temporal field using the proposed discretized PINN framework, given partial temporal data available at specific time instances. This data could originate from either experimental measurements or numerical simulations; however, in this study, we focus solely on numerical data. In addition, the proposed framework is employed to predict the solution field at future time instances or for different values of the viscosity parameter $\nu$. Now, we first assess the effectiveness of our proposed physics-based loss term, $\mathrm{L}_{\text{Dis, Mom}}$ and $\mathrm{L}_{\text{Dis, pressure}}$ as defined in \autoref{eq:LDIS_ROM1} and \autoref{eq:LDIS_ROM2} instead of using $\mathrm{L}_{\text{eqn, mom}}$ and $\mathrm{L}_{\text{eqn, pressure}}$ as shown in \autoref{eq:leqn_mom} and \autoref{eq:leqn_pressure} within the backpropagation algorithm of the neural network. While carrying out this assessment, for the data-driven loss term, we used numerical datasets corresponding to $5$ time steps $[0, 250, 500, 750, 1000]$ and the parameters $\nu = 0.005$ and $\nu = 0.006$ $\text{m}^2/\text{s}$. \autoref{fig:MSE_Momentum} illustrates the decay of the MSE loss term associated with the reduced formulation of the momentum equation, i.e., \autoref{eq:red_system1}, while \autoref{fig:MSE_Preesure} depicts the MSE loss term associated with the reduced formulation of the pressure equation, i.e., \autoref{eq:red_system2}. Both figures compare the performance of the loss terms $\mathrm{L}_{\text{Dis}}$ (denoted as with correction) and $\mathrm{L}_{\text{eqn}}$ (denoted as without correction). As shown, the MSE loss terms decay more rapidly for both equations while using the corrected loss terms, i.e., $\text{L}_{\text{Dis, Mom}}$ and $\text{L}_{\text{Dis, pressure}}$, ensuring effective neural network training. 
  
\begin{figure}[ht]
           \centering

        % First Row
        \begin{subfigure}[b]{1\textwidth}
            \centering
            \subfloat[t = $50.5s$]{\label{fig:FPC_benchmark_50}\includegraphics[width=.25\linewidth]{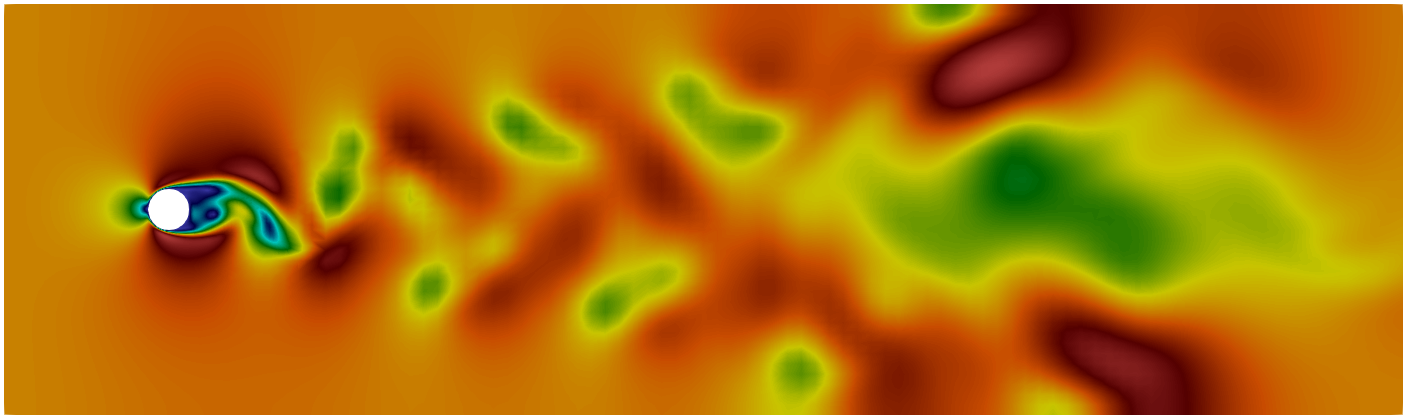}}
            \subfloat[t = $55s$]{\label{fig:FPC_benchmark_500}\includegraphics[width=.25\linewidth]{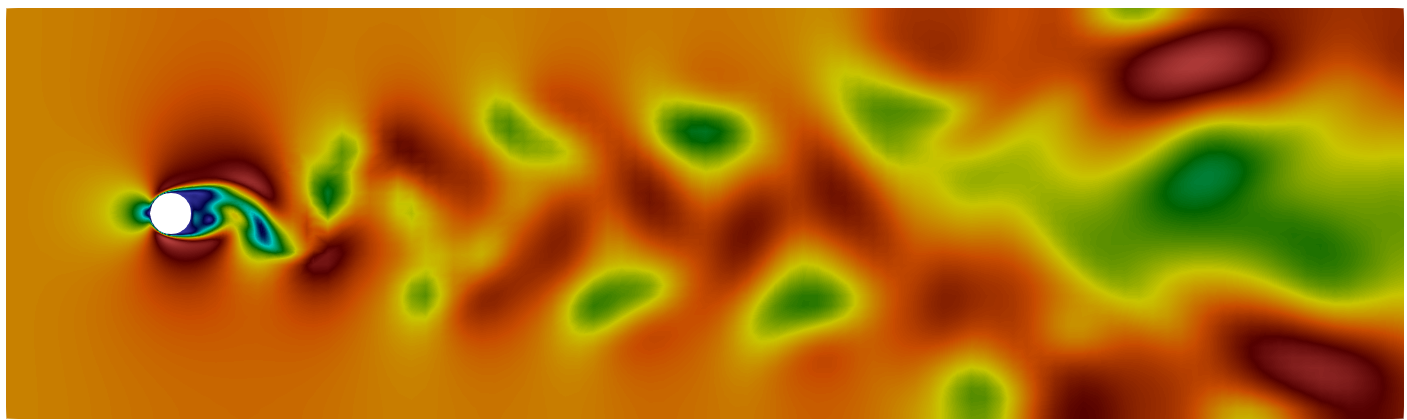}}
            \subfloat[t = $58s$]{\label{fig:FPC_benchmark_800}\includegraphics[width=.25\linewidth]{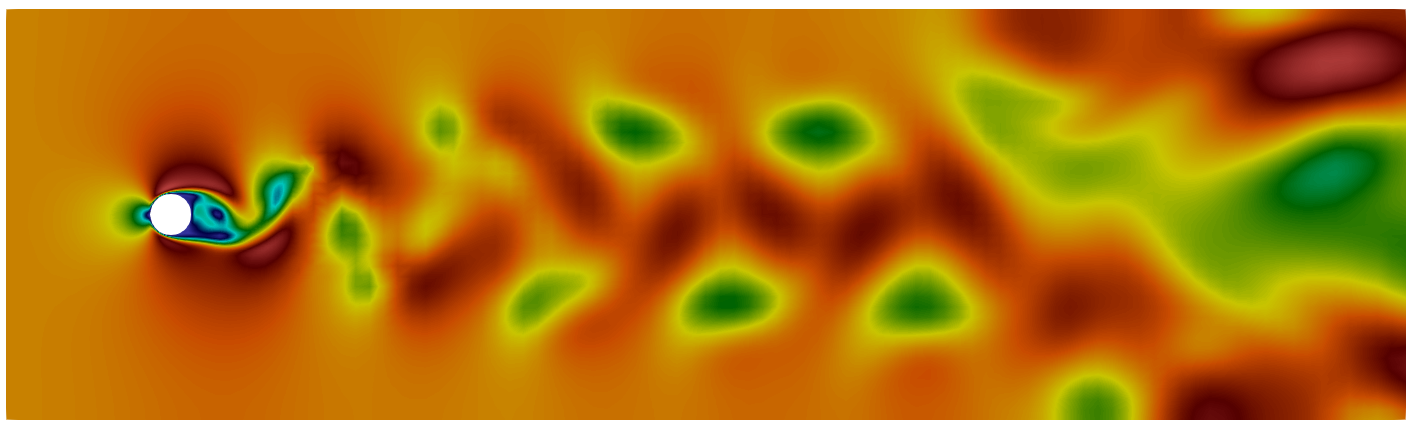}}
        \end{subfigure}

        \vspace{0.3cm} % Space between rows

        % Second Row
        \begin{subfigure}[b]{1\textwidth}
            \centering
            \subfloat[t = $50.5s$]{\label{fig:FPC_prediction_50}\includegraphics[width=.25\linewidth]{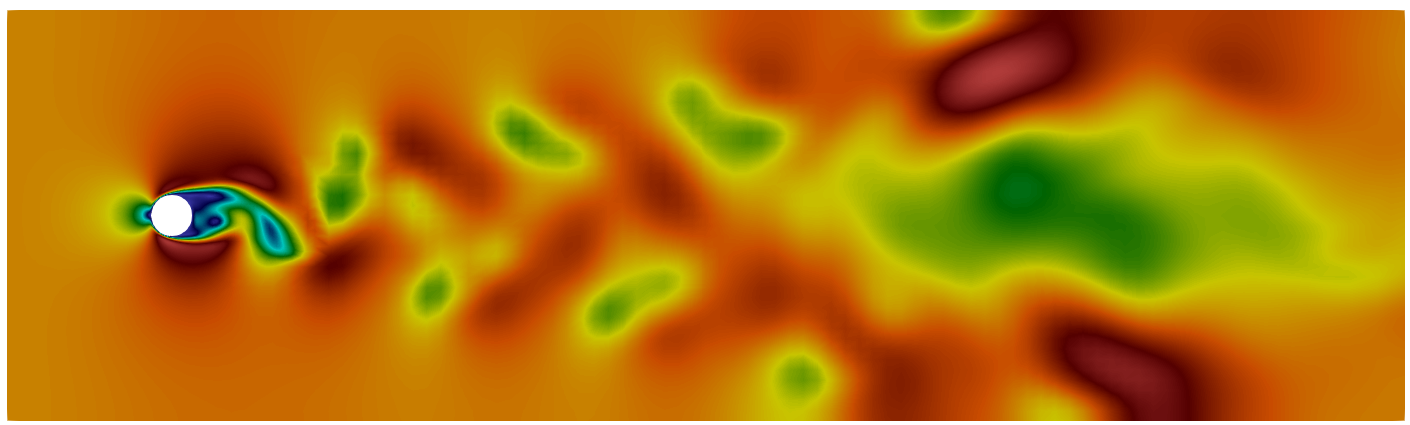}}
            \subfloat[t = $55s$]{\label{fig:FPC_prediction_500}\includegraphics[width=.25\linewidth]{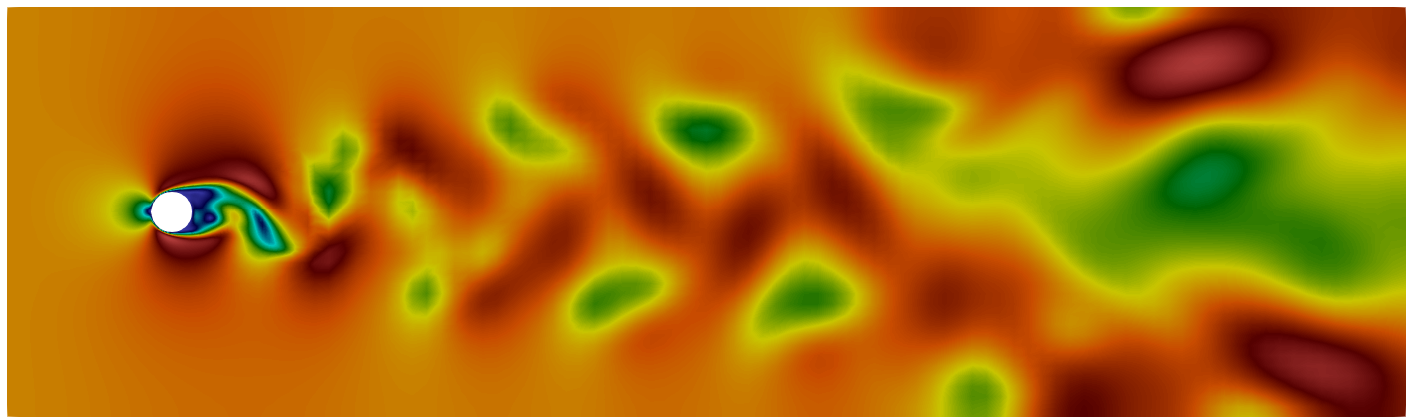}}
            \subfloat[t = $58s$]{\label{fig:FPC_prediction_800}\includegraphics[width=.25\linewidth]{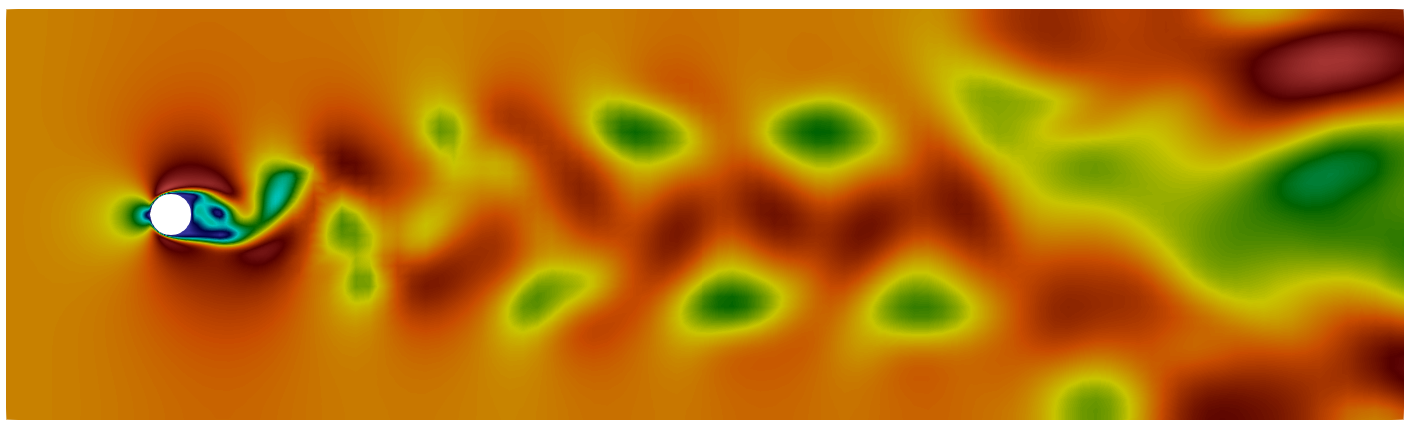}}
        \end{subfigure}

    \begin{subfigure}[b]{1\textwidth}
            \centering
{\label{fig:FPC_benchmark_50_4}\includegraphics[width=0.25\linewidth]{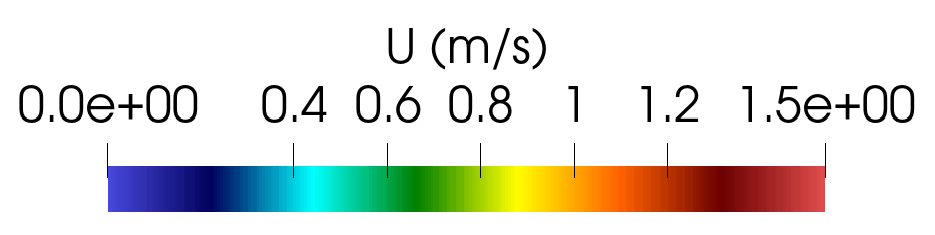}}          
        \end{subfigure}

    \medskip
    \caption{Comparison of the velocity field between (a) benchmark CFD results and (b) ANN-DisPINN prediction results of the reduced order incompressible viscous flow past a cylinder at $\nu=0.005$ and $t = $ $50.5$ $s$, $55$ $s$ and $58$ $s$}.
    \label{fig:FPC_u}
\end{figure}

\begin{figure}[ht]
\centering
\begin{subfigure}[b]{0.90\textwidth}
\centering
\subfloat[relative error, $\nu=0.005$ $\text{m}^2/\text{s}$ ]{\label{fig:L2_Error_005}\includegraphics[width=.50\linewidth]{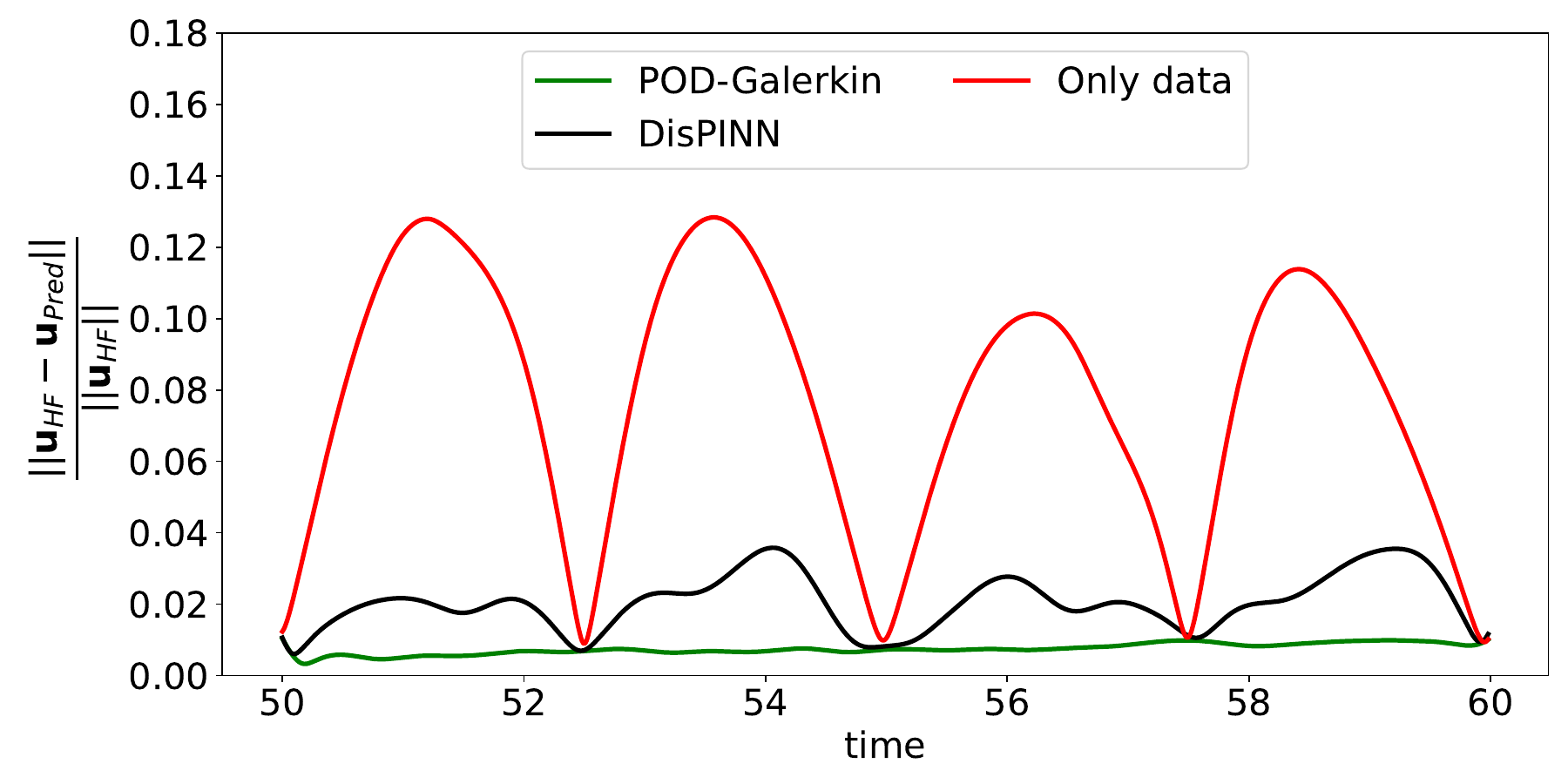}}
\subfloat[relative error, $\nu=0.0055$ $\text{m}^2/\text{s}$]{\label{fig:L2_Error_0055}\includegraphics[width=.50\linewidth]{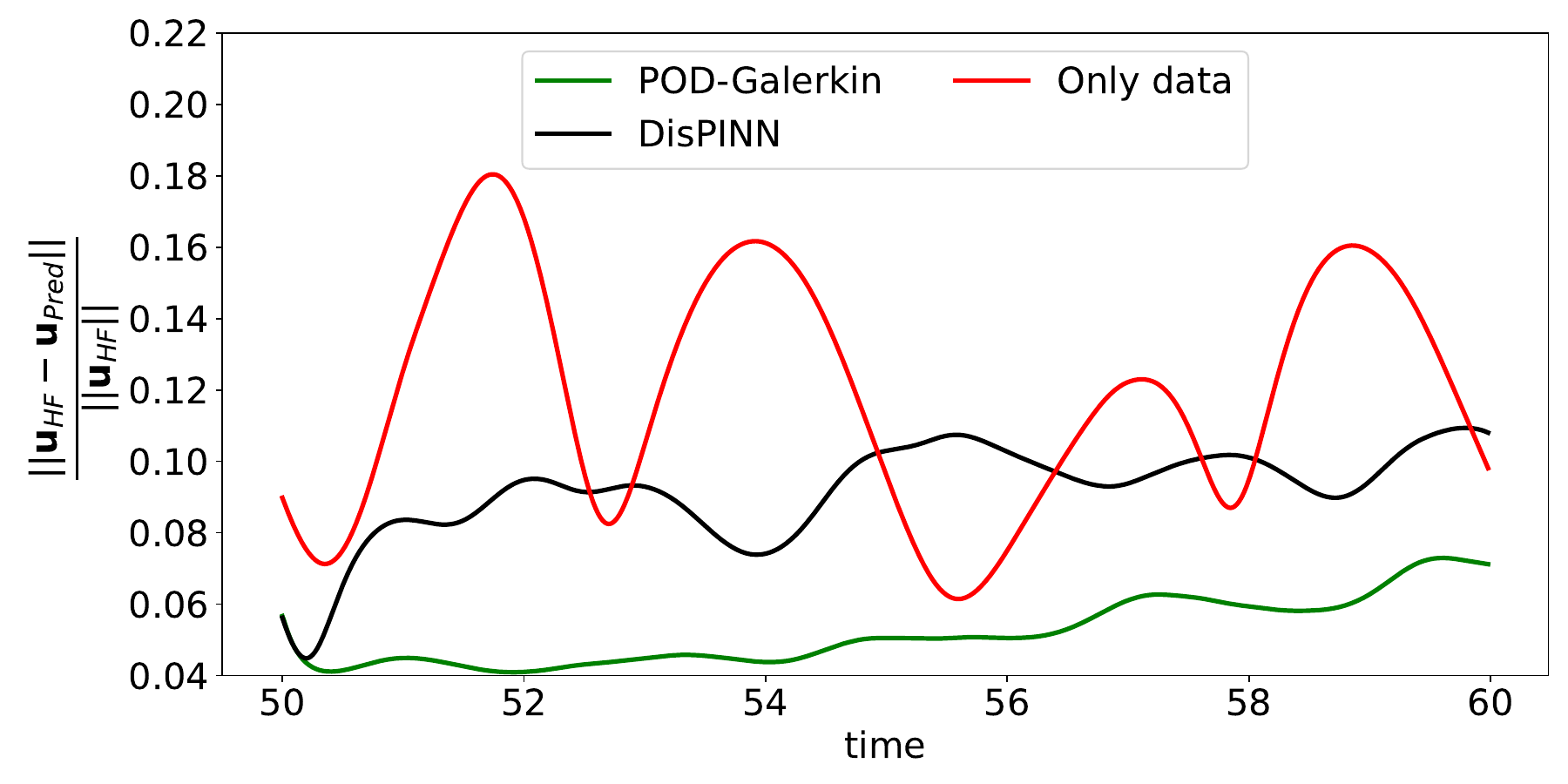}}
\end{subfigure}
\begin{subfigure}[b]{0.90\textwidth}
\centering
\subfloat[relative error, $\nu=0.005$ $\text{m}^2/\text{s}$]{\label{fig:L2_Error_005_pressure}\includegraphics[width=.50\linewidth]{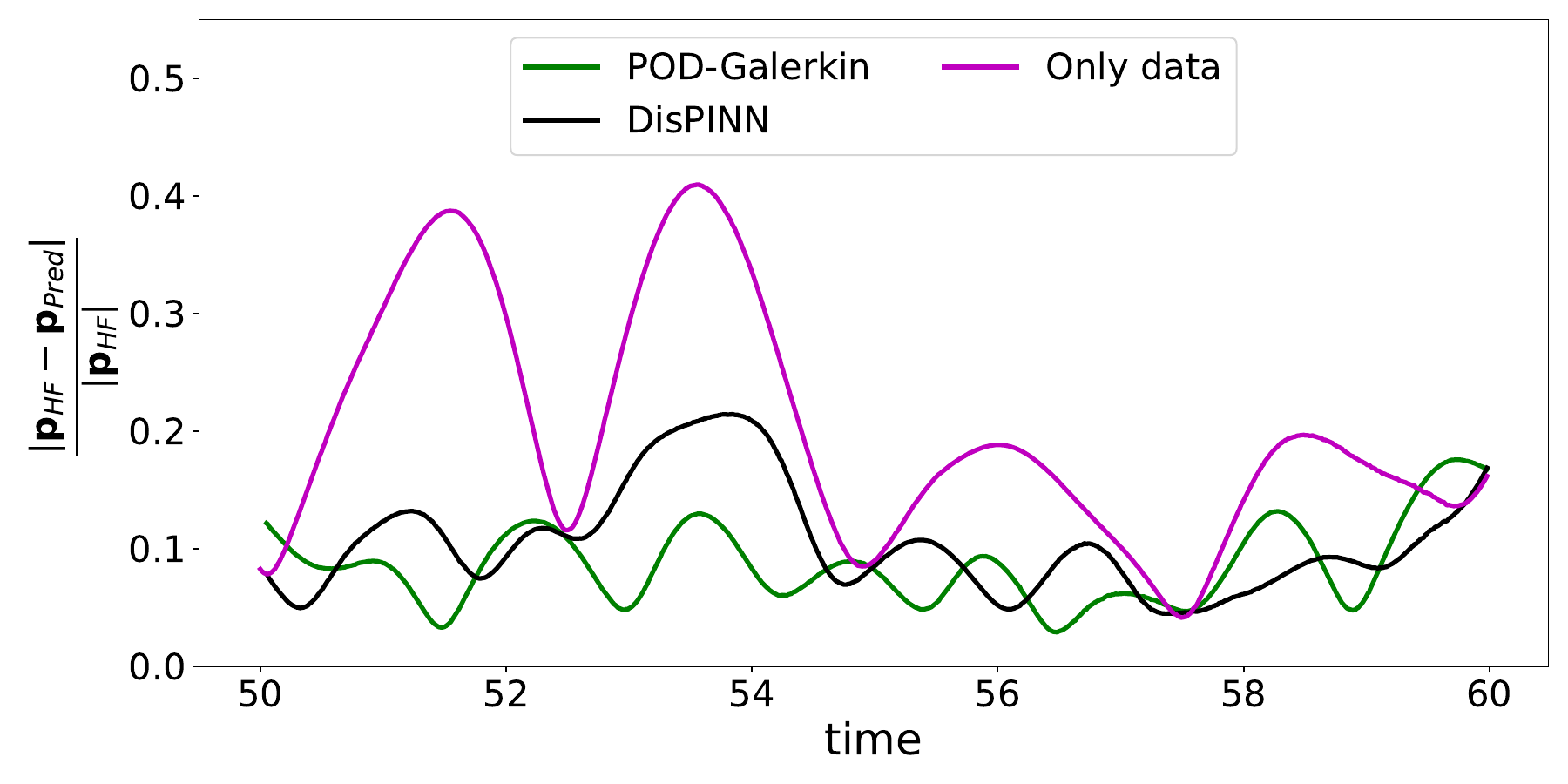}}
\subfloat[relative error, $\nu=0.0055$ $\text{m}^2/\text{s}$]{\label{fig:L2_Error_0055_pressure}\includegraphics[width=.50\linewidth]{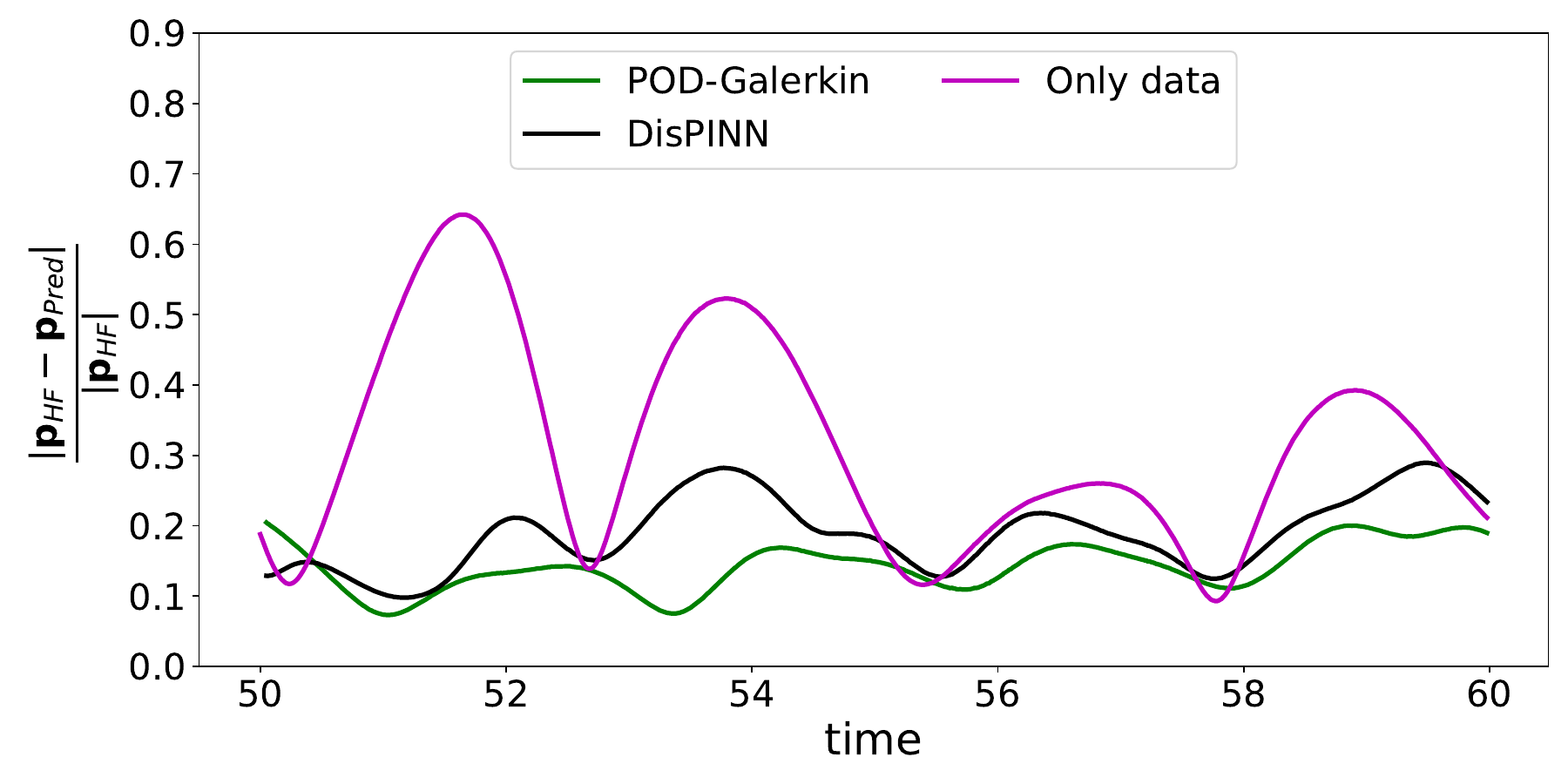}}
\end{subfigure}
\medskip
\caption{Comparison of the relative L2 errors (between the benchmark CFD results and predicted results from DisPINN) of the velocity field and pressure field considering three case scenarios: purely data-driven approach, POD-Galerkin projection based approach and current DisPINN approach}
\label{fig:Error_cylinder}
\end{figure}

At $\nu = 0.005$, we compare the reconstruction of the velocity field obtained from the current discretized PINN with benchmark CFD results in \autoref{fig:FPC_u}. The flow physics around the cylinder predicted by the proposed DisPINN closely matches the benchmark results shown in three time instances: $50.5$ s, $55$ s, and $58$ s. To demonstrate the relative L2 error (computed using \autoref{eq:L2error}) between the predicted results using the current discretized PINN and benchmark CFD results in \autoref{fig:L2_Error_005}, three cases are considered. The first case uses only the data-driven loss term in the neural network, without incorporating the physics-based loss term from the external CFD solver, resulting in a significantly high error at those time instants where no numerical data is available. The maximum relative error in this case is $12.83\%$. The second case, employing the POD-Galerkin-based ROM as discussed in \cite{STABILE2018273}, shows the best prediction accuracy among the three approaches, with a maximum error of $1\%$. The hybrid approach, using the DisPINN algorithm, results in a maximum relative error of $3.5\%$. We like to remark here that POD-Galerkin projection-based ROM, however, offers the best prediction accuracy, but the prediction time is not real-time, and this approach encompasses several numerical stability issues (as mentioned in \cite{STABILE2018273}). On the contrary, the prediction obtained from the DisPINN algorithm is real-time and numerically less challenging, and therefore offers a good balance between prediction accuracy, prediction time and numerical stability.      

\begin{figure}[ht]
           \centering

        % First Row
        \begin{subfigure}[b]{1\textwidth}
            \centering
            \subfloat[t = $50.5s$]{\label{fig:p_benchmark_50}\includegraphics[width=.25\linewidth]{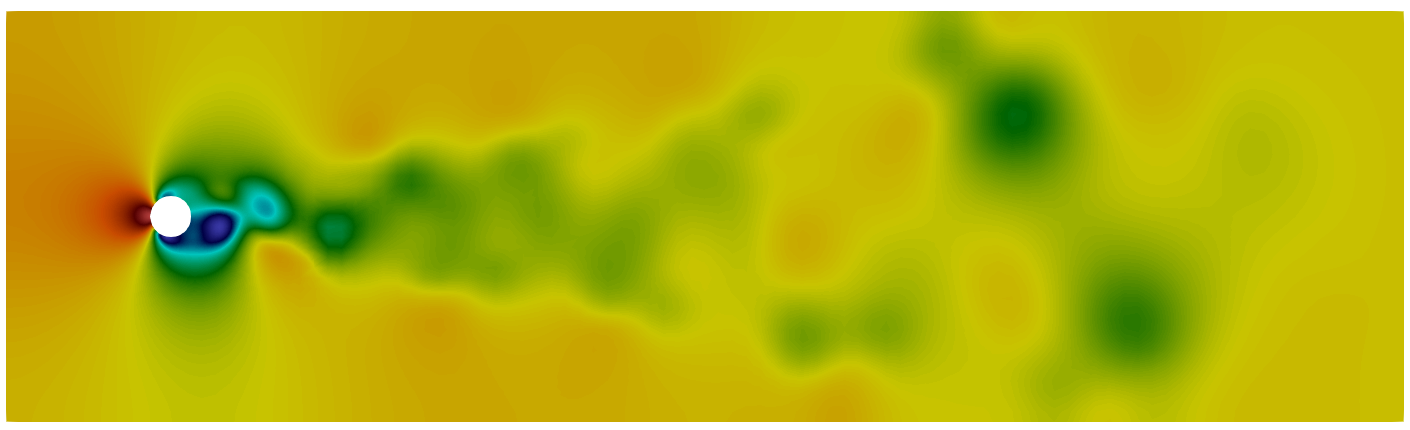}}
            \subfloat[t = $55s$]{\label{fig:p_benchmark_500}\includegraphics[width=.25\linewidth]{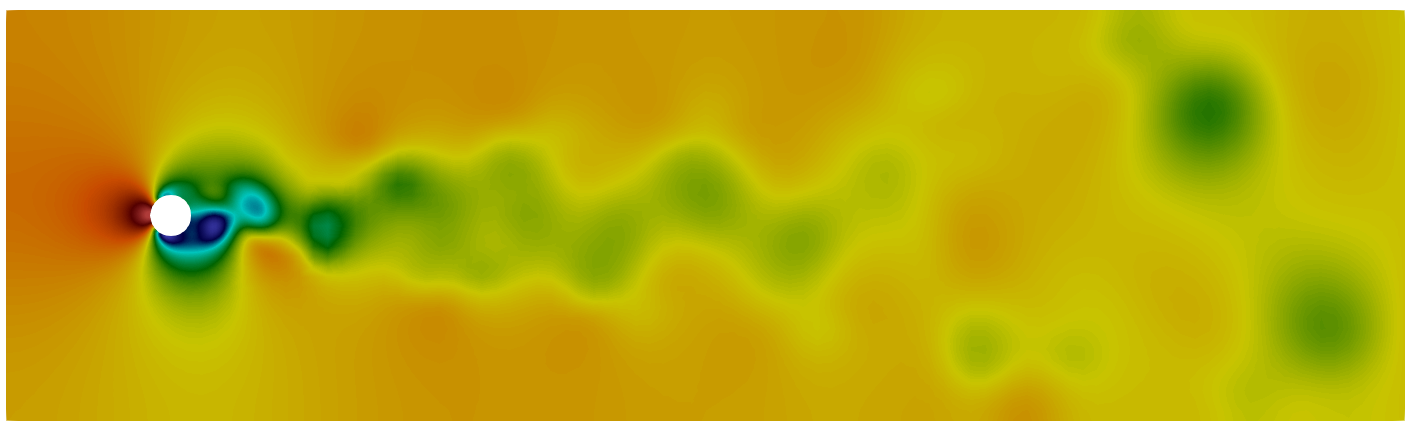}}
            \subfloat[t = $58s$]{\label{fig:p_benchmark_800}\includegraphics[width=.25\linewidth]{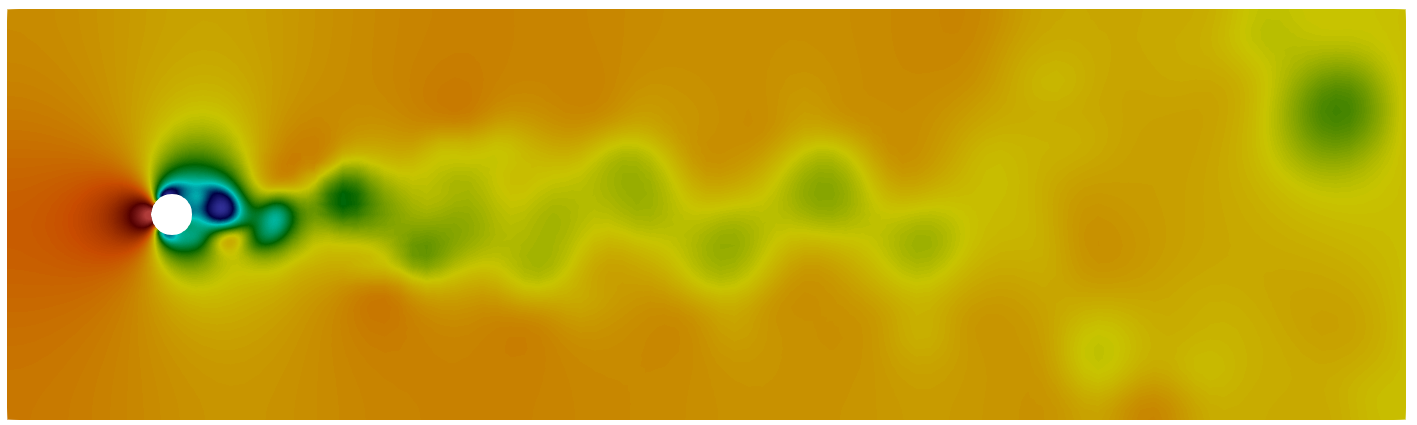}}
        \end{subfigure}

        \vspace{0.3cm} % Space between rows

        % Second Row
        \begin{subfigure}[b]{1\textwidth}
            \centering
            \subfloat[t = $50.5s$]{\label{fig:p_prediction_50}\includegraphics[width=.25\linewidth]{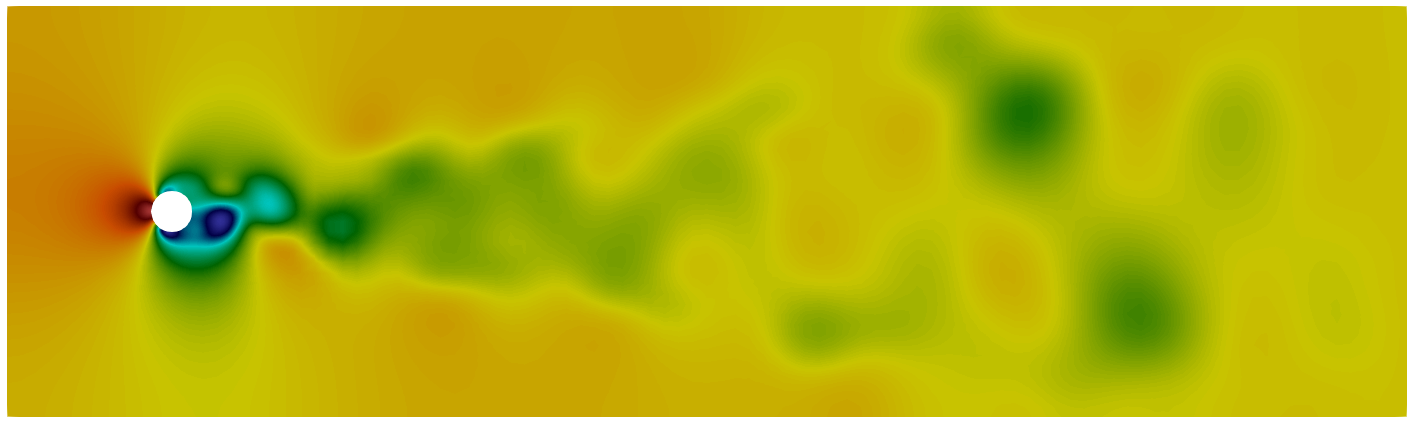}}
            \subfloat[t = $55s$]{\label{fig:p_prediction_500}\includegraphics[width=.25\linewidth]{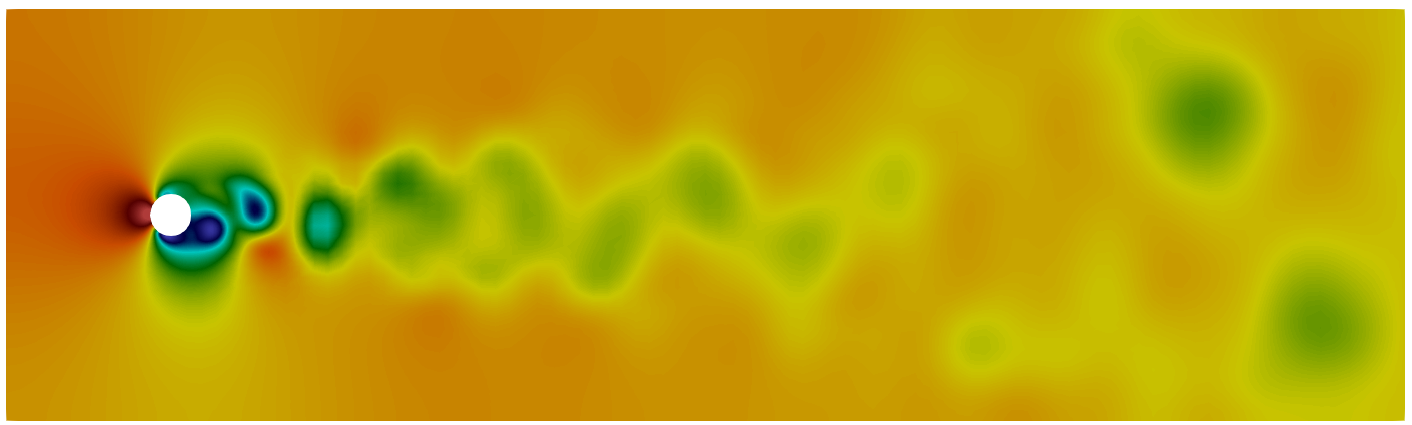}}
            \subfloat[t = $58s$]{\label{fig:p_prediction_800}\includegraphics[width=.25\linewidth]{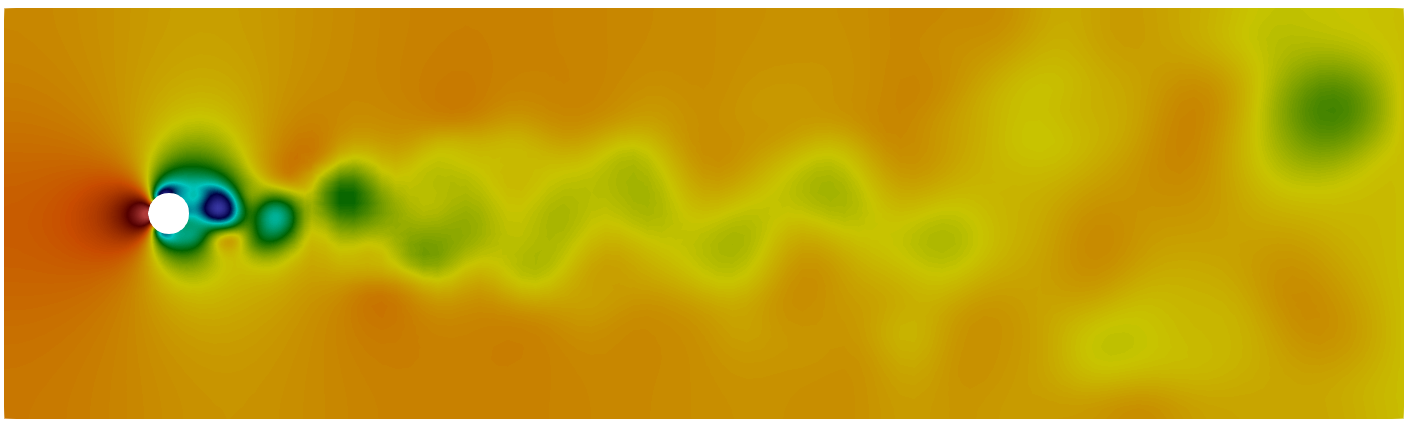}}
        \end{subfigure}

    \begin{subfigure}[b]{1\textwidth}
            \centering
{\label{fig:FPC_benchmark_50_3}\includegraphics[width=0.25\linewidth]{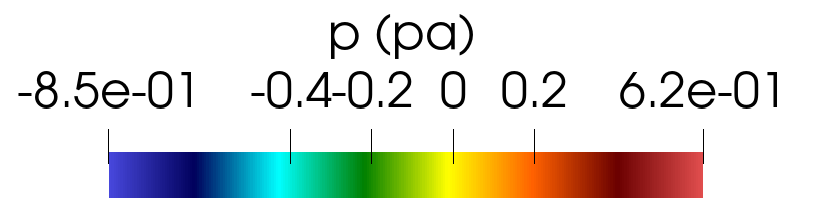}}          
        \end{subfigure}

    \medskip
    \caption{Comparison of the pressure field between (a) benchmark CFD results and (b) ANN-DisPINN prediction results of the reduced order incompressible viscous flow past a cylinder at $\nu=0.005$ and $t = $ $50.5$ $s$, $55$ $s$ and $58$ $s$}.
    \label{fig:FPC_p}
\end{figure}

For the reconstruction of the pressure field at $\nu = 0.005$ $\text{m}^2/\text{s}$, \autoref{fig:FPC_p} compares the predicted pressure field from the DisPINN with the benchmark results. Meanwhile, \autoref{fig:L2_Error_005_pressure} illustrates the variation of the relative error over time. When only the data-driven loss term is considered, the maximum error in the pressure field reconstruction is $40\%$. Using the POD-Galerkin approach, the error decreases to $18\%$, while the DisPINN approach, keeping a balance between the classical method and the neural network-based approach, results in a maximum error of $21\%$. 

We now aim to predict the velocity and pressure fields for a new parameter value, $\nu = 0.0055$ $\text{m}^2/\text{s}$, for which no spatio-temporal dataset corresponding to this physical parameter was included during the training of the neural network. \autoref{fig:L2_Error_0055} and \autoref{fig:L2_Error_0055_pressure} display the relative error associated with the velocity and pressure fields, respectively, using different data-driven and physics-based approaches. Three approaches are considered: First, using only the data-driven loss term, the prediction error is significantly high, with an $18\%$ error for the velocity field and a $64\%$ error for the pressure field. The DisPINN approach results in a maximum prediction error of $11\%$ for the velocity field and $28\%$ for the pressure field. These errors are higher compared to the reconstruction at $\nu = 0.005$. The POD-Galerkin approach also incurs significantly higher errors, with $7\%$ for the velocity field and $20\%$ for the pressure field. To investigate the reason for the higher errors in the flow physics predictions at this $\nu$ value, we recall that the time window considered was from $50$ s to $60$ s, which is not within the periodic regime for $\nu = 0.005$, whereas the flow physics at $\nu = 0.006$ has become periodic within this time window. This anomaly in the flow physics contributes to the higher errors in both the POD-Galerkin and DisPINN predictions. Additionally, the prediction error is expected to decrease if training datasets with a higher number of physical parameter values for $\nu$ are included.

\section{Concluding Remarks}\label{sec:Conclusion}

In this work, we have demonstrated several algorithms for coupling any external solver with a physics-informed neural network, in those cases where the PINN and CFD solver belong to separate environments. The prediction of the solution field of the nonlinear transport equation with very sparse numerical data is significantly good with the benchmark results. The proposed method demonstrates accurate prediction in the future time steps outside the training regime as well. The current work also suggests the relevance of reduced-order models in the context of PINN training, as it reduces the dimensionality of the system, which results in efficient optimization. The novel correction term introduced to the physics-based loss term ensures accurate gradient computation of the loss term with respect to the network hyperparameter. The behaviours of the MSE loss term decay with respect to the epoch number in the context of full order and reduced order systems, coupled with PINN, draw particular interest to the author. The future research directions of the authors will include variations of the deep learning algorithms and different optimization approaches in the training of DisPINN. 
\section*{\bf{Disclosure Statement}}\label{sec:DS}
The authors report no potential conflict of interest.

\section*{Acknowledgements}
The authors gratefully acknowledge the financial support under the scope of the COMET program within the K2 Center “Integrated Computational Material, Process and Product Engineering (IC-MPPE)” (Project No 886385). This program is supported by the Austrian Federal Ministries for Climate Action, Environment, Energy, Mobility, Innovation and Technology (BMK) and for Labour and Economy (BMAW), represented by the Austrian Research Promotion Agency (FFG), and the federal states of Styria, Upper Austria, and Tyrol. The authors also acknowledge financial support from MUR PRIN project $\text{NA\_FROM\_PDEs}$, INDAM GNCS and MUR PRIN project ROMEU. Giovanni Stabile acknowledges the financial support under the National Recovery and Resilience Plan (NRRP), Mission 4, Component 2, Investment 1.1, Call for tender No. 1409 published on 14.9.2022 by the Italian Ministry of University and Research (MUR), funded by the European Union – NextGenerationEU– Project Title ROMEU – CUP P2022FEZS3 - Grant Assignment Decree No. 1379 adopted on 01/09/2023 by the Italian Ministry of Ministry of University and Research (MUR) and acknowledges the financial support by the European Union (ERC, DANTE, GA-101115741). Views and opinions expressed are however those of the author(s) only and do not necessarily reflect those of the European Union or the European Research Council Executive Agency. Neither the European Union nor the granting authority can be held responsible for them.
\clearpage
\newpage
\bibliographystyle{plain}
\bibliography{mybib_clean}
\end{document}